\documentclass[11pt,reqno,a4paper]{amsart}

\usepackage[left=3.5cm,right=3.5cm,top=3cm,bottom=3cm]{geometry}
\usepackage{amsmath}
\usepackage{amsfonts,amssymb}
\usepackage{enumerate}

\usepackage{amsthm}

\usepackage{dsfont}

\usepackage{array}

\usepackage{color}

\usepackage[hyperindex, colorlinks=true]{hyperref}

\newcommand{\eqn}{\begin{eqnarray}}
\newcommand{\een}{\end{eqnarray}}

\def\eps{\epsilon}

\newtheorem{theorem}{Theorem}[section]

\newtheorem{lemma}{Lemma}[section]
\newtheorem{prop}[theorem]{Proposition}

\theoremstyle{definition}

\newtheorem{remark}{Remark}

\numberwithin{equation}{section}

\begin{document}

\date{}
\author{ Francisco Gancedo and  Omar Lazar  }
\title{Global well-posedness for the 3D Muskat problem in the critical Sobolev space } 
\maketitle

\begin{abstract}
We prove that the 3D stable Muskat problem is globally well-posed in the critical Sobolev space $\dot H^2 \cap \dot W^{1,\infty}$ provided that the semi-norm $\Vert f_0 \Vert_{\dot H^{2}}$ is small enough. Consequently, this allows the Lipschitz semi-norm to be arbitrarily large. The proof is based on a new formulation of the 3D  Muskat problem that allows to capture the hidden oscillatory nature of the problem. The latter formulation allows  to prove the $\dot H^{2}$ {\emph{a priori}} estimates. In the literature, all the known global existence results  for the 3D Muskat problem are for small slopes (less than 1). This is the first arbitrary large slope theorem for the 3D stable Muskat problem.\end{abstract}
    \bfseries
    \hypersetup{linkcolor=black}
  \tableofcontents

    \normalfont
\maketitle

\bibliographystyle{plain}

\section{{{Introduction}} }

In this article, we study the 3D Muskat problem which  models the dynamics of two incompressible and immiscible fluids with different densities and viscosities separated by a porous media (see \cite{Muskat}). This problem, initiated by Morris Muskat in the early '30, has appeared in the first time in the study of science of geophysics mainly for petroleum engineering applications (\cite{Naras}). His main contributions has been to introduce a mathematical concepts to the knowledge of flow of oil and gas in sands. Since then,  many other applications  such as in civil engineering or in modern biology have been studied  (see e.g. \cite{Kim}).  Since the fluids are immiscible and separated by a porous media, they therefore lie in two different time dependent domains. Set $\Omega_1(t)$ and $\Omega_{2}(t)$ these two different fluid regions. We assume that $\rho_{i}$ is the density of the fluid in the moving region $\Omega_i(t)$ and that the two fluids have the same viscosity (see {\it{e.g.}} \cite{GGJPS} for the viscosity jump case). The velocity $v_i$ in the fluid domain $\Omega_i(t)$ for $i=1,2$,  is given by the following  so-called Darcy's \cite{Darcy} law as follows
\begin{eqnarray} \label{eq1}
\frac{\mu}{\kappa}v_{i}&=&(0,0,g\rho_{i})-\nabla{P_{i}},\\
\nabla\cdot v_{i}&=&0. \label{eq2}
\end{eqnarray}
Where $g$ is the gravity, $\kappa$ is the permeability of the porous media, $\mu$ is the viscosity. Since $g$, $\kappa$ and $\mu$ are fixed constants, without loss for generality,  we may assume that there are all equal to 1 for simplicity. The second identity means that the two fluids are incompressible. Recall that  $P_i$ is the pressure on the different fluid domains, while on the interface $\partial \Omega_1 (t)=\partial \Omega_2 (t)$ the pressure are equal that is $P_1=P_2$. Lastly, since the density $\rho_i$ is transported by the flow, it obeys the following equation
\begin{eqnarray} \label{eq3}
\partial_t \rho_i + v_i\cdot\nabla \rho_i=0.
\end{eqnarray}

The coupling of equations \eqref{eq1},  \eqref{eq2},  \eqref{eq3} is the incompressible porous media equation  (\cite{Muskat}).  Note that all those physical quantities namely $v_i, \rho_i, P_i$ are functions of $(x,t) \in \mathbb{R}^3 \times [0,\infty)$. In particular, since the two fluids have different densities, $\rho_i$ is a step function, that is $$\rho(x,t)=\rho_{1} \mathds{1}_{\Omega_1(t)}(x)+\rho_{2} \mathds{1}_{\Omega_2(t)}(x).$$   This problem is analogous to the so-called Hele-Shaw equation \cite{HS1,HS2}. We refer to \cite{RT, GL, G, CCG} for a complete picture of this analogy and to \cite{AB,AMS}  for some recent mathematical developments on this equation and related models. \\
Since $\rho_1 \ne \rho_2$, we may assume that $\rho_1<\rho_2$. In that case, the word "stable" Muskat problem means that  $\Omega_2(t)$ corresponds to the heavier fluid domain which lies below $\Omega_1(t)$ which is the lighter fluid domain. This physical structure is preserved for any time as long as the interface is a graph of a regular enough function and this is the case as long as the Raleygh-Taylor condition is satisfied (see \cite{RT}). Indeed, a common assumption when studying the moving fluid domains is to parametrize the interface as being the graph of a sufficiently regular function. In this case the Rayleigh-Taylor simplifies to  $\rho_2-\rho_1>0$. By using classical tools from potential theory, it was shown in \cite{CG} that the interface obeys a nice contour equation which is both nonlocal (unlike its Eulerian version) and nonlinear.  This formulation gives a closed equation which is fully determined only by the dynamics of the interface itself. The dynamics of this moving interface is a function $f$ which depends of the position $x\in \mathbb R^{2}$ and time $t\geq0$,   This gives rise to an evolution equation which is called the Muskat problem. We shall further assume that we are dealing with an interface which is flat at infinity and that there is no surface tension.

 In this paper we shall focus on the 3D case.  The 3D Muskat problem reads as follows
\begin{equation}  \label{f1}
\ (\mathcal{M}_1) \ : \\\left\{
\aligned
& f_{t} (t,x) = \frac{\rho}{2\pi} P.V.\int \nabla_{x} \Delta_{y} f\cdot \frac{y}{\vert y \vert^2} \frac{1}{(1+\Delta^{2}_{y} f )^{3/2}}  \ dy, \ \ (x,t) \in \mathbb R^2 \times [0,T] 
\\
& f(0,x)=f_{0}(x), \nonumber
\endaligned
\right.
\end{equation}
where $\rho=\rho_2-\rho_1>0$ and  the operator $\Delta_{y}f(x,t)=\frac{f(x,t)-f(x-y,t)}{\vert y \vert}$. Note that the $p.v.$ is mainly needed  when $y$ approaches 0, some models have been studied taking into account this fact  (see  e.g. \cite{RGS}). Local existence for this equation in subcritical spaces either in $2D$ or $3D$ has been studied in several articles. Local existence in the Sobolev space in $H^{k}$, $k\geq3$ and illposedness results in the unstable regime have been shown in \cite{CG}. In \cite{CBS}, Chang, Granero-Bellinch\'on and Shkoller proved local well-posedness  in $H^{2}$ provided the norm $H^{3/2+\eps}$, $\eps \in (0,1/2)$  is small enough.  In \cite{CGSV}, Constantin, Gancedo, Schydkoy and Vicol were able to prove that the Muskat problem is locally-well posed in $\dot W^{2,p}, p>1$. They also proved a regularity criteria in terms of the uniform continuity of the bounded slope  (see also \cite{GL} where a very weak regularity criteria is proved). The later result has been recently extended in \cite{AM} to the 3D case and to the wider class of subcritical Sobolev spaces $W^{s,p}$ where $s\in (1+1/p,2)$ and $p\in(1,\infty).$  In \cite{Matioc2}, Matioc proved local-wellposedness in the subcritical Sobolev space $H^{3/2+\eps}$, $\eps \in (0,1/2)$. By using a purely paradifferential approach, Nguyen and Pausader  \cite{NP} were able to prove that the Muskat problem is locally-well posed in $H^{s}$, $s>1+d/2$ regardless of the characteristic of the fluids. In the 2D case, the homogeneous version of the result in \cite{NP}  has been obtained by Alazard and the second author \cite{AL} using a paralinearization formula of the Muskat equation \cite{AL}. The latter allows to identify the most important terms in the study of the Cauchy problem.  \\

Similarly, up to an integration by parts (see \cite{CG}), the 3D Muskat problem may be written as 
\begin{equation}\label{f2}
\ (\mathcal{M}_1) \ : \\\left\{
\aligned
& f_{t} (t,x) =  \frac{\rho}{2\pi} P.V.\int \frac{\nabla f(x)\cdot y-(f(x,t)-f(x-y,t))}{(1+\Delta^{2}_{y}f)^{3/2}}  \frac{dy}{\vert y \vert^3},
\\ \nonumber
& f(0,x)=f_{0}(x).
\endaligned
\right.
\end{equation}
The latter formulation is well adapted when dealing with the Cauchy problem for the Muskat equation with data in the Lipschitz class. Indeed, it has been used for instance in the recent work by Cameron \cite{Cam2}  to prove global regularity for small slopes for the 3D Muskat problem. Besides being a physically relevant quantity when dealing with the geometry of the moving interface, the Lipschitz semi-norm is also a fundamental quantity in the Muskat problem (see the survey \cite{GL, G}).  \\

Importantly, the Muskat equation has a scaling. Namely, if $f$ is a solution to 3D Muskat problem with initial data $f_0$ so does the whole family $\lambda^{-1}f(\lambda x, \lambda t)$ with initial data $\lambda^{-1}f_{0}(\lambda x)$, where $\lambda>0$. Recall that a space is called critical if its norm (or semi-norm) is left invariant by the scaling of the equation. In the case of the 3D Muskat problem, it is not difficult to observe that the Lipschitz space, the Wiener space studied in \cite{CCGS},  the homogeneous Sobolev space $\dot H^2$ or the homogeneous Besov space $\dot B^{1}_{\infty,\infty}$ are examples of critical spaces for the 3D Muskat problem. To get a first idea of the structure of the equation a classical idea consists in linearizing around the trivial solution. By doing so, one may check that the equation reduces to
$$
\partial_{t} f(x,t)=\frac{\rho}{2\pi} \Lambda f
$$  
where in 2D, $$\Lambda f(x,t)= \frac{P.V.}{\pi} \int \frac{f(x,t)-f(x-y,t)}{\vert y\vert^{3}} \ dy$$

This linearization shows that one needs $\rho>0$ in order to ensure existence of a local solution to the "half" heat equation.   \\

The Cauchy problem for equation $\mathcal{M}_1$ in the critical setting is delicate, even if one assumes smallness of the initial data. Indeed, the Muskat problem is not a fully parabolic PDE since regular enough solutions may blow-up as it has been show by Castro, C\'ordoba, Gancedo and Fefferman in \cite{CCFG1,CCFG2}. Indeed, they proved that there exists a class of smooth initial data which fails to be $C^4$ regular after sometimes and after a later time becomes a non-graph (see also \cite{GG}).  The instablity of the Cauchy problem associated to regular enough initial data is also very well described in a series of papers by  C\'ordoba, G\'omez-Serrano and Zlato\v{s} (\cite{shift1,shift2}). They were able to show some special dynamical scenarios are possible {\it{e.g.}} solutions passing from stable regime to unstable regime and finally go back to stable regime. Another kind of singularity are the so-called splash singularity (the curve self intersect in a point) or splat singularity (the curve self intersect in set a of Lebesgue measure $>0$) while its regularity is preserved. For the Muskat problem, both splash \cite{GS} and splat singularities \cite{CG2010,CP} have been ruled out. In the one phase Muskat problem problem  splash singularities are possible as it was shown by Castro, C\'ordoba, Fefferman,  Gancedo and L\'opez-Fern\'andez in \cite{CCCFM}). These kind of singularities have been shown to exist or ruled out for water waves and related fluid equations (see \cite{CCFGG, FIL, CS1, CS2, GS}). Note that the Muskat can be seen as the "parabolic" version of the water-waves equations (see e.g. \cite{GL}). \\

 All the singularity results known require initial data which are sufficiently regular and with sufficiently high slope. Global existence results for very small slopes have been obtained by Constantin and Pugh \cite{CP} or Escher-Matioc \cite{EM} they were able to ruled out turnover scenario. Actually, if one assumes that the initial data is sufficiently small in the critical Lipschitz space $\dot W^{1,\infty}$, then the Muskat problem turns out to be more stable. More precisely, there is a maximum principle for the slope (\cite{CG2009}) in the sense that, if the Lipschitz semi-norm is initially smaller than 1 so do the solutions for all time. In \cite{CCGRPS}, Constantin, C\'ordoba, Gancedo,  and Rodriguez-Piazza and Strain were able to prove that if the initial data is at least $H^{3}$ (to ensure local existence \cite{CG}) and if the initial  data is smaller than 1/3 in the Lipschitz class, then the 3D Muskat problem is globally well-posed. We refer also to \cite{PSt} where decay estimates are obtained. Recently, Cameron \cite{Cam2} was able to construct global unique solution for initial data $\Vert \nabla_{x} f_{0} \Vert_{L^{\infty}} <\frac{1}{\sqrt{5}}$. The unique solution can be unbounded provided that it grows sublinearly.  However, unlike his result in the 2D case (\cite{Cam1}),  the main results in the 3D case deals with small slopes only. \\
 
  While arbitrary large slope results have been shown to exist globally for the 2D Muskat problem in : \\
  
 \noindent - Deng, Lei, and  Lin \cite{Lin} (under a monotonicity assumption) \\
- Cameron \cite{Cam1} (under the condition that  $\sup f'_0(x) \times  \sup -f'_0(y)<1$) \\
   - C\'ordoba and the second author \cite{CL} (small data in the critical $\dot H^{3/2}$ space),\\
   
   no large solutions in Lipschitz are known to exist for the 3D Muskat problem. In terms of the geometry of the interface, the condition of very small slopes $(<1)$ is  quite restrictive. \\
 
 The aim of this article is to show that the 3D Muskat problem is globally well-posed for any large initial data in Lipschitz. Indeed, we shall only assume smallness in the critical $\dot H^{2}$ semi-norm. So the slope can be arbitrarily large, this is the first result of large slope solutions for the 3D Muskat problem. \\
 
 Besides being mathematically challenging to prove global results without any smallness assumption on the Lipschitz semi-norm, it is also physically relevant since it would show that the interface can  be highly oscillating in an arbitrarily short time. This is obviously impossible to observe if the slope is small. Also, allowing the slope to be arbitrarily large shows that there exist solutions which can be arbitrarily close to the turnover phenomena observed by Castro, C\'ordoba, Fefferman,  Gancedo and L\'opez-Fern\'andez in \cite{CCCFM}) but without never reaching it.  \\

When dealing with the Cauchy problem for data in the critical $\dot H^2$ space, both aforementioned formulations  give rise to severe difficulties to close the {\it{a priori}} estimates for the most singular terms. This motivate the introduction of a new formulation to treat the Cauchy problem \eqref{f1} for initial data in $\dot H^{2}$. The idea behind this new formulation in terms of oscillatory integral was pioneered in an article by C\'ordoba and the second author \cite{CL} were they studied the Cauchy problem for 2D Muskat equation with regular enough data and small $\dot H^{3/2}$ semi-norm.  However, the 3D case (2D interface) is not only more nonlinear than the 2D case (1D interface) but also more technical because of the fact that one has to deal with directional derivatives.  The fact that the rational function in 
$\Delta_{y}f$ appearing in the Muskat equation cannot be seen as the restriction of the Fourier transform of some well chosen $L^1$ function (in the same spirit as \cite{CL}) generates some technical difficulties. Also, one of the most difficult term is $S_{2,2}$. This is mainly because one looses the nice symmetry of the Hilbert transform which gives rise to nice controlled commutators in the case of the 2D Muskat problem \cite{CL}. In higher dimension, we get the Riesz transforms but due to the fact that the critical space becomes $\dot H^2$ it seems not possible to get some nice cancellations and symmetries.  One has to guess which decompositions will give the desired control to close the energy estimates.

\section{{{Main result}} }

\begin{theorem}

Let $\mathcal{F}(x)=C({1+x^{2}})^{-3/2}$, where $C>0$ is a fixed constant.  For any initial data $f_{0} \in \dot H^{2} \cap \dot W^{1,\infty}$ with $\Vert f_{0} \Vert_{\dot H^{2}}< \mathcal{F}(\Vert f_{0}\Vert_{\dot W^{1,\infty}})$  small enough, then, there exists a unique global solution $f$ to the 3D Muskat problem such that $f \in L^{\infty}([0,T], \dot H^2 \cap \dot W^{1,\infty} ) \cap L^2([0,T];\dot H^{5/2})$ for all $T>0$.
\end{theorem}
\begin{remark} This theorem allows the interface to be arbitrarily large in $\dot W^{1,\infty}$ which is the first result of this kind in the 3D case. Note that the smallness is only assumed on the critical $\dot H^2$ Sobolev semi-norm.  Besides, this theorem is fully dealing with the critical setting in the sense that both the initial data and the smallness lie in critical spaces.
\end{remark}
\begin{remark} The proof of the {\it{a priori}} estimates in $\dot H^2$  is based on a series of decomposition of the terms together with estimates on homogeneous Besov spaces. This $\dot H^2$ control shows that there is a regularizing effect of order $L^{2}\dot H^{5/2}$. The control of the slope by means of the $L^{2}\dot H^{5/2}$ semi-norm is obtained thanks to a combination of the study of the evolution of the extrema (justified thanks to Rademacher's theorem)  together with Besov estimates. 
\end{remark}
\begin{remark} When performing $\dot H^{2}$ {\it{a priori}} estimates, the dissipation one hope for is  of fractional order. This amount to take fractional derivatives into the nonlinear term.  One would need to use multilinear estimates of singular integral operators together with estimates of composition functions (\cite{BM}). This may lead to tedious computations. However, the strategy to get the 
$\dot H^{2}$ {\it{a priori}} estimates presented in this paper avoid this difficulty.  
\end{remark}

The plan of the paper is the following, in the next section we shall introduce a new formulation of the 3D Muskat problem in terms of oscillatory integrals. In the second section, we shall give the definition of the functional spaces together with notations of some operators that will be used throughout the article. The third section, which is the central part of the article, is devoted to the proof of the $\dot H^{2}$ {\emph{a priori}} estimates. The fourth section contains the Sobolev energy inequality. The fifth section is the control of the slope together with a boostrap argument to close the estimates with respect to critical quantities only. {  {The sixth and last section is the proof the uniqueness.}}

\section{{{A new formulation of the 3D Muskat problem}}}

Let us recall that the  Muskat  equation in $\mathbb R^{3}$ in the stable case and when the interface is parametrized as a graph is given by the following 2D equation
\begin{equation} 
\ (\mathcal{M}_1) \ : \\\left\{
\aligned
& f_{t} (t,x) = \frac{\rho}{2\pi} P.V.\int \nabla_{x} \Delta_{y} f. \frac{y}{\vert y \vert^2} \frac{1}{(1+\Delta^{2}_{y} f )^{3/2}}  \ dy \hspace{2cm} 
\\ \nonumber
& f(0,x)=f_{0}(x).
\endaligned
\right.
\end{equation}
In this section we shall prove the following Proposition which gives an equivalent formulation of the 3D Muskat  in terms of oscillatory integrals.

\begin{prop} \label{formulation}
Consider the following Cauchy problem
\begin{equation} 
\ (\mathcal{M}_2) \ : \\\left\{
\aligned
& f_{t} (t,x) = \frac{\rho}{2\pi} P.V.\int \nabla_{x} \Delta_{y} f. \frac{y}{\vert y \vert^2} \cos(\arctan(\Delta_{y} f )) \int_{0}^{\infty} e^{-k} \cos(k\Delta_{y} f ) \ dk  \ dy, \hspace{2cm} 
\\ \nonumber
& f(0,x)=f_{0}(x).
\endaligned
\right.
\end{equation}
Then,
\begin{center}
$(\mathcal{M}_1) \Longleftrightarrow (\mathcal{M}_2).$
\end{center}
\end{prop}

\noindent {\bf{Proof of Proposition \ref{formulation}}}  {  {One may easily check that, for any $x \in \mathbb{R}$ 
$$
\frac{1}{1+x^2}= \int_{0}^{\infty} e^{k} \cos(kx) \ dk \ \ {\text{and}} \ \ \cos(\arctan(x))=\frac{1}{\sqrt{1+x^2}}
$$
Hence, for any $x \in \mathbb{R}$
$$
\frac{1}{(1+x^2)^{3/2}}=\int_{0}^{\infty} e^{k} \cos(k x) \cos(\arctan( x)) \ dk.$$
In particular, for $x=\Delta_{y} f$, one gets the identity
$$
\frac{1}{(1+\Delta^{2}_y f)^{3/2}}=\int_{0}^{\infty} e^{k} \cos(k \Delta_{y}f) \cos(\arctan( \Delta_{y}f)) \ dk.
$$
Hence $\mathcal{M}_{1} \Leftrightarrow \mathcal{M}_{2}.$
}}
\qed
\section{{{Functional setting and notations}}}

As usually, for $s>0,$ $\dot H^{s}$ denotes the homogeneous Sobolev space endowed with the semi-norm
$$
\Vert f \Vert_{\dot H^{s}} = \Vert \Lambda^{s} f \Vert_{L^{2}}
$$

The definition of the homogeneous Besov spaces that we shall use have been introduced by Oleg Vladimirovich Besov in \cite{Besov}. Let $(p,q,s) \in [1,\infty]^2 \times \mathbb R^2$,  a tempered distribution $f$ (we assume that its Fourier transform is locally integrable near 0)   belongs to the homogeneous Besov space $\dot B^{s}_{p,q}(\mathbb R^2)$ if and only if the following  semi-norm is finite
$$
\Vert f \Vert_{\dot B^{s}_{p,q}} =   \left\Vert \frac{\Vert \mathds{1}_{ ]0,1[}(s) \delta_{y}f + \mathds{1}_{[1,2[}(s) (\delta_{y}f+\bar \delta_{y} f)\Vert_{L^{p}}}{ \vert{y}\vert^{s} } \right\Vert_{L^{q}(\mathbb R^2, \vert y \vert^{-2} dy)} <\infty,
$$
where $\delta_{y} f(x)=f(x)-f(x-y)$ and $\bar \delta_{y} f(x)=f(x)-f(x+y)$. \\

We have the following embedding between homogeneous Besov spaces   (see e.g. \cite{BCD}, \cite{PGLR}, \cite{Peetre}). For all $(p_1,p_2,r) \in [1,\infty]^3$ such that $p_1 \leq p_2$ and $q_1 \leq q_2$ we have
$$
\dot B^{s_1}_{p_1,r}(\mathbb R^2) \hookrightarrow \dot B^{s_2}_{p_2,r}(\mathbb R^2),
$$
where $(s_1,s_2) \in \mathbb R^2$ are so that  $s_1+\frac{2}{p_2}=s_2+\frac{2}{p_1}$. We also have for all $(p_1, s_1) \in [2,\infty] \times \mathbb R$,
$$
\dot B^{s_1}_{p_1,r_1}(\mathbb R^2) \hookrightarrow \dot B^{s_1}_{p_1,r_2}(\mathbb R^2),
$$
for all $(r_1,r_2) \in \ ]1,\infty]^2$ such that $r_1 \leq r_2$. \\


Throughout the article, we shall use the following notations:
\begin{itemize}
\item $\Delta_{y}f(x,t)=\frac{f(x,t)-f(x-y,t)}{\vert y \vert}$
\item $\delta_{y}f(x,t)=f(x,t)-f(x-y,t)$\\
\item $\bar\Delta_{y}f(x,t)=\frac{f(x,t)-f(x+y,t)}{\vert y \vert}$
\item $\bar \delta_{y} f(x,t)=f(x,t)-f(x+y,t)$\\
\item $S_{y}f(x,t)=\frac{2f(x,t)-f(x-y,t)-f(x+y,t)}{\vert y \vert}$
\item $s_{y}f(x,t)={2f(x,t)-f(x-y,t)-f(x+y,t)}$\\
\item $D_{y}f(x,t)=\frac{f(x+y,t)-f(x-y,t)}{\vert y \vert}$
\item $d_{y}f(x,t)=f(x+y,t)-f(x-y,t)$\\
\end{itemize}

For the sake of readability, we shall not write the time dependence. \\

The notation $\nabla_{i}$ will denote the gradient vector with respect to the variable $i\in \mathbb R^2$. The operator $\Delta$ will always mean the classical Laplacian with respect to $x$. \\

As well, $A \lesssim B$ means that there exists a fixed constant $C>0$ such that $A \leq C B$.

\section{{{{\it{A priori}} estimates in $\dot H^2$}}}

We shall use an energy method. That is, we shall do $\dot H^{2}$ {\it{a priori}} estimates which allows us to get enough compactness to pass to the limit in a regularized equation. Without loss of generality we may assume that $\rho=2\pi$. {  { To prove the existence of solution we use the following regularized Muskat equation which was introduced by Alazard and Hung (see \cite{AN2}). Let $\phi>0$ be a smooth bounded even function whose integral over $\mathbb R^2$ is 1 and such that $\phi(x)=1$ in $B_1$ (the ball of radius 1 centered at the origin) and 0 outside $B_{2}$. Let $\eps \in (0,1]$ and set $\phi_{\eps}(x):=\eps^{-1}\phi(\eps^{-1} \vert x \vert ).$ 
\begin{equation} \label{musk}
\ (\mathcal{M}_\eps):  \\\left\{
\aligned
&  \partial_t f_\eps (t,x) = \int\nabla_{x} \Delta_{y} f_\eps . \frac{y}{\vert y \vert^2} \cos(\arctan(\Delta_{y} f_\eps )) \\
 & \ \ \ \ \  \times \int_{0}^{\infty}  e^{-k} \cos(k\Delta_{y} f_\eps ) (1- \phi_{\eps}(y)) \ dk  \ dy \\
& \hspace{1,8cm} + \vert \log(\eps ) \vert^{-1} \Delta f_\eps
\\ \nonumber
& f_{\eps}(0,x)=f_{0}(x)*\phi_{\epsilon}(x). 
\endaligned
\right.
\end{equation}

Then, using section 2.6 in \cite{AN2} we know that for all $\eps \in (0,1]$ and all data in $H^{2}(\mathbb R^2)$  the regularized Muskat equation admits a unique global solution $f_{\eps} \in C^{1}([0,\infty), H^{\infty}(\mathbb R^2))$.}} The aim will be to prove that the associated solution to the Cauchy problem $(\mathcal{M}_\eps)$ will converge (as $\eps$ goes to 0) in some  Banach spaces (assuming  that the solution is further $L^2$ for the sake of simplicity). The strong compactness in $(L^2L^2)_{loc}$ will be obtained in the usual way thanks to the Rellich compactness theorem (see {\it{e.g.}} \cite{PGLR2}). To avoid redundancy, the details will be omitted since the arguments are classical. One may prove uniqueness by  using the same technics to estimate the difference of two solutions and we shall omit the details. In the sequel we assume that the solution  is from this regularized equation but we will omit to write the parameter $\eps$. \\

  The main effort will be devoted to the proof of the {\it{a priori}} estimates in the critical space $\dot H^{2}(\mathbb R^2)$. By taking the Laplacian of the Muskat equation  and multiplying by $\Delta f$ and finally integrating in the space variable, one finds
\begin{eqnarray*}
\frac{1}{2}\partial_{t} \Vert f \Vert^{2}_{\dot H^{2}}=\int \Delta f  \Delta \left(\int \nabla_{x} \Delta_{y} f. \frac{y}{\vert y \vert^2} \cos(\arctan(\Delta_{y} f )) \int_{0}^{\infty} e^{-k} \cos(k\Delta_{y} f )   dk  \ dy\right)   dx 
\end{eqnarray*}
Then, by using classical formulas for the differential operator $\Delta$, we find
\begin{eqnarray*}
\frac{1}{2}\partial_{t} \Vert f \Vert^{2}_{\dot H^{2}}&=& \int \Delta f \ \int \nabla_{x} \Delta_{y} \Delta f. \frac{y}{\vert y \vert^2} \cos(\arctan(\Delta_{y} f )) \int_{0}^{\infty} e^{-k} \cos(k\Delta_{y} f )  \ dk  \ dy  \ dx \\
&+& 2 \int \Delta f \ \int \Delta_{y} \Delta f \frac{y}{\vert y \vert^2}. \nabla_{x} \left( \cos(\arctan(\Delta_{y} f )) \int_{0}^{\infty} e^{-k} \cos(k\Delta_{y} f )  \ dk \right) \ dy   \ dx \\
&+& \int \Delta f \ \int \Delta_{y} \nabla_x f .\frac{y}{\vert y \vert^2} \Delta \left( \cos(\arctan(\Delta_{y} f )) \int_{0}^{\infty} e^{-k} \cos(k\Delta_{y} f )  \ dk \right) \ dy   \ dx \\
\end{eqnarray*}
hence, we obtain
\begin{eqnarray*}
\frac{1}{2}\partial_{t} \Vert f \Vert^{2}_{\dot H^{2}}&=&\underbrace{\int \Delta f \ \int \nabla_{x} \Delta_{y} \Delta f. \frac{y}{\vert y \vert^2} \cos(\arctan(\Delta_{y} f )) \int_{0}^{\infty} e^{-k} \cos(k\Delta_{y} f )  \ dk  \ dy  \ dx}_{:=\, \mathcal{S} (most \ singular \ term)} \\ 
&+& 2 \int \Delta f \ \int \Delta_{y} \Delta f \frac{y}{\vert y \vert^2}. \nabla_{x} \left( \cos(\arctan(\Delta_{y} f )) \right) \int_{0}^{\infty} e^{-k} \cos(k\Delta_{y} f )  \ dk  \ dy   \ dx \\
&+& 2 \int \Delta f \ \int \Delta_{y} \Delta f \frac{y}{\vert y \vert^2}. \cos(\arctan(\Delta_{y} f )) \nabla_{x} \left(  \int_{0}^{\infty} e^{-k} \cos(k\Delta_{y} f ) \ dk \right)  \ dy   \ dx \\
&+&\int \Delta f  \int \Delta_{y} \nabla_x f .\frac{y}{\vert y \vert^2} \Delta \left( \cos(\arctan(\Delta_{y} f )) \right) \int_{0}^{\infty} e^{-k} \cos(k\Delta_{y} f )  \ dk \ dy   \ dx \\
&+&2 \int \Delta f  \int \Delta_{y} \nabla_x f .\frac{y}{\vert y \vert^2} \nabla_{x} \left( \cos(\arctan(\Delta_{y} f )) \right). \nabla_{x}\left( \int_{0}^{\infty} e^{-k} \cos(k\Delta_{y} f ) \ dk \right) \ dy \ dx \\
&+& \int \Delta f \int \Delta_{y} \nabla_x f .\frac{y}{\vert y \vert^2} \cos(\arctan(\Delta_{y} f ))  \Delta \left( \int_{0}^{\infty} e^{-k} \cos(k\Delta_{y} f ) \right) \ dk \ dy   \ dx \\
&:=& \mathcal{S}+ \sum_{i}^{5} \mathcal{T}_i
\end{eqnarray*}
Our aim will be to control $\frac{1}{2}\partial_{t} \Vert f \Vert^{2}_{\dot H^{2}}$, we shall actually prove that $\frac{1}{2}\partial_{t} \Vert f \Vert^{2}_{\dot H^{2}}<0$ if the $\dot H^{2}$ is sufficiently small and the Lipschitz semi-norm does not blow-up. This, combining with the control of the Lipschitz semi-norm will give the main result by using a bootstrap argument.

\section{{{Estimates of the less singular term:}}  $\mathcal{T}=\displaystyle\sum_{i=1}^{5} \mathcal{T}_i$} \label{ti}

 {  {To estimate the less singular terms, one does not have to introduce any symmetrizations since the spatial derivatives will be well balanced. Indeed, these terms come from the differentiation of the oscillatory parts.  More precisely, we will prove the following Lemma.
\begin{lemma} The less singular terms can be estimated as follows
\begin{eqnarray} \label{flop}
\sum_{i=1}^{5}\mathcal{T}_{i}&\lesssim& \Vert f \Vert^{2}_{\dot H^{5/2}}\left(\Vert f  \Vert_{\dot H^{2}} +\Vert f  \Vert^2_{\dot H^{2}} \right)
\end{eqnarray}
\end{lemma}

\noindent{\bf{Proof of Lemma 6.1 }} The estimates of this terms do not require to use technical decompositions since it would be easy to balance the regularity in $x$ and in $y$. We shall estimate each $\mathcal{T}_{i}$ for $i=1,...,5$ separately.

}}

\subsection{Estimate of $\mathcal{T}_1$}

We start by estimating $\mathcal{T}_1 $, that is
\begin{eqnarray*}
\mathcal{T}_1 &=& 2 \int \Delta f \ \int \Delta_{y} \Delta f \frac{y}{\vert y \vert^2}. \nabla_{x} \left( \cos(\arctan(\Delta_{y} f )) \right) \int_{0}^{\infty} e^{-k} \cos(k\Delta_{y} f )  \ dk  \ dy   \ dx \\
&\lesssim&2 \int \vert \Delta f \vert  \ \int \vert \Delta_{y} \Delta f \vert \frac{1}{\vert y \vert} \left \vert \nabla_{x} \left( \cos(\arctan(\Delta_{y} f )) \right)  \right \vert  \ dy   \ dx. \\
\end{eqnarray*}
Then, since an easy computation gives 
$\left \vert \nabla_{x} \left( \cos(\arctan(\Delta_{y} f )) \right)  \right \vert \lesssim \vert \nabla_{x} \Delta_{y} f\vert$ one finds

\begin{eqnarray} \label{class}
\mathcal{T}_1  &\lesssim& \int \vert \Delta f \vert  \ \int \vert \Delta_{y} \Delta f \vert \frac{1}{\vert y \vert} \left \vert \nabla_{x} \Delta_{y} f   \right \vert  \ dy   \ dx \\ \nonumber
&\lesssim& \Vert f \Vert_{\dot H^{2}} \int \frac{\Vert \Delta \delta_{y} f \Vert_{L^{2}}}{\vert y \vert^{3/2}} \frac{\Vert \nabla_{x}\delta_{y}f \Vert_{L^{\infty}}}{\vert y \vert^{3/2}} \ dy \\ \nonumber
&\lesssim& \Vert f \Vert_{\dot H^{2}} \Vert \Delta f \Vert_{\dot B^{1/2}_{2,2}} \Vert \nabla_{x} f \Vert_{\dot B^{1/2}_{\infty,2}}  \nonumber \\
&\lesssim& \Vert f \Vert^{2}_{\dot H^{5/2}} \Vert f \Vert_{\dot H^{2}} \nonumber
\end{eqnarray}

\subsection{Estimate of $\mathcal{T}_2$}

Recall that 

\begin{eqnarray*}
\mathcal{T}_2&=& 2 \int \Delta f \ \int \Delta_{y} \Delta f \frac{y}{\vert y \vert^2}. \cos(\arctan(\Delta_{y} f )) \nabla_{x} \left(  \int_{0}^{\infty} e^{-k} \cos(k\Delta_{y} f )\ dk   \right)  \ dy   \ dx \\
&\lesssim&\int \vert \Delta f \vert \ \int \left \vert \Delta_{y} \Delta f \right\vert  \frac{1}{\vert y \vert^2} \left \vert \nabla_{x}  \cos(k\Delta_{y} f ) \ dk\right \vert  \ dy   \ dx 
\end{eqnarray*}
Using that $\left \vert \nabla_{x}  \cos(k\Delta_{y} f ) \right \vert \lesssim \vert \nabla_{x}\Delta_{y} f \vert$, one finds

\begin{eqnarray*}
\mathcal{T}_2&\lesssim&\int \vert \Delta f \vert  \ \int \vert \Delta_{y} \Delta f \vert \frac{1}{\vert y \vert} \left \vert \nabla_{x} \Delta_{y} f   \right \vert  \ dy   \ dx  
\end{eqnarray*}
which is the same estimate as \eqref{class}, so we conclude as in the estimate of $\mathcal{T}_2$ that is
\begin{eqnarray*}
\mathcal{T}_2\lesssim \Vert f \Vert^{2}_{\dot H^{5/2}} \Vert f \Vert_{\dot H^{2}}
\end{eqnarray*}

\subsection{Estimate of $\mathcal{T}_{3}$}

We have
\begin{eqnarray*}
\mathcal{T}_{3}&=& \int \Delta f \ \int \Delta_{y} \nabla_x f .\frac{y}{\vert y \vert^2} \Delta \left( \cos(\arctan(\Delta_{y} f )) \right) \int_{0}^{\infty} e^{-k} \cos(k\Delta_{y} f )  \ dk \ dy   \ dx 
\end{eqnarray*}

So that,
\begin{eqnarray*}
\mathcal{T}_{3}&\lesssim& \int\vert \Delta f \vert \ \int  \frac{\vert \Delta_{y} \nabla_x f \vert}{\vert y \vert}  \left \vert \Delta \left( \cos(\arctan(\Delta_{y} f )) \right) \right\vert  \ dy   \ dx \\
\end{eqnarray*}

Then, an easy estimate on $\Delta \left( \cos(\arctan(\Delta_{y} f )) \right)$ gives

\begin{eqnarray*}
\mathcal{T}_{3}\lesssim \int\vert \Delta f \vert \ \int\frac{\vert\Delta_{y} \nabla_x f \vert}{\vert y \vert}  \vert \Delta_{y} \Delta f \vert \ dy   \ dx + \int\vert \Delta f \vert \ \int \frac{\vert \Delta_{y} \nabla_x f \vert}{\vert y \vert}  \vert \Delta_{y} \nabla_{x} f \vert^{2} \ dy   \ dx  \\
\end{eqnarray*}
Using the same step as \eqref{class} one may estimate the first term in the right hand side as $\mathcal{T}_1 $. For the second term, we observe that

\begin{eqnarray} \label{class2}
\int\vert \Delta f \vert \ \int \frac{ \vert \Delta_{y} \nabla_x f \vert}{\vert y \vert}  \vert \Delta_{y} \nabla_{x} f \vert^{2} \ dy   \ dx  &\lesssim& \Vert f \Vert_{\dot H^{2}} \int \frac{\Vert \nabla_x \delta_{y}f \Vert^{3}_{L^{6}}}{\vert y \vert^{4}}  \ dy \nonumber \\
&\lesssim& \Vert f \Vert_{\dot H^{2}} \Vert  f \Vert^3_{\dot B^{5/3}_{6,3}} 
\end{eqnarray}
then, using that $\dot H^{7/3} \hookrightarrow \dot B^{5/3}_{6,3}$ we find
\begin{eqnarray*}
\int\vert \Delta f \vert \ \int \vert \frac{\Delta_{y} \nabla_x f}{\vert y \vert}  \vert \Delta_{y} \nabla_{x} f \vert^{2} \ dy   \ dx \lesssim \Vert f \Vert_{\dot H^{2}} \Vert  f \Vert^3_{\dot H^{7/3}}
\end{eqnarray*}
and finally, using that $\dot H^{7/3}=\left[\dot H^{2},\dot H^{5/2} \right]_{\frac{1}{3},\frac{2}{3}} $ we finally find that
\begin{eqnarray*}
\int\vert \Delta f \vert \ \int  \vert\frac{\Delta_{y} \nabla_x f}{\vert y \vert}  \vert \Delta_{y} \nabla_{x} f \vert^{2} \ dy   \ dx  &\lesssim& \Vert f \Vert^{2}_{\dot H^{5/2}} \Vert f \Vert^2_{\dot H^{2}}
\end{eqnarray*}
Hence,
\begin{eqnarray*}
\mathcal{T}_{3} \lesssim\Vert f \Vert^{2}_{\dot H^{5/2}} \Vert f \Vert^2_{\dot H^{2}}
\end{eqnarray*}

\subsection{Estimate of $\mathcal{T}_{4}$}

We have
\begin{eqnarray*}
\mathcal{T}_{4}=2 \int \Delta f \ \int \Delta_{y} \nabla_x f .\frac{y}{\vert y \vert^2} \nabla_{x} \left( \cos(\arctan(\Delta_{y} f )) \right). \nabla_{x}\left( \int_{0}^{\infty} e^{-k} \cos(k\Delta_{y} f ) \ dk \right) \ dy   \ dx 
\end{eqnarray*}
Therefore,
\begin{eqnarray*}
\mathcal{T}_{4}\lesssim  \int \vert \Delta f \vert \ \int \vert \Delta_{y} \nabla_x f  \vert \frac{1}{\vert y \vert} \left\vert \nabla_{x} \left( \cos(\arctan(\Delta_{y} f )) \right) \right\vert \left\vert \nabla_{x}\left( \int_{0}^{\infty} e^{-k} \cos(k\Delta_{y} f ) \ dk \right) \right\vert  \ dy   \ dx 
\end{eqnarray*}
Using that $$\left\vert \nabla_{x} \left( \cos(\arctan(\Delta_{y} f )) \right) \right\vert \lesssim \vert \Delta_{y} \nabla_{x} f \vert,$$ and that $$\left\vert \nabla_{x}\left( \int_{0}^{\infty} e^{-k} \cos(k\Delta_{y} f ) \ dk \right) \right\vert \lesssim \vert \Delta_{y} \nabla_{x} f \vert,$$
one finds,
\begin{eqnarray*}
\mathcal{T}_{4}\lesssim  \Vert f \Vert_{\dot H^{2}} \int \frac{\Vert \nabla_x \delta_{y}f \Vert^{3}_{L^{6}}}{\vert y \vert^{4}}  \ dy \lesssim \Vert f \Vert_{\dot H^{2}} \Vert  f \Vert^3_{\dot B^{5/3}_{6,3}}. 
\end{eqnarray*}

This is the same estimate as \eqref{class2}, hence following exactly the same step as the control of $\mathcal{T}_{2}$ we finally find that 
\begin{eqnarray*}
\mathcal{T}_{4} \lesssim\Vert f \Vert^{2}_{\dot H^{5/2}} \Vert f \Vert^2_{\dot H^{2}} 
\end{eqnarray*}

\subsection{Estimate of $\mathcal{T}_{5}$}

We write
\begin{eqnarray*}
\mathcal{T}_{5}&=&\int \Delta f \ \int \Delta_{y} \nabla_x f .\frac{y}{\vert y \vert^2} \cos(\arctan(\Delta_{y} f )) \Delta \left( \int_{0}^{\infty} e^{-k} \cos(k\Delta_{y} f ) \ dk \right) \ dy   \ dx  \\
&\lesssim&\int \vert \Delta f  \vert \ \int \frac{\vert \Delta_{y} \nabla_x f \vert}{\vert y \vert} \left\vert \Delta \left( \int_{0}^{\infty} e^{-k} \cos(k\Delta_{y} f ) \ dk \right) \right \vert \ dy   \ dx \\
\end{eqnarray*}

Using the fact that $\vert \Delta \left( \int_{0}^{\infty} e^{-k} \cos(k\Delta_{y} f ) \ dk \right) \vert \lesssim \vert \Delta_{y} \nabla_{x}f \vert^{2}+ \vert \Delta_{y} \Delta f \vert $ we find

\begin{eqnarray*}
\mathcal{T}_{5}&\lesssim& \Vert f  \Vert_{\dot H^{2}} \left(\int \frac{\Vert \nabla_x \delta_{y} f\Vert^3_{L^{6}}}{\vert y \vert^{4}} dy +  \int \frac{\Vert \nabla_x \delta_{y} f\Vert_{L^{\infty}}}{\vert y \vert^{3/2}} \frac{\Vert \Delta \delta_{y}f \Vert_{L^{2}}}{\vert y \vert^{3/2}} dy \right)  \\
\end{eqnarray*}

Hence, following the same step as \eqref{class} and \eqref{class2} one finds that
\begin{eqnarray*}
\mathcal{T}_{5}&\lesssim& \Vert f \Vert^{2}_{\dot H^{5/2}}\left(\Vert f  \Vert_{\dot H^{2}} +\Vert f  \Vert^2_{\dot H^{2}} \right)
\end{eqnarray*}

We have therefore obtain that all the less singular terms $\mathcal{T}_{i}$ for any $i=1,...,5$ are controlled as follows
\begin{eqnarray} \label{flop}
\sum_{i=1}^{5}\mathcal{T}_{i}&\lesssim& \Vert f \Vert^{2}_{\dot H^{5/2}}\left(\Vert f  \Vert_{\dot H^{2}} +\Vert f  \Vert^2_{\dot H^{2}} \right)
\end{eqnarray}

This ends the estimates of the less singular term and the proof of Lemma \ref{t1} is complete.

\qed \\

 In the next section, we shall estimate the more singular term. The analysis of the singular term requires much more effort, the first part consist in symmetrizing in a tricky way.

\section{Symmetrization and useful identities}

Throughout the article, we shall need to use some identities involving second finite differences. We collect all those identities in the following lemma

\begin{lemma} \label{astuce} Set $K_{y}f:=\frac{1}{1+\Delta^{2}_{y}{f}}$ and $\bar K_{y}f:=\frac{1}{1+\bar\Delta^{2}_{y}{f}}.$ The following equalities hold.
\begin{eqnarray}\label{formuleA+}
\hspace{-0,3cm}\nabla_{y}\left\{\arctan(\Delta_{y} f)+\arctan(\bar\Delta_{y} f) \right\} &=&-\frac{1}{2} S_y f D_y f K_y f \bar K_y f\nabla_{y}D_{y} f \nonumber\\
&+&  \frac{1}{2}\left(K_y f + \bar K_{y} f\right)\nabla_{y}S_{y} f
\end{eqnarray}
analogously, 
\begin{eqnarray} \label{formuleA-}
\hspace{-0,3cm}\nabla_{y}\left\{\arctan(\Delta_{y} f)-\arctan(\bar\Delta_{y} f) \right\} &=&-\frac{1}{2} S_y f D_y f K_y f \bar K_y f\nabla_{y}S_{y} f  \nonumber\\
&+&  \frac{1}{2}\left(K_{y}f+ \bar K_{y}f\right){\nabla_{y}D_{y} f}
\end{eqnarray}
\end{lemma}

\noindent {\bf{Proof of Lemma}} \ref{astuce} Set $A(x):=\nabla_{y}\left\{\arctan(\Delta_{y} f)+\arctan(\bar\Delta_{y} f) \right\}$. One may write that,
\begin{eqnarray} \label{formule1+}
A(x) &=& \frac{\nabla_{y}\Delta_{y} f}{1+\Delta^{2}_{y} f}\underbrace{-\frac{\nabla_{y}\Delta_{y} f}{1+\bar\Delta^{2}_{y} f}+\frac{\nabla_{y}\Delta_{y} f}{1+\bar\Delta^{2}_{y} f}
}_{=0} + \frac{\nabla_{y}\bar\Delta_{y} f}{1+\bar\Delta^{2}_{y} f} \nonumber \\
&=&\nabla_{y}\Delta_{y} f\frac{(\Delta_{y} f +\bar\Delta_{y} f)(\bar\Delta_{y} f -\Delta_{y} f )}{(1+\Delta^2_{y}f)(1+\bar\Delta^2_{y}f)}+\frac{\nabla_{y}S_{y} f}{1+\bar\Delta^{2}_{y} f} \nonumber \\
&=&-\nabla_{y}\Delta_{y} f\frac{S_y f D_y f}{(1+\Delta^2_{y}f)(1+\bar\Delta^2_{y}f)}+\frac{\nabla_{y}S_y f}{1+\bar\Delta^{2}_{y} f}
\end{eqnarray}
On the other hand,
\begin{eqnarray} \label{formule2+}
A(x) &=& \frac{\nabla_{y}\Delta_{y} f}{1+\Delta^{2}_{y} f}\underbrace{+\frac{\nabla_{y} \bar\Delta_{y} f}{1+\Delta^{2}_{y} f}-\frac{\nabla_{y} \bar\Delta_{y} f}{1+\Delta^{2}_{y} f}}_{=0} + \frac{\nabla_{y}\bar\Delta_{y} f}{1+\bar\Delta^{2}_{y} f} \nonumber \\
&=& \frac{\nabla_{y}S_{y} f}{1+\Delta^{2}_{y} f} + {\nabla_{y}\bar\Delta_{y} f}\left(\frac{1}{1+\bar\Delta^{2}_{y} f}-\frac{1}{1+\Delta^{2}_{y} f}\right) \nonumber \\
&=&\frac{\nabla_{y}S_y f}{1+\Delta^{2}_{y} f}+\nabla_{y} \bar\Delta_{y} f\frac{ S_{y} f D_{y}f}{(1+\Delta^2_{y}f)(1+\bar\Delta^2_{y}f)}.
\end{eqnarray}
Therefore, by combining \eqref{formule1+} and \eqref{formule2+} one gets \eqref{formuleA+}.  Analogously, set $B(x):=\nabla_{y}\left\{\arctan(\Delta_{y} f)-\arctan(\bar\Delta_{y} f) \right\}$, then we write that
\begin{eqnarray} 
 B(x)&=& \frac{\nabla_{y}\Delta_{y} f}{1+\Delta^{2}_{y} f}\underbrace{-\frac{\nabla_{y}\Delta_{y} f}{1+\bar\Delta^{2}_{y} f}+\frac{\nabla_{y}\Delta_{y} f}{1+\bar\Delta^{2}_{y} f}}_{=0} \underbrace{-\frac{\nabla_{y}\Delta_{y} g}{1+\bar\Delta^{2}_{y} g}+\frac{\nabla_{y}\Delta_{y} g}{1+\bar\Delta^{2}_{y} g}}_{=0} -\frac{\nabla_{y}\Delta_{y} g}{1+\Delta^{2}_{y} g} \nonumber \\
 \end{eqnarray}
 
\begin{eqnarray} \label{formule1}
 B(x)&=& \frac{\nabla_{y}\Delta_{y} f}{1+\Delta^{2}_{y} f}\underbrace{-\frac{\nabla_{y}\Delta_{y} f}{1+\bar\Delta^{2}_{y} f}+\frac{\nabla_{y}\Delta_{y} f}{1+\bar\Delta^{2}_{y} f}}_{=0} - \frac{\nabla_{y}\bar\Delta_{y} f}{1+\bar\Delta^{2}_{y} f} \nonumber \\
&=& {\nabla_{y}\Delta_{y} f}\left(\frac{1}{1+\Delta^{2}_{y} f}-\frac{1}{1+\bar\Delta^{2}_{y} f}\right)+\frac{\nabla_{y}\Delta_{y} f-\nabla_{y}\bar\Delta_{y} f}{1+\bar\Delta^{2}_{y} f} \nonumber\\
&=&-\nabla_{y}\Delta_{y} f\frac{(\Delta_{y} f +\bar\Delta_{y} f)(\Delta_{y} f -\bar\Delta_{y} f )}{(1+\Delta^2_{y}f)(1+\bar\Delta^2_{y}f)}+\frac{\nabla_{y}\Delta_{y} f-\nabla_{y}\bar\Delta_{y} f}{1+\bar\Delta^{2}_{y} f} \nonumber \\
&=&-\nabla_{y}\Delta_{y} f\frac{ S_{y} f D_{y}f}{(1+\Delta^2_{y}f)(1+\bar\Delta^2_{y}f)}+\frac{\nabla_{y}D_y f}{1+\bar\Delta^{2}_{y} f}
\end{eqnarray}
On the other hand we may write that
\begin{eqnarray} \label{formule2}
B(x) &=& \frac{\nabla_{y}\Delta_{y} f}{1+\Delta^{2}_{y} f}\underbrace{-\frac{\nabla_{y} \bar\Delta_{y} f}{1+\Delta^{2}_{y} f}+\frac{\nabla_{y} \bar\Delta_{y} f}{1+\Delta^{2}_{y} f}}_{=0} - \frac{\nabla_{y}\bar\Delta_{y} f}{1+\bar\Delta^{2}_{y} f} \nonumber \\
&=& \frac{\nabla_{y}\Delta_{y} f-\nabla_{y}\bar\Delta_{y} f}{1+\Delta^{2}_{y} f} + {\nabla_{y}\bar\Delta_{y} f}\left(\frac{1}{1+\Delta^{2}_{y} f}-\frac{1}{1+\bar\Delta^{2}_{y} f}\right) \nonumber \\
&=&\frac{\nabla_{y}D_y f}{1+\Delta^{2}_{y} f}-\nabla_{y} \bar\Delta_{y} f\frac{ S_{y} f D_{y}f}{(1+\Delta^2_{y}f)(1+\bar\Delta^2_{y}f)}.
\end{eqnarray}
Hence, combining \eqref{formule1} and \eqref{formule2} we get \eqref{formuleA-}. 

\qed

We shall need to compute gradients with respect to $y$ of the operators $S_y$ and $D_y$. The following lemma collects the main identities that we shall use.

\begin{lemma} \label{L1}We have
\begin{eqnarray} \label{diff}
\noindent D_{y} f=\frac{y}{\vert y \vert}.\left(\int_{0}^{1} \nabla\left( f(x+(r-1)y)  + f(x-(r-1)y)   -2 f(x) \right)  dr\right) +2 \nabla f . \frac{y}{\vert y \vert}
\end{eqnarray}

Moreover,
\begin{eqnarray} \label{ddiff}
y.\nabla_{y}D_{y} f&=& \frac{1}{\vert y \vert}\int_{0}^{1} y. s_{(r-1)y}  \nabla_{x}f \ dr +  \frac{y}{\vert y \vert} .{\nabla_x s_{y}f}
\end{eqnarray}

\begin{eqnarray} \label{som}
y.\nabla_{y}S_{y} f&=&\frac{1}{\vert y \vert}s_y f(x)  + \frac{y}{\vert y \vert}.\nabla_{x}\bar \delta_{y}f- \frac{y}{\vert y \vert}.\nabla_{x}\delta_{y}f
\end{eqnarray}

\end{lemma}

{\noindent{\bf{Proof of Lemma \ref{L1}}}} 
In order to prove \eqref{diff}, we first recall that since we have $D_{y} f=\frac{1}{\vert y \vert}(f(x+y)-f(x-y))$ one may readily check that
\begin{eqnarray*}
D_{y} f&=&\frac{1}{\vert y \vert}\int_{0}^{1} \left(\nabla f(x+(r-1)y).y  +\nabla f(x-(r-1)y).y   -2 \nabla f . y \right) \ dr +2 \nabla f . \frac{y}{\vert y \vert}
\end{eqnarray*}

The proof of \eqref{diff} is obtained as follows. First, we write that
\begin{eqnarray*}
\nabla_{y}D_{y} f&=&\nabla_{y} \frac{1}{\vert y \vert}(f(x+y)-f(x-y))+\frac{1}{\vert y \vert} \nabla_{y}(f(x+y)-f(x-y))\\
&=&\nabla_{y} \frac{1}{\vert y \vert}(f(x+y)-f(x-y))+\frac{1}{\vert y \vert} \nabla_{x}f(x+y)+\nabla_{x}f(x-y)\\
&=&\nabla_{y} \left(\frac{1}{\vert y \vert}\right) \int_{0}^{1} \nabla f(x+(r-1)y).y + \nabla f(x-(r-1)y).y   -2 \nabla f (x). y  \ dr\\
&+&  \frac{\nabla_x f(x+y) + \nabla_{x} f(x-y)}{\vert y \vert}+  \nabla_{y} \left(\frac{2}{\vert y \vert}\right)(\nabla f (x). y)
\end{eqnarray*}
Using that $\nabla_{y} \frac{1}{\vert y \vert}.y=-\frac{1}{\vert y \vert}$ and recalling that $s_{y}$ denotes the second finite difference operator, we immediately find that \\
 \begin{eqnarray*}
  y.\nabla_{y}D_{y} f&=&- \frac{1}{\vert y \vert}\int_{0}^{1} (\nabla f(x+(r-1)y).y) + (\nabla f(x-(r-1)y).y)   -{2} \nabla f (x).y \ dr \nonumber \\
&& \ + \  \frac{\nabla_x f(x+y) + \nabla_{x} f(x-y)-2\nabla_x f (x)}{\vert y \vert}.y \\
&=&  \frac{1}{\vert y \vert}\int_{0}^{1} y. s_{(r-1)y}  \nabla_{x}f \ dr 
+  \frac{y.\nabla_x s_{y}f}{\vert y \vert},
\end{eqnarray*}
which is the desired identity \eqref{ddiff}. Let us prove \eqref{som}. We observe that
\begin{eqnarray*} 
\nabla_{y}S_{y} f&=&-\nabla_{y} \frac{1}{\vert y \vert} (f(x+y)+f(x-y)-2f(x))-\frac{1}{\vert y \vert} \nabla_{y}f(x+y)+\frac{1}{\vert y \vert} \nabla_{y}f(x-y)\\
&=&-\nabla_{y} \frac{1}{\vert y \vert} (f(x+y)+f(x-y)-2f(x))-\frac{1}{\vert y \vert} \nabla_{x}f(x+y)+\frac{1}{\vert y \vert} \nabla_{x}f(x-y) \\
&=&-\nabla_{y} \frac{1}{\vert y \vert} (f(x+y)+f(x-y)-2f(x))-\frac{1}{\vert y \vert} \nabla_{x}(f(x+y)-f(x))\\
&&+ \ \frac{1}{\vert y \vert} \nabla_{x}(f(x-y)-f(x))\\
&=&s_y f\nabla_{y} \frac{1}{\vert y \vert}  +\frac{1}{\vert y \vert} \nabla_{x}\bar \delta_{y}f -\frac{1}{\vert y \vert} \nabla_{x}\delta_{y}f\\
&=&s_y f\nabla_{y} \frac{1}{\vert y \vert}  -\nabla_{x} D_{y}f.
\end{eqnarray*}
Hence,
\begin{eqnarray} \label{dsy}
\nabla_{y}S_{y} f&=&s_y f\nabla_{y} \frac{1}{\vert y \vert}  -\nabla_{x} D_{y}f. 
\end{eqnarray}
Therefore,
\begin{eqnarray*} 
y.\nabla_{y}S_{y} f&=&\frac{1}{\vert y \vert}s_y f(x)  + \frac{y}{\vert y \vert}\nabla_{x}\bar \delta_{y}f- \frac{y}{\vert y \vert}\nabla_{x}\delta_{y}f
\end{eqnarray*}

Which is the wanted identity \eqref{som}. This ends the proof of Lemma \ref{astuce}.

\qed

Finally, we state an easy lemma that will be systematically used throughout the article.

\begin{lemma} \label{facile}
Let $r>0$, we have
  \begin{equation} 
     \left\vert \nabla. \frac{x}{\vert x \vert^{r}} \right\vert \lesssim \frac{1}{\vert x \vert^{r}} 
      \end{equation}
\end{lemma}

\noindent {\bf{Proof of Lemma}} \ref{facile}   A direct computation leads to the estimate.
\qed

\section{Estimates of the most singular term : $\mathcal{S}=\displaystyle\sum_{i=1}^4 S_i$}  \label{sing}

{  {In order to control the most singular terms, that is when the Laplacian opeator falls onto the non-oscillatory term, one has to balance the regularity in both $x$ and $y$. This is mainly because of the fact that if we only balance the derivatives in the spatial variable then this amounts to control terms whose regularity in Sobolev or Besov spaces are higher than 1. Recall that controlling such terms require to have second finite order difference. The main goal of the next Lemma is to force the appearance of these terms, in other words, we need to symmetrize the terms.}}

\subsection{{{Algebraic decomposition of the most singular term}} : $\mathcal{S}$} 
 Set $D_y f := \Delta_{y} f-\bar\Delta_{y} f$ and $S_y f := \Delta_{y} f+\bar\Delta_{y} f$. We shall prove the following Lemma

\begin{lemma} \label{dec} (symmetrization of the singular term) We have the following decomposition
\begin{eqnarray*}
 \mathcal{S}&=&\frac{1}{2} \int \nabla_{x}\Delta D_{y} f . \frac{y}{\vert y \vert^2} \sin(\frac{1}{2}(\arctan(\Delta_{y} f )+\arctan(\bar\Delta_{y} f )))     \times \\
 &&\sin(\frac{1}{2}(\arctan(\Delta_{y} f )-\arctan(\bar\Delta_{y} f )))\int_{0}^{\infty} e^{-k} \sin(\frac{k}{2} S_y f) \sin(\frac{k}{2} D_y f)    \ dk  \ dy  \\
  &+& \int \nabla_{x} \Delta\Delta_{y} f . \frac{y}{\vert y \vert^2} \sin(\frac{1}{2}(\arctan(\Delta_{y} f )+\arctan(\bar\Delta_{y} f )))    \times \\
&&\sin(\frac{1}{2}(\arctan(\Delta_{y} f )-\arctan(\bar\Delta_{y} f )))\int_{0}^{\infty} e^{-k} \cos(\frac{k}{2} S_y f) \cos(\frac{k}{2} D_y f)   \ dk  \ dy  \\
&+&\int \nabla_{x} \Delta\bar\Delta_{y} f . \frac{y}{\vert y \vert^2} 
\cos(\frac{1}{2}(\arctan(\bar\Delta_{y} f + \arctan(\Delta_{y} f )) \times \\
&&\cos(\frac{1}{2}(\arctan(\bar\Delta_{y} f )- \arctan(\Delta_{y} f )))\int_{0}^{\infty} e^{-k} \sin(\frac{k}{2} S_y f) \sin(\frac{k}{2} D_y f)  \ dk  \ dy \\
&+& \frac{1}{2} \int  \nabla_{x}\Delta D_{y} f . \frac{y}{\vert y \vert^2} \cos(\frac{1}{2}(\arctan(\Delta_{y} f )+\arctan(\bar\Delta_{y} f )))   \times  \\
&& \cos(\frac{1}{2}(\arctan(\Delta_{y} f )-\arctan(\bar\Delta_{y} f )))\int_{0}^{\infty} e^{-k} (\cos(\frac{k}{2}D_y f )(\cos(\frac{k}{2}S_y f )  \ dk  \ dy  \\
&:=& \sum_{i=1}^4 \mathcal{S}_{i}
\end{eqnarray*}
\end{lemma}

\noindent {\bf{Proof of Lemma \ref{dec}  }}  {  {We start by symmetrizing the non-oscillatory part, that is we write}}
\begin{eqnarray*}
 \mathcal{S} &=&\int   (\Delta\nabla_{x}\Delta_{y}  f - \Delta\nabla_{x}\bar\Delta_{y} f). \frac{y}{\vert y \vert^2} \cos(\arctan(\Delta_{y} f ))\int_{0}^{\infty} e^{-k} \cos(k\Delta_{y} f ) \ dk  \ dy  \\
  &-& \int \nabla_{x}\Delta\Delta_{y} f . \frac{y}{\vert y \vert^2} \cos(\arctan(\bar\Delta_{y} f )) \int_{0}^{\infty} e^{-k} \cos(k\bar\Delta_{y} f )   \ dk  \ dy  
  \end{eqnarray*}
  Then, by doing a change of variable ($y\rightarrow -y$), one finds 
  \begin{eqnarray*}
\mathcal{S} &=&   \int   (\Delta\nabla_{x}\Delta_{y}  f - \Delta\nabla_{x}\bar\Delta_{y} f). \frac{y}{\vert y \vert^2} \cos(\arctan(\Delta_{y} f )) \int_{0}^{\infty} e^{-k} \cos(k\Delta_{y} f ) \ dk  \ dy  \\
  &-& \int \nabla_{x}\Delta\Delta_{y} f . \frac{y}{\vert y \vert^2} \left(\cos(\arctan(\bar\Delta_{y} f ))-\cos(\arctan(\Delta_{y} f ))\right) \int_{0}^{\infty} e^{-k} \cos(k\bar\Delta_{y} f )   \ dk  \ dy  \\
    &-& \int \nabla_{x}\Delta\Delta_{y} f . \frac{y}{\vert y \vert^2} \cos(\arctan(\Delta_{y} f )) \int_{0}^{\infty} e^{-k} \cos(k\bar\Delta_{y} f )   \ dk  \ dy.  \\
    \end{eqnarray*}
    Then, we find that
    \begin{eqnarray*}
   \mathcal{S} &=&   \int   (\Delta\nabla_{x}\Delta_{y}  f - \Delta\nabla_{x}\bar\Delta_{y} f). \frac{y}{\vert y \vert^2} \cos(\arctan(\Delta_{y} f ))\int_{0}^{\infty} e^{-k} \cos(k\Delta_{y} f ) \ dk  \ dy  \\
  &-& \int \nabla_{x}\Delta\Delta_{y} f . \frac{y}{\vert y \vert^2} \left(\cos(\arctan(\bar\Delta_{y} f ))-\cos(\arctan(\Delta_{y} f ))\right) \int_{0}^{\infty} e^{-k} \cos(k\bar\Delta_{y} f )   \ dk  \ dy  \\
    &-& \int \nabla_{x}\Delta\Delta_{y} f . \frac{y}{\vert y \vert^2} \cos(\arctan(\Delta_{y} f ))\int_{0}^{\infty} e^{-k} (\cos(k\bar\Delta_{y} f )-\cos(k\Delta_{y} f ) )  \ dk  \ dy  \\
    &-&\int \nabla_{x}\Delta\Delta_{y} f . \frac{y}{\vert y \vert^2} \cos(\arctan(\Delta_{y} f )) \int_{0}^{\infty} e^{-k} \cos(k\Delta_{y} f )   \ dk  \ dy
     \end{eqnarray*}
    Noticing that the last term is nothing but $-\mathcal{S}(t)$ one finds
    
\begin{eqnarray*}
 \mathcal{S}&=& \frac{1}{2}  \int   (\nabla_{x}\Delta\Delta_{y} f-\nabla_{x}\Delta \bar\Delta_{y} f) . \frac{y}{\vert y \vert^2} \cos(\arctan(\Delta_{y} f )) \int_{0}^{\infty} e^{-k} \cos(k\Delta_{y} f )  \ dk  \ dy  \\
  &-&\frac{1}{2}\int \nabla_{x}\Delta\bar\Delta_{y} f . \frac{y}{\vert y \vert^2} \cos(\arctan(\bar\Delta_{y} f )) \int_{0}^{\infty} e^{-k} (\cos(k\bar\Delta_{y} f )-\cos(k\Delta_{y} f ))   \ dk  \ dy  \\
&-&\frac{1}{2}\int \nabla_{x}\Delta\bar\Delta_{y} f . \frac{y}{\vert y \vert^2} (\cos(\arctan(\bar\Delta_{y} f ))-\cos(\arctan(\Delta_{y} f )))\int_{0}^{\infty} e^{-k} \cos(k\Delta_{y} f )   \ dk  \ dy  \\
&=& \mathcal{O}_1 + \mathcal{O}_2 + \mathcal{O}_3.
\end{eqnarray*}
Then, one observes that, \\
\begin{eqnarray*}
&& \hspace{-1cm} \mathcal{S}= \frac{1}{2} \int (\nabla_{x}\Delta\Delta_{y} f-\nabla_{x}\Delta \bar\Delta_{y} f) . \frac{y}{\vert y \vert^2} \cos(\arctan(\Delta_{y} f )) \int_{0}^{\infty} e^{-k} \cos(k\Delta_{y} f )  \ dk  \ dy \\
  &&-\frac{1}{2}\int   \nabla_{x}\Delta\bar\Delta_{y} f . \frac{y}{\vert y \vert^2} \left(\cos(\arctan(\bar\Delta_{y} f ))+\cos(\arctan(\Delta_{y} f ))\right) \times \\
  &&\int_{0}^{\infty} e^{-k} (\cos(k\bar\Delta_{y} f )-\cos(k\Delta_{y} f ))  \ dk  \ dy \\
  &&+\frac{1}{2}\int  \nabla_{x}(\Delta\bar\Delta_{y} f-\Delta\Delta_{y} f) . \frac{y}{\vert y \vert^2} \cos(\arctan(\Delta_{y} f )) \\
  && \int_{0}^{\infty} e^{-k} (\cos(k\bar\Delta_{y} f )-\cos(k\Delta_{y} f ))   \ dk  \ dy  \\
  &&+\underbrace{\frac{1}{2}\int \nabla_{x}\Delta\Delta_{y} f . \frac{y}{\vert y \vert^2} \cos(\arctan(\Delta_{y} f )) \int_{0}^{\infty} e^{-k} (\cos(k\bar\Delta_{y} f )-\cos(k\Delta_{y} f ))   \ dk  \ dy }_{=-\mathcal{O}_2} \\
&&-\frac{1}{2} \int \nabla_{x}(\Delta\bar\Delta_{y} f-\Delta\Delta_{y} f) . \frac{y}{\vert y \vert^2} (\cos(\arctan(\bar\Delta_{y} f ))-\cos(\arctan(\Delta_{y} f )))\\
&& \times \int_{0}^{\infty} e^{-k} (\cos(k\Delta_{y} f )   \ dk  \ dy  \\
&&-\frac{1}{2} \int \nabla_{x}\Delta\Delta_{y} f . \frac{y}{\vert y \vert^2} (\cos(\arctan(\bar\Delta_{y} f ))-\cos(\arctan(\Delta_{y} f ))) \\
&& \times \int_{0}^{\infty} e^{-k} (\cos(k\Delta_{y} f + \cos(k\bar\Delta_{y} f) \ dk  \ dy  \\
&&+\underbrace{\frac{1}{2} \int \nabla_{x}\Delta\Delta_{y} f . \frac{y}{\vert y \vert^2} (\cos(\arctan(\bar\Delta_{y} f ))-\cos(\arctan(\Delta_{y} f )))\int_{0}^{\infty} e^{-k}  \cos(k\bar\Delta_{y} f)  \ dk  \ dy }_{=-\mathcal{O}_3} 
\end{eqnarray*}
Therefore,
\begin{eqnarray*}
\mathcal{S}&=& \frac{1}{2} \int  (\nabla_{x}\Delta\Delta_{y} f-\nabla_{x}\Delta \bar\Delta_{y} f) . \frac{y}{\vert y \vert^2} \cos(\arctan(\Delta_{y} f )) \int_{0}^{\infty} e^{-k} \cos(k\Delta_{y} f )  \ dk  \ dy \ dx \\
  &-&\frac{1}{4}\int \Delta f \int \nabla_{x}\Delta\bar\Delta_{y} f . \frac{y}{\vert y \vert^2} \left(\cos(\arctan(\bar\Delta_{y} f ))+\cos(\arctan(\Delta_{y} f ))\right) \\
  && \times \int_{0}^{\infty} e^{-k} (\cos(k\bar\Delta_{y} f )-\cos(k\Delta_{y} f ))   \ dk  \ dy  \\
  &+&\frac{1}{4} \int \nabla_{x}(\Delta\bar\Delta_{y} f-\Delta\Delta_{y} f) . \frac{y}{\vert y \vert^2} \cos(\arctan(\Delta_{y} f )) \\
  && \times \int_{0}^{\infty} e^{-k} (\cos(k\bar\Delta_{y} f )-\cos(k\Delta_{y} f ))   \ dk  \ dy  \\
&-&\frac{1}{4}\int \Delta f \int \nabla_{x}(\Delta\bar\Delta_{y} f-\Delta\Delta_{y} f) . \frac{y}{\vert y \vert^2} (\cos(\arctan(\bar\Delta_{y} f ))-\cos(\arctan(\Delta_{y} f )))\\
&& \times \int_{0}^{\infty} e^{-k} (\cos(k\Delta_{y} f )  \ dk  \ dy  \\
&-&\frac{1}{4} \int \Delta\nabla_{x}\Delta_{y} f . \frac{y}{\vert y \vert^2} (\cos(\arctan(\bar\Delta_{y} f ))-\cos(\arctan(\Delta_{y} f ))) \\
&\times& \int_{0}^{\infty} e^{-k} (\cos(k\Delta_{y} f + \cos(k\bar\Delta_{y} f)   \ dk  \ dy.  
\end{eqnarray*}
Then, by noticing that the third and fourth terms cancel out, one finds that
\begin{eqnarray*}
 \mathcal{S}&=& \frac{1}{2}  \int \nabla_{x}\Delta\Delta_{y} f-\nabla_{x} \Delta\bar\Delta_{y} f) . \frac{y}{\vert y \vert^2} \cos(\arctan(\Delta_{y} f )) \int_{0}^{\infty} e^{-k} \cos(k\Delta_{y} f ) \ dk  \ dy  \\
  &-&\frac{1}{4}\int  \nabla_{x}\Delta\bar\Delta_{y} f . \frac{y}{\vert y \vert^2} \left(\cos(\arctan(\bar\Delta_{y} f ))+\cos(\arctan(\Delta_{y} f ))\right) \\
  &\times&\int_{0}^{\infty} e^{-k} (\cos(k\bar\Delta_{y} f )-\cos(k\Delta_{y} f ))  \ dk  \ dy  \\
&-&\frac{1}{4} \int \Delta\nabla_{x}\Delta_{y} f . \frac{y}{\vert y \vert^2} (\cos(\arctan(\bar\Delta_{y} f ))-\cos(\arctan(\Delta_{y} f ))) \\
&&\times \int_{0}^{\infty} e^{-k} (\cos(k\Delta_{y} f + \cos(k\bar\Delta_{y} f)   \ dk  \ dy.
\end{eqnarray*}
Then, one observes that the first term, namely
\begin{eqnarray} \label{Lint}
L:=\frac{1}{2} \int (\nabla_{x} \Delta\Delta_{y} f-\nabla_{x} \Delta \bar\Delta_{y} f) . \frac{y}{\vert y \vert^2} \cos(\arctan(\Delta_{y} f )) \int_{0}^{\infty} e^{-k} \cos(k\Delta_{y} f )  \ dk  \ dy, \\
\end{eqnarray}
may be rewritten as
\begin{eqnarray*}
L&=&  \frac{1}{8}  \int  (\nabla_{x} \Delta\Delta_{y} f-\nabla_{x} \Delta \bar\Delta_{y} f) . \frac{y}{\vert y \vert^2} (\cos(\arctan(\Delta_{y} f ))- \cos(\arctan(\bar\Delta_{y} f ))) \\
  &\times& \int_{0}^{\infty} e^{-k} (\cos(k\Delta_{y} f )-\cos(k\bar\Delta_{y} f ) ) \ dk  \ dy  \\
  &+& \frac{1}{8}  \int (\nabla_{x} \Delta\Delta_{y} f-\nabla_{x} \Delta \bar\Delta_{y} f) . \frac{y}{\vert y \vert^2}  \left(\cos(\arctan(\Delta_{y} f ))+\cos(\arctan(\bar\Delta_{y} f ))\right) \\
  &\times& \int_{0}^{\infty} e^{-k}( \cos(k\bar\Delta_{y} f )+\cos(k\Delta_{y} f ) ) \ dk  \ dy.  \\
\end{eqnarray*}
To prove this idendity, the idea is to try to symmetrize the integral \eqref{Lint} . To this end, one writes 
\begin{eqnarray*}
 L &=&\frac{1}{2}\int (\nabla_{x} \Delta\Delta_{y} f-\nabla_{x} \Delta \bar\Delta_{y} f) . \frac{y}{\vert y \vert^2} (\cos(\arctan(\Delta_{y} f ))- \cos(\arctan(\bar\Delta_{y} f ))) \\
  && \times \int_{0}^{\infty} e^{-k} \cos(k\Delta_{y} f )  \ dk  \ dy  \\
  &+& \frac{1}{2} \int (\nabla_{x} \Delta\Delta_{y} f-\nabla_{x} \Delta \bar\Delta_{y} f) . \frac{y}{\vert y \vert^2}  \cos(\arctan(\Delta_{y} f ))\\ && \times \int_{0}^{\infty} e^{-k}( \cos(k\bar\Delta_{y} f )+\cos(k\Delta_{y} f ) ) \ dk  \ dy \\
  &-&\frac{1}{2}  \int (\nabla_{x} \Delta\Delta_{y} f-\nabla_{x} \Delta \bar\Delta_{y} f) . \frac{y}{\vert y \vert^2}  \cos(\arctan(\Delta_{y} f )) \int_{0}^{\infty} e^{-k} \cos(k\Delta_{y} f )  \ dk  \ dy.  
  \end{eqnarray*}
 Noticing that the last integral is equal to $-L$ one may symmetrize and find that

  \begin{eqnarray*}
  L&=& \frac{1}{8}  \int  (\nabla_{x} \Delta\Delta_{y} f-\nabla_{x} \Delta \bar\Delta_{y} f) . \frac{y}{\vert y \vert^2} (\cos(\arctan(\Delta_{y} f ))- \cos(\arctan(\bar\Delta_{y} f ))) \\
  &\times& \int_{0}^{\infty} e^{-k} (\cos(k\Delta_{y} f )-\cos(k\bar\Delta_{y} f ) ) \ dk  \ dy  \\
  &+& \frac{1}{8}  \int (\nabla_{x} \Delta\Delta_{y} f-\nabla_{x} \Delta \bar\Delta_{y} f) . \frac{y}{\vert y \vert^2}  \left(\cos(\arctan(\Delta_{y} f ))+\cos(\arctan(\bar\Delta_{y} f ))\right) \\
  &\times& \int_{0}^{\infty} e^{-k}( \cos(k\bar\Delta_{y} f )+\cos(k\Delta_{y} f ) ) \ dk  \ dy  \
  \end{eqnarray*}

Hence,
\begin{eqnarray} \label{dec}
\mathcal{S}&=& \nonumber\frac{1}{8}  \int  (\nabla_{x} \Delta\Delta_{y} f-\nabla_{x} \Delta \bar\Delta_{y} f) . \frac{y}{\vert y \vert^2} (\cos(\arctan(\Delta_{y} f ))- \cos(\arctan(\bar\Delta_{y} f ))) \\
 &\times& \int_{0}^{\infty} e^{-k} (\cos(k\Delta_{y} f )-\cos(k\bar\Delta_{y} f ) ) \ dk  \ dy  \\
\nonumber  &+& \frac{1}{8}  \int (\nabla_{x} \Delta\Delta_{y} f-\nabla_{x} \Delta \bar\Delta_{y} f) . \frac{y}{\vert y \vert^2}  \left(\cos(\arctan(\Delta_{y} f ))+\cos(\arctan(\bar\Delta_{y} f ))\right) \\
 \nonumber &\times& \int_{0}^{\infty} e^{-k}( \cos(k\bar\Delta_{y} f )+\cos(k\Delta_{y} f ) ) \ dk  \ dy  \\
  \nonumber&-&\frac{1}{4} \int \nabla_{x} \Delta\bar\Delta_{y} f . \frac{y}{\vert y \vert^2} \left(\cos(\arctan(\bar\Delta_{y} f ))+\cos(\arctan(\Delta_{y} f ))\right) \times \\
 \nonumber && \int_{0}^{\infty} e^{-k} (\cos(k\bar\Delta_{y} f )-\cos(k\Delta_{y} f ))  \ dk  \ dy \\
\nonumber&-&\frac{1}{4}\int \nabla_{x} \Delta\Delta_{y} f . \frac{y}{\vert y \vert^2} (\cos(\arctan(\bar\Delta_{y} f ))-\cos(\arctan(\Delta_{y} f )))\times\\
\nonumber &&\int_{0}^{\infty} e^{-k} (\cos(k\Delta_{y} f) + \cos(k\bar\Delta_{y} f)) \ dk  \ dy  
\end{eqnarray}
Finally, by denoting $D_y f = \Delta_{y} f-\bar\Delta_{y} f$ and $S_y f = \Delta_{y} f+\bar\Delta_{y} f$ along with the use of trigonometry identities, we obtain the desired decompostion 
\begin{eqnarray*}
 \mathcal{S}&=&\frac{1}{2} \int \nabla_{x}\Delta D_{y} f . \frac{y}{\vert y \vert^2} \sin(\frac{1}{2}(\arctan(\Delta_{y} f )+\arctan(\bar\Delta_{y} f )))    \times \\
 &&\sin(\frac{1}{2}(\arctan(\Delta_{y} f )-\arctan(\bar\Delta_{y} f ))) \int_{0}^{\infty} e^{-k} \sin(\frac{k}{2} S_y f) \sin(\frac{k}{2} D_y f)    \ dk  \ dy  \\
 &+& \int \nabla_{x} \Delta\Delta_{y} f . \frac{y}{\vert y \vert^2} \sin(\frac{1}{2}(\arctan(\Delta_{y} f )+\arctan(\bar\Delta_{y} f )))    \times \\
&&\sin(\frac{1}{2}(\arctan(\Delta_{y} f )-\arctan(\bar\Delta_{y} f )))\int_{0}^{\infty} e^{-k} \cos(\frac{k}{2} S_y f) \cos(\frac{k}{2} D_y f)   \ dk  \ dy  \\
  &+&\int \nabla_{x} \Delta\bar\Delta_{y} f . \frac{y}{\vert y \vert^2} 
\cos(\frac{1}{2}(\arctan(\bar\Delta_{y} f) + \arctan(\Delta_{y} f ))\times \\
&&\cos(\frac{1}{2}(\arctan(\bar\Delta_{y} f )- \arctan(\Delta_{y} f )) \int_{0}^{\infty} e^{-k} \sin(\frac{k}{2} S_y f) \sin(\frac{k}{2} D_y f)  \ dk  \ dy \\
&+& \frac{1}{2} \int  \nabla_{x}\Delta D_{y} f . \frac{y}{\vert y \vert^2} \cos(\frac{1}{2}(\arctan(\Delta_{y} f )+\arctan(\bar\Delta_{y} f ))) \int_{0}^{\infty} e^{-k}   \\
&&\cos(\frac{1}{2}(\arctan(\Delta_{y} f )-\arctan(\bar\Delta_{y} f )))  (\cos(\frac{k}{2}(D_y f ))(\cos(\frac{k}{2}(S_y f ))  \ dk  \ dy\\
&:=& \sum_{i=1}^4 \mathcal{\sigma}_{i}(t)
\end{eqnarray*}
This ends the proof the Lemma 8.1 \\
\qed

In the sequel we shall use the notation $\mathcal{S}_{i}:=(\Delta f ,\mathcal{\sigma}_{i})$

\subsection{Estimate of $\mathcal{S}_{1}$} \ {  {In this subsection, we are going to prove the following control of $\mathcal{S}_{1}.$ The main idea will be to transfer the regularity from the singular non-oscillatory term onto the oscillatory terms by using the regularity in $x$ and then write it in terms of $y$. In this subsection, we are going to proof the following estimate for  $\mathcal{S}_{1}$.
\begin{lemma} \label{s1} The term $\mathcal{S}_1$ is estimated as follows
\begin{eqnarray} 
\mathcal{S}_1=\sum_{i=1}^{7} \mathcal{S}_{1,i} \lesssim \Vert  f \Vert^2_{\dot H^{5/2}} \Vert f \Vert_{\dot H^{2}}
\end{eqnarray}
\end{lemma}

\noindent{\bf{Proof of Lemma}} \ref{s1} \\
}}

To estimate $\mathcal{S}_{1}$ one first needs to integrate by parts 
\begin{eqnarray*}
\mathcal{S}_{1} &=& \frac{1}{2}\int \Delta f \int \frac{\nabla_{x}(\delta_{y} \Delta f-\bar\delta_{y} \Delta f)}{\vert y \vert} . \frac{y}{\vert y \vert^2} \sin(\frac{1}{2}(\arctan(\Delta_{y} f )+\arctan(\bar\Delta_{y} f )))    \\
  &&     \sin(\frac{1}{2}(\arctan(\Delta_{y} f )-\arctan(\bar\Delta_{y} f ))) \int_{0}^{\infty} e^{-k}
   \sin(\frac{k}{2}S_y f) \sin(\frac{k}{2}D_y f)   \ dk  \ dy \ dx. \\
  \end{eqnarray*}

In order to make appear the more favorable second finite order differences it suffices to observe for instance that $\nabla_{x}(\delta_{y} f-\bar\delta_{y} f)=-\nabla_{y}(\delta_{y} f+ \bar\delta_{y} f)$. Hence, we may integrate by parts (in $y$) and we find that

 \begin{eqnarray} \label{t1}
\mathcal{S}_{1} &=& \frac{1}{2}\int \Delta f \int {\Delta s_y f}  \nabla_y . \left( \frac{y}{\vert y \vert^3} \sin(\frac{1}{2}(\arctan(\Delta_{y} f )+\arctan(\bar\Delta_{y} f )))  \right.  \\
  &&  \left.   \sin(\frac{1}{2}(\arctan(\Delta_{y} f )-\arctan(\bar\Delta_{y} f ))) \int_{0}^{\infty} e^{-k}
   \sin(\frac{k}{2}S_y f) \sin(\frac{k}{2}D_y f)   \ dk  \ dy \ dx \right)\nonumber\\ \nonumber
  \end{eqnarray}

Then, by using the identities \eqref{formuleA+} and \eqref{formuleA-} we  find that 

\begin{eqnarray*}
\mathcal{S}_{1} &=&\frac{1}{2}\int \Delta f \int {\Delta s_{y} f} \left(\nabla.\frac{y}{\vert y \vert^3}\right)  \sin(\frac{1}{2}(\arctan(\Delta_{y} f )+\arctan(\bar\Delta_{y} f )))     \\
  &&     \sin(\frac{1}{2}(\arctan(\Delta_{y} f )-\arctan(\bar\Delta_{y} f ))) \int_{0}^{\infty} e^{-k}
    \sin(\frac{k}{2}S_y f) \sin(\frac{k}{2}D_y f)  \ dk  \ dy \ dx \\
  &+&\frac{1}{8}\int \Delta f \int \frac{\Delta s_{y} f}{\vert y \vert^3}   S_y f D_y f K_y f \bar K_y f \ y.\nabla_{y}D_{y} f  \cos(\frac{1}{2}(\arctan(\Delta_{y} f )+\arctan(\bar\Delta_{y} f )))      \\
  &&     \sin(\frac{1}{2}(\arctan(\Delta_{y} f )-\arctan(\bar\Delta_{y} f ))) \int_{0}^{\infty} e^{-k}
    \sin(\frac{k}{2}S_y f) \sin(\frac{k}{2}D_y f)   \ dk  \ dy \ dx \\
  &-&\frac{1}{8}\int \Delta f \int \frac{\Delta s_{y} f}{\vert y \vert^3} \left(K_y f + \bar K_{y} f\right)y.\nabla_{y}S_{y} f \cos(\frac{1}{2}(\arctan(\Delta_{y} f )+\arctan(\bar\Delta_{y} f )))     \\
  &&    \sin(\frac{1}{2}(\arctan(\Delta_{y} f )-\arctan(\bar\Delta_{y} f ))) \int_{0}^{\infty} e^{-k}
    \sin(\frac{k}{2}S_y f) \sin(\frac{k}{2}D_y f) \ dk  \ dy \ dx \\
   &+&\frac{1}{8}\int \Delta f \int \frac{\Delta s_{y} f}{\vert y \vert^3}   S_y f D_y f K_y f \bar K_y f \ y.\nabla_{y}S_{y} f  \sin(\frac{1}{2}(\arctan(\Delta_{y} f )+\arctan(\bar\Delta_{y} f )))      \\
  &&     \cos(\frac{1}{2}(\arctan(\Delta_{y} f )-\arctan(\bar\Delta_{y} f ))) \int_{0}^{\infty} e^{-k}
    \sin(\frac{k}{2}S_y f) \sin(\frac{k}{2}D_y f)   \ dk  \ dy \ dx \\
  &-&\frac{1}{8}\int \Delta f \int \frac{\Delta s_{y} f}{\vert y \vert^3} \left(K_y f + \bar K_{y} f\right)y.\nabla_{y}D_{y} f \sin(\frac{1}{2}(\arctan(\Delta_{y} f )+\arctan(\bar\Delta_{y} f )))     \\
  &&    \cos(\frac{1}{2}(\arctan(\Delta_{y} f )-\arctan(\bar\Delta_{y} f ))) \int_{0}^{\infty} e^{-k}
    \sin(\frac{k}{2}S_y f) \sin(\frac{k}{2}D_y f) \ dk  \ dy \ dx \\
 &+& \frac{1}{4}\int \Delta f \int  \frac{\Delta s_{y} f}{\vert y \vert^3}  \sin(\frac{1}{2}(\arctan(\Delta_{y} f )+\arctan(\bar\Delta_{y} f )))  y.\nabla_{y}S_y f    \\
  &&  \cos(\frac{1}{2}(\arctan(\Delta_{y} f )-\arctan(\bar\Delta_{y} f )))     \int_{0}^{\infty} k e^{-k} \cos(\frac{k}{2} S_y f) \sin(\frac{k}{2} D_y f)     \ dk  \ dy \ dx \\
  &+& \frac{1}{4}\int \Delta f \int  \frac{\Delta s_{y} f}{\vert y \vert^3}  \sin(\frac{1}{2}(\arctan(\Delta_{y} f )+\arctan(\bar\Delta_{y} f )))    \\
  &&     \sin(\frac{1}{2}(\arctan(\Delta_{y} f )-\arctan(\bar\Delta_{y} f )))  y.\nabla_{y}D_y f \int_{0}^{\infty} k e^{-k} \sin(\frac{k}{2} S_y f) \cos(\frac{k}{2} D_y f)   \\
  && \ dk  \ dy \ dx \\
  &:=& \sum_{i}^7 \mathcal{S}_{1,i}
  \end{eqnarray*}
 
\subsubsection{{\bf{Estimate of}} $\mathcal{S}_{1,1}$ }
We have
\begin{eqnarray*}
\mathcal{S}_{1,1} &=&\frac{1}{2}\int  \Delta f \int \left( \Delta s_{y}f \right)\left(\nabla.\frac{y}{\vert y \vert^3}\right) \sin(\frac{1}{2}(\arctan(\Delta_{y} f )+\arctan(\bar\Delta_{y} f )))     \\
  &&     \sin(\frac{1}{2}(\arctan(\Delta_{y} f )-\arctan(\bar\Delta_{y} f ))) \int_{0}^{\infty} e^{-k} \sin(\frac{k}{2}S_y f) \sin(\frac{k}{2}D_y f)  \ dk  \ dy \ dx. \\
  \end{eqnarray*}
  Since $\left\vert \nabla.\frac{y}{\vert y \vert^3} \right \vert \lesssim \frac{1}{\vert y \vert^3} $ then we find  that,
  \begin{eqnarray*}
\mathcal{S}_{1,1}  &\leq& \frac{\Gamma(2)}{4}\Vert f \Vert_{\dot H^{2}} 
\int \frac{\Vert \Delta s_{y}f \Vert_{L^{2}}}{\vert y \vert^{3/2}} \frac{\Vert s_{y} f \Vert_{L^{\infty}}}{\vert y \vert^{5/2}} \ dy \\
  &\leq& \frac{1}{2} \Vert f \Vert_{\dot H^{2}} \Vert \Delta f\Vert_{\dot B^{1/2}_{2,2}}
   \Vert f \Vert_{\dot B^{3/2}_{\infty,2}}
\end{eqnarray*}
Then, since $\dot H^{5/2} \hookrightarrow \dot B^{3/2}_{\infty,2}$ one finds
\begin{eqnarray*}
\mathcal{S}_{1,1} \lesssim  \Vert f\Vert^2_{\dot H^{5/2}} \Vert f \Vert_{\dot H^{2}}
  \end{eqnarray*}
  
  \subsubsection{{\bf{Estimate of}} $\mathcal{S}_{1,2}$ }
We have
\begin{eqnarray} \label{relou}
\mathcal{S}_{1,2} &=&\frac{1}{8}\int \Delta f \int \frac{s_y \Delta f}{\vert y \vert^3}    {S_y f D_y f} \ K_y f \bar K_y f \ y.\nabla_{y}D_{y} f \nonumber \\
&&\times \cos(\frac{1}{2}(\arctan(\Delta_{y} f )+\arctan(\bar\Delta_{y} f ))) \nonumber \\ 
  &&\sin(\frac{1}{2}(\arctan(\Delta_{y} f )-\arctan(\bar\Delta_{y} f ))) \nonumber \\
 && \times      \int_{0}^{\infty} e^{-k} \sin(\frac{k}{2}S_y f) \sin(\frac{k}{2} D_y f)    \ dk  \ dy \ dx \nonumber \\
  \end{eqnarray}
  It not really difficult to observe that 
  
  \begin{eqnarray}
 \mathcal{S}_{1,2} &\lesssim& \int \vert \Delta f  \vert \int \frac{\vert s_y \Delta f \vert}{\vert y \vert^3} \vert y.\nabla_{y}D_{y} f \vert \left \vert  {S_y f D_y f} \ K_y f \bar K_y f  \right \vert    \ dy \ dx
  \end{eqnarray}
  From the inequality $\left\vert \frac{a^2-b^2}{(1+a^2)(1+b^2)} \right \vert\leq2$ valid for any $(a,b) \in \mathbb{R}^2$ one gets that 
  \begin{equation*} 
  \left \vert\frac{ S_{y}f D_y f}{(1+\Delta^2_{y}f)(1+\bar\Delta^2_{y}f)} \right \vert \leq2.
  \end{equation*} Which means that 
  \begin{equation} \label{borne2}
 \left \vert  {S_y f D_y f} \ K_y f \bar K_y f  \right \vert \leq 2.
  \end{equation} Therefore,
  \begin{equation} \label{type1}
\mathcal{S}_{1,2} \lesssim  \int \vert \Delta f  \vert \int \frac{\vert s_y \Delta f \vert}{\vert y \vert^3} \vert y.\nabla_{y}D_{y} f \vert \   dy \ dx
 \end{equation}
 Then, using equality \eqref{ddiff} and a classical scaling argument
  \begin{eqnarray*} 
\mathcal{S}_{1,2} &\lesssim&  \Vert \Delta f  \Vert_{L^{2}} \int_{0}^{1}\int \frac{\Vert s_y \Delta f \Vert_{L^{2}}}{\vert y \vert^3} \Vert s_{(r-1)y}  \nabla_{x}f \Vert_{L^{\infty}}  \ dy \ dr  \\
&& + \ \Vert \Delta f  \Vert_{L^{2}} \int \frac{\Vert s_y \Delta f \Vert_{L^{2}}}{\vert y \vert^3} \Vert \nabla_x s_{y}f  \Vert_{L^{\infty}}     \ dy \ dx 
 \end{eqnarray*} 

Hence,
 \begin{eqnarray*} 
 \mathcal{S}_{1,2} &\lesssim&\Vert \Delta f  \Vert_{L^{2}} \int_{0}^{1}\left(\int \frac{\Vert s_y \Delta f \Vert^{2}_{L^{2}}}{\vert y \vert^3} \right)^{1/2} \left(\int \frac{\Vert s_{(r-1)y}  \nabla_{x}f \ dy\Vert^{2}_{L^{\infty}}}{\vert y \vert^{3}} \ dy\right)^{1/2} \ dr   \\
 && + \ \Vert \Delta f  \Vert_{L^{2}} \left(\int \frac{\Vert s_y \Delta f \Vert^{2}_{L^{2}}}{\vert y \vert^{3}} \ dy \right)^{1/2} \left(\frac{\Vert \nabla_x s_{y}f  \Vert^{2}_{L^{\infty}}}{\vert y \vert^{3}}     \ dy \right)^{1/2}  
 \end{eqnarray*}

 So that one finally finds
 \begin{eqnarray*} 
 \mathcal{S}_{1,2}&\lesssim& \Vert f \Vert_{\dot H^{2}} \Vert \Delta f \Vert_{\dot H^{1/2}}\Vert  \nabla_{x}f \Vert_{\dot B^{1/2}_{\infty,2}} \\
 &\lesssim&\Vert  f \Vert^{2}_{\dot H^{5/2}} \Vert f \Vert_{\dot H^{2}} 
 \end{eqnarray*}

   \subsubsection{{\bf{Estimate of}} $\mathcal{S}_{1,3}$ }
We have
  \begin{eqnarray*}
 \mathcal{S}_{1,3}&=&-\frac{1}{8}\int \Delta f \int \frac{\Delta s_{y} f}{\vert y \vert^3} \left(K_y f + \bar K_{y} f\right)y.\nabla_{y}S_{y} f  \\&& \times \cos(\frac{1}{2}(\arctan(\Delta_{y} f )+\arctan(\bar\Delta_{y} f )))\sin(\frac{1}{2}(\arctan(\Delta_{y} f )-\arctan(\bar\Delta_{y} f )))     \\
  && \times    \int_{0}^{\infty} e^{-k}
    \sin(\frac{k}{2}S_y f) \sin(\frac{k}{2}D_y f) \ dk  \ dy \ dx  \\
  \end{eqnarray*}
  Using that,
 $$ \left \vert \int_{0}^{\infty} e^{-k} \sin(\frac{k}{2}S_y f) \sin(\frac{k}{2}D_y f)    \ dk \right \vert \leq \frac{\Gamma(2)}{2} \vert S_y f \vert, $$
 as well as the following bound
\begin{equation} \label{borne1}
 \vert K_{y}f+\bar K_{y}f \vert \leq 2
 \end{equation}
  one has
  \begin{eqnarray*}
   \mathcal{S}_{1,3}&\lesssim& \int \vert \Delta f  \vert \int \frac{ \vert s_y \Delta f \vert}{\vert y \vert^3} {\vert y.\nabla_{y}S_y f \vert}   \ dy \ dx.  \\
     \end{eqnarray*}
Now, we use the identity \eqref{som} that is 
$$y.\nabla_{y}S_{y} f=\frac{1}{\vert y \vert}s_y f(x)  + \frac{y}{\vert y\vert} \nabla_{x}\bar \delta_{y}f- \frac{y}{\vert y\vert}\nabla_{x}\delta_{y}f$$
So that,

  \begin{eqnarray} \label{type2}
   \mathcal{S}_{1,3}&\lesssim& \int \vert \Delta f  \vert \int \frac{ \vert s_y \Delta f \vert}{\vert y \vert^3} {\vert y.\nabla_{y}S_y f \vert}   \ dy \ dx  \\
   &\lesssim&\int \vert \Delta f  \vert \int \frac{ \vert s_y \Delta f \vert}{\vert y \vert^4} {\vert s_y f \vert}   \ dy \ dx \nonumber\\
   && + \int \vert \Delta f  \vert \int \frac{ \vert \Delta f \vert}{\vert y \vert^3} \left({\vert \nabla_{x} \bar\delta_y f \vert} + {\vert \nabla_{x}\delta_y f \vert}\right)   \ dy \ dx \nonumber\\
   &\lesssim& \Vert f  \Vert_{\dot H^{2}} \int \frac{ \Vert \Delta s_y f \Vert_{L^{2}}}{\vert y \vert^4} \Vert s_y f \Vert_{L^{\infty}}   \ dy   
   \nonumber\\ && + \ \Vert f  \Vert_{\dot H^{2} } \int \frac{ \Vert \Delta s_y f \Vert_{L^{2}}}{\vert y \vert^3} \left(\Vert \nabla_{x} \bar\delta_{y} f \Vert_{L^{\infty}} + \Vert \nabla_{x}\delta_y f \Vert_{L^{\infty}}\right) \ dy \nonumber\\
   &\lesssim&\Vert f  \Vert_{\dot H^{2}} \int \frac{ \Vert s_y \Delta f \Vert_{L^{2}}}{\vert y \vert^4} \Vert s_y f \Vert_{L^{\infty}}   \ dy \nonumber\\ 
    && + \ \Vert f  \Vert_{\dot H^{2} } \left(\int \frac{ \Vert \Delta s_y f \Vert^{2}_{L^{2}}}{\vert y \vert^3} \ dy \int \frac{\Vert \nabla_{x} \bar\delta_{y} f \Vert^{2}_{L^{\infty}}}{\vert y \vert^{3}} \ dy \right)^{1/2} \nonumber\\
   && + \ \Vert f  \Vert_{\dot H^{2} } \left(\int \frac{ \Vert \Delta s_y f \Vert^{2}_{L^{2}}}{\vert y \vert^3}\int \frac{\Vert \nabla_{x}\delta_{y} f \Vert^{2}_{L^{\infty}}}{\vert y \vert^{3}} \ dy \right)^{1/2} \nonumber\\
   &\lesssim&\Vert f  \Vert_{\dot H^{2}} \left(\int \frac{ \Vert s_y \Delta f \Vert^{2}_{L^{2}}}{\vert y \vert^3} \ dy \int \frac{\Vert s_y f \Vert^{2}_{L^{\infty}}}{\vert y \vert^{5}} \ dy\right)^{1/2}   \ dy   \nonumber\\
   && + \ \Vert f  \Vert_{\dot H^{2} } \left(\int \frac{ \Vert \Delta s_y f \Vert^{2}_{L^{2}}}{\vert y \vert^3} \ dy \int \frac{\Vert \nabla_{x} \bar\delta_{y} f \Vert^{2}_{L^{\infty}}}{\vert y \vert^{3}} \ dy \right)^{1/2} \nonumber\\
   &&  + \ \Vert f  \Vert_{\dot H^{2} } \left(\int \frac{ \Vert \Delta s_y f \Vert^{2}_{L^{2}}}{\vert y \vert^3}\int \frac{\Vert \nabla_{x}\delta_{y} f \Vert^{2}_{L^{\infty}}}{\vert y \vert^{3}} \ dy \right)^{1/2}  \nonumber\\
   &\lesssim& \Vert f \Vert_{\dot H^{2}} \left(\Vert \Delta f \Vert_{\dot H^{1/2}} \Vert f \Vert_{\dot B^{3/2}_{\infty,2}}+\Vert \Delta f \Vert_{\dot H^{1/2}}\Vert  \nabla_{x}f \Vert_{\dot B^{1/2}_{\infty,2}} \right) \nonumber
   \end{eqnarray}
   By using classical Besov embeddings, one finally finds that
   \begin{eqnarray*}
   \mathcal{S}_{1,3}\lesssim \Vert  f \Vert^{2}_{\dot H^{5/2}} \Vert f \Vert_{\dot H^{2}}
   \end{eqnarray*}

  \subsubsection{{\bf{Estimate of}} $\mathcal{S}_{1,4}$ }
We have using the identity \eqref{ddiff} together with the bound \eqref{borne2} that
 \begin{eqnarray*} 
\mathcal{S}_{1,4}&= &\frac{1}{8}\int \Delta f \int \frac{\Delta s_{y} f}{\vert y \vert^3}   S_y f D_y f K_y f \bar K_y f \ y.\nabla_{y}S_{y} f  \cos(\frac{1}{2}(\arctan(\Delta_{y} f )+\arctan(\bar\Delta_{y} f )))      \\
  &&     \sin(\frac{1}{2}(\arctan(\Delta_{y} f )-\arctan(\bar\Delta_{y} f ))) \int_{0}^{\infty} e^{-k}
    \sin(\frac{k}{2}S_y f) \sin(\frac{k}{2}D_y f)   \ dk  \ dy \ dx \nonumber \\
   &\lesssim&\int \vert \Delta f \vert \int  \frac{\vert s_y \Delta f \vert}{\vert y \vert^4}  \vert y.\nabla_{y}S_{y} f \vert    \left \vert{ S_{y}D_y f K_{y} f \bar K_{y}f}  \right\vert    \ dy \ dx \nonumber \\
   &\lesssim&\int \vert \Delta f \vert \int  \frac{\vert s_y \Delta f \vert}{\vert y \vert^3}  \vert y.\nabla_{y}S_{y} f  \vert     \ dy \ dx.
   \end{eqnarray*}
   One observes that this last estimate is exactly the same as \eqref{type2}
     \begin{eqnarray*} 
\mathcal{S}_{1,2} &\lesssim& \Vert  f \Vert^{2}_{\dot H^{5/2}} \Vert f \Vert_{\dot H^{2}}
 \end{eqnarray*}

\subsubsection{{\bf{Estimate of}} $\mathcal{S}_{1,5}$ }
We have
\begin{eqnarray*}
\mathcal{S}_{1,5}&=&-\frac{1}{8}\int \Delta f \int \frac{\Delta s_{y} f}{\vert y \vert^3} \left(K_y f + \bar K_{y} f\right)y.\nabla_{y}D_{y} f \sin(\frac{1}{2}(\arctan(\Delta_{y} f )+\arctan(\bar\Delta_{y} f )))\\     &&    \cos(\frac{1}{2}(\arctan(\Delta_{y} f )-\arctan(\bar\Delta_{y} f ))) \int_{0}^{\infty} e^{-k} \sin(\frac{k}{2}S_y f) \sin(\frac{k}{2}D_y f) \ dk  \ dy \ dx \\
  &\lesssim&\int \vert \Delta f \vert \int  \frac{\vert \Delta s_y f \vert}{\vert y \vert^2}  \vert \nabla_{y}D_{y} f \vert     \ dy \ dx.
  \end{eqnarray*}
   One notices that this last estimate is exactly the same as \eqref{type2}. Therefore, one directly infers that
  \begin{eqnarray*} 
\mathcal{S}_{1,5} \lesssim \Vert  f \Vert^2_{\dot H^{5/2}} \Vert f \Vert_{\dot H^{2}} 
\end{eqnarray*}

\subsubsection{{\bf{Estimate of}} $\mathcal{S}_{1,6}$ }
Recall that
\begin{eqnarray*}
\mathcal{S}_{1,6}&=&  \frac{1}{4}\int \Delta f \int  \frac{s_{y} \Delta f}{\vert y \vert^3} \sin(\frac{1}{2}(\arctan(\Delta_{y} f )+\arctan(\bar\Delta_{y} f )))       \sin(\frac{1}{2}(\arctan(\Delta_{y} f )-\arctan(\bar\Delta_{y} f ))) \\
 && y.\nabla_{y}S_y f  \int_{0}^{\infty} k e^{-k} \cos(\frac{k}{2}S_y f) \sin(\frac{k}{2} D_y f)    \ dk  \ dy \ dx \\
 &\lesssim& \int \vert \Delta f \vert \int \frac{\vert s_y \Delta f \vert}{\vert y \vert^{3}} \vert y.\nabla_{y} S_{y}f \vert  \ dy \ dx
\end{eqnarray*}
Then, it suffice to notice that this term may be estimated by means of \eqref{type1}, so that  
\begin{eqnarray*}
\mathcal{S}_{1,6} \lesssim  \Vert f\Vert^{2}_{\dot H^{5/2}} \Vert f\Vert_{\dot H^{2}} 
 \end{eqnarray*}
  
  \subsubsection{{\bf{Estimate of}} $\mathcal{S}_{1,7}$ }

  \begin{eqnarray*}
 \mathcal{S}_{1,7} &=&  \frac{1}{4}\int \Delta f \int \frac{s_y \Delta f}{\vert y \vert^3} \sin(\frac{1}{2}(\arctan(\Delta_{y} f )+\arctan(\bar\Delta_{y} f )))      \\
  && \sin(\frac{1}{2}(\arctan(\Delta_{y} f )-\arctan(\bar\Delta_{y} f ))) y.\nabla_{y}D_y f \int_{0}^{\infty} k e^{-k} \sin(\frac{k}{2} S_y f) \cos(\frac{k}{2} D_y f)   \ dk  \ dy \ dx  \\
  &\lesssim&\int \vert \Delta f \vert \int  \frac{\vert s_y \Delta f \vert}{\vert y \vert^2}  \vert y.\nabla_{y}D_{y} f  \vert      \ dy \ dx
  \end{eqnarray*}
  Then, the analysis done for \eqref{type1} allows one to get the same control as $S_{1,4}$, that is
  \begin{eqnarray*} 
\mathcal{S}_{1,7} &\lesssim& \Vert  f \Vert^2_{\dot H^{5/2}} \Vert f \Vert_{\dot H^{2}}
\end{eqnarray*}
    
Finally, collecting all the estimates we have obtained that
\begin{eqnarray} \label{t2}
\mathcal{S}_1 \lesssim \Vert  f \Vert^2_{\dot H^{5/2}} \Vert f \Vert_{\dot H^{2}}
\end{eqnarray}

{  {This ends the proof of Lemma \ref{s1}

\qed

}}

      \subsection{\bf{Estimate of $\mathcal{S}_2$}}
      
      {  {The term $\mathcal{S}_2$ will be decomposed into several terms which involve the second finite order differences. The goal will be to prove the following estimate.   
      
      \begin{lemma} \label{s2}
      The term $\mathcal{S}_{2}$ is controlled as follows
      \begin{eqnarray*}
      \mathcal{S}_{2} \lesssim \Vert f \Vert^{2}_{\dot H^{5/2}} (\Vert f \Vert_{\dot H^{2}}+\Vert f \Vert^2_{\dot H^{2}})
\end{eqnarray*}
      \end{lemma}
      
      \noindent{{\bf{Proof of Lemma}}} \ref{s2}
      }}
      Recall that 
      
      \begin{eqnarray*}
    \mathcal{S}_2  &=&\int \Delta f \int \nabla_{x}{\Delta\Delta_{y} f . \frac{y}{\vert y \vert^2}} \sin(\frac{1}{2}(\arctan(\Delta_{y} f )+\arctan(\bar\Delta_{y} f )))   \times \\
&& \sin(\frac{1}{2}(\arctan(\Delta_{y} f )-\arctan(\bar\Delta_{y} f ))) \int_{0}^{\infty} e^{-k} \cos(\frac{k}{2} S_y f) \cos(\frac{k}{2} D_y f)   \ dk  \ dy \ dx  \\
\end{eqnarray*}
{  {This term is too singular, we cannot estimate it directly. The idea is to try to balance the regularity in the space variable.  More precisely, we write that}}
 \begin{eqnarray*}
\mathcal{S}_2&=& \int \Delta f\int \nabla_{x} \Delta f . \frac{y}{\vert y \vert^3} \sin(\frac{1}{2}(\arctan(\Delta_{y} f )+\arctan(\bar\Delta_{y} f )))  \\
&& \times \sin(\frac{1}{2}(\arctan(\Delta_{y} f )-\arctan(\bar\Delta_{y} f ))) \int_{0}^{\infty} e^{-k} \cos(\frac{k}{2} S_y f) \cos(\frac{k}{2} D_y f)   \ dk  \ dy \ dx  \\
&-&\int \Delta f\int \nabla_{x}\Delta f(x-y) . \frac{y}{\vert y \vert^3} \sin(\frac{1}{2}(\arctan(\Delta_{y} f )+\arctan(\bar\Delta_{y} f )))   \times \\
&&\sin(\frac{1}{2}(\arctan(\Delta_{y} f )-\arctan(\bar\Delta_{y} f )))  \int_{0}^{\infty} e^{-k} \cos(\frac{k}{2} S_y f) \cos(\frac{k}{2} D_y f)   \ dk  \ dy \ dx.  \\
\end{eqnarray*}
Using that $\Delta f \nabla_{x} \Delta f= \frac{1}{2} \nabla_{x} (\Delta f)^2$, we may integrate by parts in $x$ and get that 
\begin{eqnarray*}
\mathcal{S}_2&=&-\frac{1}{2} \int \Delta f\int  \Delta f  \frac{y}{\vert y \vert^3} . \nabla_{x}\left(\sin(\frac{1}{2}(\arctan(\Delta_{y} f )+\arctan(\bar\Delta_{y} f ))) \right.   \\
&& \left. \sin(\frac{1}{2}(\arctan(\Delta_{y} f )-\arctan(\bar\Delta_{y} f ))) \right. \times  \left.\int_{0}^{\infty} e^{-k} \cos(\frac{k}{2} S_y f) \cos(\frac{k}{2} D_y f)  \right) \ dk  \ dy \ dx  \\
&-&\int \Delta f\int \nabla_{x}\Delta f(x-y) . \frac{y}{\vert y \vert^3} \sin(\frac{1}{2}(\arctan(\Delta_{y} f )+\arctan(\bar\Delta_{y} f )))   \times \\
&&\sin(\frac{1}{2}(\arctan(\Delta_{y} f )-\arctan(\bar\Delta_{y} f )))  \int_{0}^{\infty} e^{-k} \cos(\frac{k}{2} S_y f) \cos(\frac{k}{2} D_y f)   \ dk  \ dy \ dx  \\
\end{eqnarray*}
Then we come back to the more favorable finite difference, that is, we write that
\begin{eqnarray} \label{A2}
\mathcal{S}_2&=&-\frac{1}{2} \int \Delta f\int  \Delta  \delta_{y}f 
\frac{y}{\vert y \vert^3} .\nabla_{x}\left(\sin(\frac{1}{2}(\arctan(\Delta_{y} f )+\arctan(\bar\Delta_{y} f )))    \right. \times \nonumber \\
&& \left.\sin(\frac{1}{2}(\arctan(\Delta_{y} f )-\arctan(\bar\Delta_{y} f )))\int_{0}^{\infty} e^{-k} \cos(\frac{k}{2} S_y f) \cos(\frac{k}{2} D_y f)  \right) \ dk  \ dy \ dx \nonumber \\
&&-\frac{1}{2} \int \Delta f\int  \Delta f(x-y)  \frac{y}{\vert y \vert^3} .\nabla_{x}\left(\sin(\frac{1}{2}(\arctan(\Delta_{y} f )+\arctan(\bar\Delta_{y} f )))    \right. \times \nonumber \\
&& \left. \sin(\frac{1}{2}(\arctan(\Delta_{y} f )-\arctan(\bar\Delta_{y} f )))\int_{0}^{\infty} e^{-k} \cos(\frac{k}{2} S_y f) \cos(\frac{k}{2} D_y f) \right)  \ dk  \ dy \ dx \nonumber  \\
&&-\int \Delta f\int \nabla_{x}\Delta f(x-y) . \frac{y}{\vert y \vert^3} \sin(\frac{1}{2}(\arctan(\Delta_{y} f )+\arctan(\bar\Delta_{y} f )))   \times \nonumber \\
&&\sin(\frac{1}{2}(\arctan(\Delta_{y} f )-\arctan(\bar\Delta_{y} f )))  \int_{0}^{\infty} e^{-k} \cos(\frac{k}{2} S_y f) \cos(\frac{k}{2} D_y f)   \ dk  \ dy \ dx \nonumber  \\
&=& \mathcal{S}_{2,1} + \mathcal{S}_{2,2} + \mathcal{S}_{2,3}.   
\end{eqnarray}
      
    \subsubsection{{\bf{Estimate of $\mathcal{S}_{2,1}$}} } 
    
    In order to control $\mathcal{S}_{2,1}$, one observes that, by setting
 \begin{eqnarray} \label{not}
   T(f):=\nabla_{x}\left(\sin(\frac{1}{2}(\arctan(\Delta_{y} f )+\arctan(\bar\Delta_{y} f )))  \int_{0}^{\infty} e^{-k} \cos(\frac{k}{2} S_y f)    \right. \nonumber\\
 \left. \sin(\frac{1}{2}(\arctan(\Delta_{y} f )-\arctan(\bar\Delta_{y} f )))  \cos(\frac{k}{2} D_y f)   \ dk  \ dy \right),
\end{eqnarray} 

one may easily notice that, 

   \begin{eqnarray} \label{cas1}
 \vert T(f) \vert  \lesssim  \left\vert\frac{\nabla_{x}\left(f(x)-f(x\pm y)\right)}{\vert y \vert} \right \vert \left\vert  R(f) \right\vert,
        \end{eqnarray}

 where  the operator $R(f)$ is uniformly bounded by a fixed constant. \\
 
  Now, set $\delta^{\pm}_{y} f:=f(x)-f(x\pm y)$,  then 
  \begin{eqnarray*}
  \mathcal{S}_{2,1} &\lesssim& \Vert \Delta f \Vert_{L^{4}} \int \frac{\Vert \Delta \delta_{y} f \Vert_{L^{2}}}{\vert y \vert^{3/2}} \frac{\Vert \nabla \delta^{\pm}_{y} f \Vert_{L^{4}}}{\vert y \vert^{3/2}} \ dy  \\
  &\lesssim& \Vert f \Vert_{\dot H^{2}}  \Vert \Delta f \Vert_{\dot B^{1/2}_{2,2}} \Vert \nabla f \Vert_{\dot B^{1/2}_{4,2}} \\
  &\lesssim& \Vert f \Vert^{2}_{\dot H^{5/2}} \Vert f \Vert_{\dot H^{2}},  
   \end{eqnarray*} 
   
   where we used the  Sobolev embedding $\dot H^{1/2} \hookrightarrow L^{4}$ and the fact that $\dot H^1 \hookrightarrow \dot B^{1/2}_{4,2}$. 
    
 \subsubsection{{\bf{Estimate of}} $\mathcal{S}_{2,2}$ } 
 Recall that 
 \begin{eqnarray*} 
 \mathcal{S}_{2,2}&=&-\frac{1}{2} \int \Delta f\int  \Delta f(x-y)  \frac{y}{\vert y \vert^3} .\nabla_{x}\left(\sin(\frac{1}{2}(\arctan(\Delta_{y} f )+\arctan(\bar\Delta_{y} f )))    \right. \times \\
&& \left.  \hspace{-0.5cm}\sin(\frac{1}{2}(\arctan(\Delta_{y} f )-\arctan(\bar\Delta_{y} f )))\int_{0}^{\infty} e^{-k} \cos(\frac{k}{2} S_y f) \cos(\frac{k}{2} D_y f)   \right) \ dk  \ dy \ dx  \\
&=&-\frac{1}{2} \int \Delta f\int  \Delta f(x-y)  \frac{y}{\vert y \vert^3} . \sin(\frac{1}{2}(\arctan(\Delta_{y} f )+\arctan(\bar\Delta_{y} f ))) \times \\
&& \hspace{-1cm}\nabla_{x}\left(\sin(\frac{1}{2}(\arctan(\Delta_{y} f )-\arctan(\bar\Delta_{y} f )))\int_{0}^{\infty} e^{-k} \cos(\frac{k}{2} S_y f) \cos(\frac{k}{2} D_y f)  \right) \ dk  \ dy \ dx \\
&&-\frac{1}{2} \int \Delta f\int  \Delta f(x-y)  \frac{y}{\vert y \vert^3} .\nabla_{x}\left(\sin(\frac{1}{2}(\arctan(\Delta_{y} f )+\arctan(\bar\Delta_{y} f )))\right)  \\
&&  \sin(\frac{1}{2}(\arctan(\Delta_{y} f )-\arctan(\bar\Delta_{y} f )))\int_{0}^{\infty} e^{-k} \cos(\frac{k}{2} S_y f) \cos(\frac{k}{2} D_y f)   \ dk  \ dy \ dx \\
&=&\mathcal{S}_{2,2,1}+\mathcal{S}_{2,2,2}.
   \end{eqnarray*} 

\noindent $\bullet$ {\bf{Estimate of $\mathcal{S}_{2,2,1}$\\}}

 To estimate the term $\mathcal{S}_{2,2,1}$, it is not difficult to see that an estimate of the kind  \eqref{cas1} does not work anymore. One needs to find a slightly more refined  inequality. More precisely, we shall use the following Lemma.
 \begin{lemma} \label{cas2} The following inequality holds
 \begin{eqnarray} \label{cas22}
 \vert \mathcal{S}_{2,2,1} \vert  \lesssim \int \vert \Delta f \vert \int \frac{\vert \Delta f(x-y) \vert}{\vert y \vert^{2}} \left\vert\frac{\nabla_{x}\left(f(x)-f(x\pm y)\right)}{\vert y \vert} \right \vert \vert S_{y}f \vert \ dx \ dy.
        \end{eqnarray}
     \end{lemma}   
     \noindent {\bf{Proof of Lemma \ref{cas2}}}  Using twice the mean value theorem for instance, we have that $\left\vert\sin(\frac{1}{2}(\arctan(\Delta_{y} f )+\arctan(\bar\Delta_{y} f )))\right\vert \leq \vert S_{y} f \vert$ then if the derivative hits on one of the terms $$\sin(\frac{1}{2}(\arctan(\Delta_{y} f )-\arctan(\bar\Delta_{y} f )))\int_{0}^{\infty} e^{-k} \cos(\frac{k}{2} S_y f) \cos(\frac{k}{2} D_y f)   \ dk  \ dy,$$ it will be easily controlled by $\left\vert\frac{\nabla_{x}\left(f(x)-f(x\pm y)\right)}{\vert y \vert} \right \vert$ which proves that \eqref{cas22} holds. \\
\qed   

Using Lemma \ref{cas22} along with Sobolev embedding, we may estimate   $ \mathcal{S}_{2,2,1}$  as follows
\begin{eqnarray}  \label{s321}
 \vert \mathcal{S}_{2,2,1} \vert  &\lesssim& \Vert \Delta f \Vert_{L^{4}} \Vert \Delta f \Vert_{L^{2}} \int \frac{\Vert \nabla_{x}\left(f(x)-f(x\pm y)\right)\Vert_{L^{4}}}{\vert y \vert^{3/2}}  \frac{\Vert s_{y}f\Vert_{L^{\infty}}}{\vert y \vert^{5/2}} \ dy \nonumber\\
 &\lesssim& \Vert f \Vert_{\dot H^{5/2}} \Vert f \Vert_{\dot H^{2}}\Vert \nabla_{x}  f\Vert_{\dot B^{1/2}_{4,2}} \Vert f \Vert_{\dot B^{3/2}_{\infty,2}} \nonumber\\
 &\lesssim& \Vert f \Vert_{\dot H^{5/2}} \Vert f \Vert_{\dot H^{2}} \Vert \nabla_{x}  f\Vert_{\dot H^{1}} \Vert f \Vert_{\dot B^{3/2}_{\infty,2}} \nonumber\\
 &\lesssim& \Vert f \Vert^{2}_{\dot H^{5/2}} \Vert f \Vert^{2}_{\dot H^{2}} 
        \end{eqnarray}
        
\noindent $\bullet$ {\bf{Estimate of $\mathcal{S}_{2,2,2}$\\}}

 We now estimate the more delicate term $\mathcal{S}_{2,2,2}$, namely
      \begin{eqnarray*}
   \mathcal{S}_{2,2,2}&=&   -\frac{1}{2} \int \Delta f\int  \Delta f(x-y)  \frac{y}{\vert y \vert^3} .\nabla_{x}\left(\sin(\frac{1}{2}(\arctan(\Delta_{y} f )+\arctan(\bar\Delta_{y} f )))\right)  \\
&&  \sin(\frac{1}{2}(\arctan(\Delta_{y} f )-\arctan(\bar\Delta_{y} f )))\int_{0}^{\infty} e^{-k} \cos(\frac{k}{2} S_y f) \cos(\frac{k}{2} D_y f)   \ dk  \ dy \ dx. \\
      \end{eqnarray*}
      To do so, we shall use the fact that
       \begin{eqnarray*}
        \nabla_{x}(\sin(\frac{1}{2}(\arctan(\Delta_{y} f )+\arctan(\bar\Delta_{y} f ))))  &=& \frac{1}{2}\left(
        \frac{\nabla_{x}\Delta_{y} f}{1+\Delta^{2}_{y} f}+ \frac{\nabla_{x}\bar\Delta_{y} f}{1+\bar\Delta^{2}_{y} f} \right) \\
        &\times& \cos(\frac{1}{2}(\arctan(\Delta_{y} f )+\arctan(\bar\Delta_{y} f ))),
                  \end{eqnarray*}
     together with the fact that
      \begin{equation} \label{reste}
     \frac{\nabla_{x}\Delta_{y} f}{1+\Delta^{2}_{y} f}+ \frac{\nabla_{x}\bar\Delta_{y} f}{1+\bar\Delta^{2}_{y} f}=\underbrace{\frac{\nabla_{x}S_{y} f}{1+\Delta^{2}_{y} f}}_{delicate}+\underbrace{\nabla_{x} D_{y}f \frac{S_{y} f \ D_{y} f}{(1+\Delta^{2}_{y} f)((1+\bar\Delta^{2}_{y} f))}}_{easy}.
    \end{equation}
    Hence, this decomposition gives rise to two terms, that are
    
    \begin{eqnarray*}
   \mathcal{S}_{2,2,2}&=&   -\frac{1}{2} \int \Delta f\int  \Delta f(x-y)  \frac{y}{\vert y \vert^3} . \nabla_{x} D_{y}f \frac{S_{y} f \ D_{y} f}{(1+\Delta^{2}_{y} f)((1+\bar\Delta^{2}_{y} f))} \\
&&  \sin(\frac{1}{2}(\arctan(\Delta_{y} f )-\arctan(\bar\Delta_{y} f )))\int_{0}^{\infty} e^{-k} \cos(\frac{k}{2} S_y f) \cos(\frac{k}{2} D_y f)   \ dk  \ dy. \\
&-&\frac{1}{2} \int \Delta f\int  \Delta f(x-y)  \frac{y}{\vert y \vert^3} . \frac{\nabla_{x}S_{y} f}{1+\Delta^{2}_{y} f} \\
&&  \sin(\frac{1}{2}(\arctan(\Delta_{y} f )-\arctan(\bar\Delta_{y} f )))\int_{0}^{\infty} e^{-k} \cos(\frac{k}{2} S_y f) \cos(\frac{k}{2} D_y f)   \ dk  \ dy. \\
&=& \mathcal{S}_{2,2,2,1}+ \mathcal{S}_{2,2,2,2}
 \end{eqnarray*}

    The analysis of the first term of this last equality can be done by means of the Lemma \ref{cas2}. Indeed, since $\frac{ \vert D_{y} f \vert}{(1+\Delta^{2}_{y} f)((1+\bar\Delta^{2}_{y} f))}<1$  and since we have that $\left\vert\nabla_{x} D_{y}f \right \vert \lesssim \left\vert\frac{\nabla_{x}\left(f(x)-f(x\pm y)\right)}{\vert y \vert} \right \vert$. We find that it is estimated as $\mathcal{S}_{3,2,1},$ that is we have
  \begin{eqnarray*}
   \vert \mathcal{S}_{2,2,2,1} \vert  &\lesssim& \int \vert \Delta f \vert \int \frac{\vert \Delta f(x-y) \vert}{\vert y \vert^{2}} \left\vert\frac{\nabla_{x}\left(f(x)-f(x\pm y)\right)}{\vert y \vert} \right \vert \vert S_{y}f \vert \ dx \ dy \\
   &\lesssim& \Vert f \Vert^{2}_{\dot H^{5/2}} \Vert f \Vert^{2}_{\dot H^{2}} 
    \end{eqnarray*}
 The part involving the term $\frac{\nabla_{x}S_{y} f}{1+\Delta^{2}_{y} f}$ in equation \eqref{reste} is more delicate. The full term corresponds to $\mathcal{S}_{2,2,2,2}$. One shall us another strategy since there is an obvious lack of regularity. The idea is to try to balance the derivatives. Since the rational function in $\Delta_{y} f$ is not regular enough, one has to make appear oscillatory terms in order to avoid regulary issues. More precisely, we have the following Lemma.
\begin{lemma} \label{R} The term $\mathcal{S}_{2,2,2,2}$ may be rewritten as follows,
\begin{eqnarray*}
\mathcal{S}_{2,2,2,2}&=&\frac{1}{4} \int \Delta f\int  (\Delta f(x-y)-\Delta f(x+y))  \frac{y}{\vert y \vert^3} .{\nabla_{x}S_{y} f}  \\&&\times\sin(\frac{1}{2}(\arctan(\Delta_{y} f )-\arctan(\bar\Delta_{y} f )))  \int_{0}^{\infty} e^{-\gamma}\sin(\frac{\gamma}{2} S_y f)\\
&&  \sin(\frac{\gamma}{2} D_y f) \int_{0}^{\infty} e^{-k} \cos(\frac{k}{2} S_y f) \cos(\frac{k}{2} D_y f) \ d\gamma \ dk  \ dy \ dx \\
&&-\frac{1}{4} \int \Delta f\int   (\Delta f(x-y)+\Delta f(x+y))  \frac{y}{\vert y \vert^3} .{\nabla_{x}S_{y} f}  \\
&& \times \sin(\frac{1}{2}(\arctan(\Delta_{y} f )-\arctan(\bar\Delta_{y} f )))  \int_{0}^{\infty} e^{-\gamma} (\cos(\gamma \Delta_{y} f)+\cos(\gamma \bar\Delta_{y} f))\\
&&  \int_{0}^{\infty} e^{-k} \cos(\frac{k}{2} S_y f) \cos(\frac{k}{2} D_y f) \ d\gamma \ dk  \ dy \ dx \\
\end{eqnarray*}
\end{lemma}

\noindent {\bf{Proof of Lemma}} \ref{R}. First recall that 
\begin{eqnarray*}
\mathcal{S}_{2,2,2,2}&=&-\frac{1}{2} \int \Delta f\int  \Delta f(x-y)  \frac{y}{\vert y \vert^3} . \frac{\nabla_{x}S_{y} f}{1+\Delta^{2}_{y} f} \\
&&  \sin(\frac{1}{2}(\arctan(\Delta_{y} f )-\arctan(\bar\Delta_{y} f )))\int_{0}^{\infty} e^{-k} \cos(\frac{k}{2} S_y f) \cos(\frac{k}{2} D_y f)   \ dk  \ dy. 
\end{eqnarray*}
Then we symmetrize the non-oscillatory term, in other words, we write that 
\begin{eqnarray*}
\mathcal{S}_{2,2,2,2}&=&-\frac{1}{2} \int \Delta f\int  (\Delta f(x-y)-\Delta f(x+y))  \frac{y}{\vert y \vert^3} .{\nabla_{x}S_{y} f}  \\
&&\times \sin(\frac{1}{2}(\arctan(\Delta_{y} f )-\arctan(\bar\Delta_{y} f )))  \int_{0}^{\infty} e^{-\gamma} \cos(\gamma \Delta_{y} f)   \\
&& \int_{0}^{\infty} e^{-k} \cos(\frac{k}{2} S_y f) \cos(\frac{k}{2} D_y f) \ d\gamma \ dk  \ dy \ dx \\
&&-\frac{1}{2} \int \Delta f\int  \Delta f(x+y)  \frac{y}{\vert y \vert^3} .{\nabla_{x}S_{y} f}  \sin(\frac{1}{2}(\arctan(\Delta_{y} f )-\arctan(\bar\Delta_{y} f )))\\
&& \times \int_{0}^{\infty} e^{-\gamma} (\cos(\gamma \Delta_{y} f)+\cos(\gamma \bar\Delta_{y} f)) \\
&&\times\int_{0}^{\infty} e^{-k} \cos(\frac{k}{2} S_y f) \cos(\frac{k}{2} D_y f) \ d\gamma \ dk  \ dy \ dx \\
&&+\frac{1}{2} \int \Delta f\int  \Delta f(x+y)  \frac{y}{\vert y \vert^3} .{\nabla_{x}S_{y} f}    \sin(\frac{1}{2}(\arctan(\Delta_{y} f )-\arctan(\bar\Delta_{y} f ))) \\
&&\int_{0}^{\infty} e^{-\gamma} \cos(\gamma \bar\Delta_{y} f)\int_{0}^{\infty} e^{-k} \cos(\frac{k}{2} S_y f) \cos(\frac{k}{2} D_y f) \ d\gamma \ dk  \ dy \ dx \\
\end{eqnarray*}
By doing the change of variable $y \rightarrow -y$, one observes that the last term is equal to $-\mathcal{S}_{2,2,2,2}$ and that the two first terms may be symmetrized. More precisely, we find that
 
 \begin{eqnarray*}
\mathcal{S}_{2,2,2,2}&=&-\frac{1}{8} \int \Delta f\int  (\Delta f(x-y)-\Delta f(x+y))  \frac{y}{\vert y \vert^3} .{\nabla_{x}S_{y} f} \\
&&\hspace{-0.5cm} \times \sin(\frac{1}{2}(\arctan(\Delta_{y} f )-\arctan(\bar\Delta_{y} f ))) \int_{0}^{\infty} e^{-\gamma} (\cos(\gamma \Delta_{y} f)-\cos(\gamma \Delta_{y} f))  \\
&&  \int_{0}^{\infty} e^{-k} \cos(\frac{k}{2} S_y f) \cos(\frac{k}{2} D_y f) \ d\gamma \ dk  \ dy \ dx \\
&&-\frac{1}{4} \int \Delta f\int   (\Delta f(x-y)+\Delta f(x+y))  \frac{y}{\vert y \vert^3} .{\nabla_{x}S_{y} f}  \\
&&\hspace{-0.5cm} \times \sin(\frac{1}{2}(\arctan(\Delta_{y} f )-\arctan(\bar\Delta_{y} f )))  \int_{0}^{\infty} e^{-\gamma} (\cos(\gamma \Delta_{y} f)+\cos(\gamma \bar\Delta_{y} f))\\
&&  \int_{0}^{\infty} e^{-k} \cos(\frac{k}{2} S_y f) \cos(\frac{k}{2} D_y f) \ d\gamma \ dk  \ dy \ dx. \\
\end{eqnarray*}
This ends the proof of Lemma 8.4
\qed\\

Then, by using classical trigonometry formula and the fact that $$\Delta (f(x-y)-f(x+y))=-{\nabla_y .\nabla_{x} s_y f}$$ and 
\begin{eqnarray*}
{\Delta (f(x-y)+f(x+y))}&=&{\nabla_y .\nabla_x \left( f(x) -f(x-y)+f(x+y)-f(x)\right )}\\
&=&{\nabla_y .\nabla_x \left( \delta_y f - \bar \delta_y f \right )},
\end{eqnarray*} 
one may write  that, 
 \begin{eqnarray*}
\mathcal{S}_{2,2,2,2}&=&-\frac{1}{4} \int \Delta f\int (\nabla_y .\nabla_{x} s_y f)  \frac{y}{\vert y \vert^3} .{\nabla_{x}S_{y} f} \sin(\frac{1}{2}(\arctan(\Delta_{y} f )-\arctan(\bar\Delta_{y} f )))  \\
&&  \int_{0}^{\infty} e^{-\gamma} \sin(\frac{\gamma}{2} S_y f) \sin(\frac{\gamma}{2} D_y f) \int_{0}^{\infty} e^{-k} \cos(\frac{k}{2} S_y f) \cos(\frac{k}{2} D_y f) \ d\gamma \ dk  \ dy \ dx \\
&-&\frac{1}{2} \int \Delta f\int  \left(\nabla_y .\nabla_x d_{y} f  \right) \frac{y}{\vert y \vert^3} .{\nabla_{x}S_{y} f}  \sin(\frac{1}{2}(\arctan(\Delta_{y} f )-\arctan(\bar\Delta_{y} f )))\\
&&  \int_{0}^{\infty} e^{-\gamma} 
\cos(\frac{\gamma}{2} S_y f) \cos(\frac{\gamma}{2} D_y f) \int_{0}^{\infty} e^{-k} 
\cos(\frac{k}{2} S_y f) \cos(\frac{k}{2} D_y f) \ d\gamma \ dk  \ dy \ dx \\
\end{eqnarray*}
Hence, one finds
\begin{eqnarray*}
\mathcal{S}_{2,2,2,2}&=& -\frac{1}{4} \int \Delta f\int \left((\nabla_y .\nabla_{x} (s_y f+2d_y f\right)  \frac{y}{\vert y \vert^3} .{\nabla_{x}S_{y} f} \\&&\times\sin(\frac{1}{2}(\arctan(\Delta_{y} f )-\arctan(\bar\Delta_{y} f )))  \\
&&  \int_{0}^{\infty} e^{-\gamma} \sin(\frac{\gamma}{2} S_y f) \sin(\frac{\gamma}{2} D_y f) \int_{0}^{\infty} e^{-k} \cos(\frac{k}{2} S_y f) \cos(\frac{k}{2} D_y f) \ d\gamma \ dk  \ dy \ dx \\
\end{eqnarray*}

By integrating by parts (with respect to $y$), one finds

\begin{eqnarray*}
\mathcal{S}_{2,2,2,2}&=&\frac{1}{4} \int \Delta f\int  \nabla_x (s_{y}f+2d_y f) . \nabla_y\left(\frac{y}{\vert y \vert^3}\right) .{\nabla_{x}S_{y} f}\\
  && \times \int_{0}^{\infty} e^{-\gamma} \sin(\frac{\gamma}{2} S_y f) \sin(\frac{\gamma}{2} D_y f) \sin(\frac{1}{2}(\arctan(\Delta_{y} f )-\arctan(\bar\Delta_{y} f ))) \\
&&\times \int_{0}^{\infty} e^{-k} \cos(\frac{k}{2} S_y f) \cos(\frac{k}{2} D_y f)   \ d\gamma \ dk  \ dy \ dx \\
&+&\frac{1}{4} \int \Delta f\int  \nabla_x (s_{y}f+2d_y f) .\frac{y}{\vert y \vert^3}  \nabla_{x} .\nabla_y S_{y} f  \\
&&\times \int_{0}^{\infty} e^{-\gamma} \sin(\frac{\gamma}{2} S_y f) \sin(\frac{\gamma}{2} D_y f) \sin(\frac{1}{2}(\arctan(\Delta_{y} f )-\arctan(\bar\Delta_{y} f ))) \\
&&  \int_{0}^{\infty} e^{-k} \cos(\frac{k}{2} S_y f) \cos(\frac{k}{2} D_y f)   \ d\gamma \ dk  \ dy \ dx \\
&+& \frac{1}{8}\int \Delta f\int  \nabla_x (s_{y}f+2d_y f) . \frac{y}{\vert y \vert^3}  \nabla_{x} S_{y} f  \\
&& \times  \int_{0}^{\infty} \gamma e^{-\gamma} \nabla_{y}S_{y}f \ \cos(\frac{\gamma}{2} S_y f) \sin(\frac{\gamma}{2} D_y f)   \sin(\frac{1}{2}(\arctan(\Delta_{y} f )-\arctan(\bar\Delta_{y} f ))) \\&\times&\int_{0}^{\infty} e^{-k} \cos(\frac{k}{2} S_y f) \cos(\frac{k}{2} D_y f)   \ d\gamma \ dk  \ dy \ dx \\
 &+& \frac{1}{8}\int \Delta f\int  \nabla_x (s_{y}f+2d_y f). \frac{y}{\vert y \vert^3}  \nabla_{x} S_{y} f    \int_{0}^{\infty} \gamma e^{-\gamma} \nabla_{y}D_{y}f \\
 && \times \sin(\frac{\gamma}{2} S_y f) \cos(\frac{\gamma}{2} D_y f)\sin(\frac{1}{2}(\arctan(\Delta_{y} f )-\arctan(\bar\Delta_{y} f ))) \\
&&  \int_{0}^{\infty} e^{-k} \cos(\frac{k}{2} S_y f) \cos(\frac{k}{2} D_y f)   \ d\gamma \ dk  \ dy \ dx \\
&-&\frac{1}{8} \int \Delta f\int  \nabla_x (s_{y}f+2d_y f) .\frac{y}{\vert y \vert^3} \nabla_{x} S_{y} f \int_{0}^{\infty} e^{-\gamma} \sin(\frac{\gamma}{2} S_y f) \sin(\frac{\gamma}{2} D_y f)  \\
&& \nabla_{y}(\arctan(\Delta_{y} f )-\arctan(\bar\Delta_{y} f ))\cos(\frac{1}{2}(\arctan(\Delta_{y} f )-\arctan(\bar\Delta_{y} f ))) \\
&&\int_{0}^{\infty} e^{-k} \cos(\frac{k}{2} S_y f) \cos(\frac{k}{2} D_y f)   \ d\gamma \ dk  \ dy \ dx \\
&-& \frac{1}{8}\int \Delta f\int  \nabla_x (s_{y}f+2d_y f) .\frac{y}{\vert y \vert^3}  \nabla_{x} S_{y} f . \nabla_{y}S_{y}f \\
&&\times   \int_{0}^{\infty} e^{-\gamma}  \ \sin(\frac{\gamma}{2} S_y f) \sin(\frac{\gamma}{2} D_y f)   \sin(\frac{1}{2}(\arctan(\Delta_{y} f )-\arctan(\bar\Delta_{y} f )))\\
&&\times \int_{0}^{\infty} ke^{-k} \sin(\frac{k}{2} S_y f) \cos(\frac{k}{2} D_y f)   \ d\gamma \ dk  \ dy \ dx \\
&-& \frac{1}{8}\int \Delta f\int  \nabla_x (s_{y}f+2d_y f) . \frac{y}{\vert y \vert^3}  \nabla_{x} S_{y} f.\nabla_{y}D_{y}f    \\
&& \times \int_{0}^{\infty} e^{-\gamma}  \ \sin(\frac{\gamma}{2} S_y f) \sin(\frac{\gamma}{2} D_y f) \sin(\frac{1}{2}(\arctan(\Delta_{y} f )-\arctan(\bar\Delta_{y} f )))\\
&& \int_{0}^{\infty} ke^{-k} \cos(\frac{k}{2} S_y f) \sin(\frac{k}{2} D_y f)   \ d\gamma \ dk  \ dy \ dx \\
&=& \sum_{i}^{7} \mathcal{S}_{2,2,2,2,i}
\end{eqnarray*}
One can now start estimating $\mathcal{S}_{2,2,2,2,i}$, $i=1,...,7$. \\

$\bullet$ \ {\it{Estimate of}} $\mathcal{S}_{2,2,2,2,1}$ \\

Using Lemma \ref{facile} together with the fact that $\dot H^{5/2} \hookrightarrow \dot B^{2}_{4,2}$, one finds
\begin{eqnarray} \label{p4}
\mathcal{S}_{2,2,2,2,1}&=&\frac{1}{4} \int \Delta f\int  \nabla_x (s_{y}f+2d_y f) . \nabla_y\left(\frac{y}{\vert y \vert^3}\right) .{\nabla_{x}S_{y} f}
  \int_{0}^{\infty} e^{-\gamma} \sin(\frac{\gamma}{2} S_y f) \sin(\frac{\gamma}{2} D_y f) \nonumber\\
&&  \sin(\frac{1}{2}(\arctan(\Delta_{y} f )-\arctan(\bar\Delta_{y} f )))\int_{0}^{\infty} e^{-k} \cos(\frac{k}{2} S_y f) \cos(\frac{k}{2} D_y f)   \ d\gamma \ dk  \ dy \ dx \nonumber\\
&\lesssim& \Vert \Delta f \Vert_{L^{2}} \int \frac{\Vert \nabla_x \delta^{\pm}_{y}f \Vert_{L^{\infty}}}{\vert y \vert^{3/2}} \frac{\Vert \nabla_x s_{y}f \Vert_{L^{2}}}{\vert y \vert^{5/2}} \nonumber \ dy \\
&\lesssim& \Vert f \Vert_{\dot H^{2}} \Vert f \Vert_{\dot B^{3/2}_{\infty,2}} \Vert f \Vert_{\dot H^{5/2}} \nonumber \\
&\lesssim& \Vert f \Vert^{2}_{\dot H^{5/2}} \Vert f \Vert_{\dot H^{2}}\\ \nonumber
\end{eqnarray}

$\bullet$ \ {\it{Estimate of}} $\mathcal{S}_{2,2,2,2,2}$ \\

Using identity \eqref{dsy}, one finds
 
 \begin{eqnarray*}
\mathcal{S}_{2,2,2,2,2}&=& \frac{1}{4} \int \Delta f\int  \nabla_x (s_{y}f+2d_y f) .\frac{y}{\vert y \vert^3}  \nabla_{x} .\nabla_y S_{y} f 
 \\
 && \times \int_{0}^{\infty} e^{-\gamma} \sin(\frac{\gamma}{2} S_y f) \sin(\frac{\gamma}{2} D_y f)   \sin(\frac{1}{2}(\arctan(\Delta_{y} f )-\arctan(\bar\Delta_{y} f )))\\
 &&\times \int_{0}^{\infty} e^{-k} \cos(\frac{k}{2} S_y f) \cos(\frac{k}{2} D_y f)   \ d\gamma \ dk  \ dy \ dx \\
  \end{eqnarray*}
 Then, we write that
  \begin{eqnarray*}
\mathcal{S}_{2,2,2,2,2} &=& \frac{1}{4} \int \Delta f\int  \nabla_x (s_{y}f+2d_y f) .\frac{y}{\vert y \vert^3}  \nabla_{x} . \left(\nabla_{y}\left(\frac{1}{\vert y \vert}\right) s_y f\right) \\
 &\times&\int_{0}^{\infty} e^{-\gamma} \sin(\frac{\gamma}{2} S_y f)\sin(\frac{\gamma}{2} D_y f)   \sin(\frac{1}{2}(\arctan(\Delta_{y} f )-\arctan(\bar\Delta_{y} f )))\\
 &&\times\int_{0}^{\infty} e^{-k} \cos(\frac{k}{2} S_y f) \cos(\frac{k}{2} D_y f)   \ d\gamma \ dk  \ dy \ dx \\
&-& \frac{1}{4} \int \Delta f\int  \nabla_x (s_{y}f+2d_y f) .\frac{y}{\vert y \vert^3}  \nabla_{x} . \left(\nabla_{x} D_{y}f\right) \\
&&\times \int_{0}^{\infty} e^{-\gamma} \sin(\frac{\gamma}{2} S_y f) \sin(\frac{\gamma}{2} D_y f)  \sin(\frac{1}{2}(\arctan(\Delta_{y} f )-\arctan(\bar\Delta_{y} f ))) \\
&&\times \int_{0}^{\infty} e^{-k} \cos(\frac{k}{2} S_y f) \cos(\frac{k}{2} D_y f)   \ d\gamma \ dk  \ dy \ dx \\
&=& \mathcal{S}_{2,2,2,2,2,1}+ \mathcal{S}_{2,2,2,2,2,2,2}
\end{eqnarray*}

One observes that the estimate of $\mathcal{S}_{2,2,2,2,2,1}$ is similar to $\mathcal{S}_{2,2,2,2,1}$ (see \eqref{p4}). Indeed, we have that
 
$$
\mathcal{S}_{2,2,2,2,2,1} \lesssim  \Vert \Delta f \Vert_{L^{2}} \int \frac{\Vert \nabla_x \delta^{\pm}_{y}f \Vert_{L^{\infty}}}{\vert y \vert^{3/2}} \frac{\Vert \nabla_x s_{y}f \Vert_{L^{2}}}{\vert y \vert^{5/2}} \nonumber \ dy 
$$
hence,
\begin{equation} \label{e1}
\mathcal{S}_{2,2,2,2,2,1} \lesssim  \Vert f \Vert^{2}_{\dot H^{5/2}} \Vert f \Vert_{\dot H^{2}}
\end{equation}
As for $\mathcal{S}_{2,2,2,2,2,2}$, using Sobolev embedding and that $\dot H^{1} \hookrightarrow B^{1/2}_{4,2}$, we find

\begin{eqnarray} \label{e2}
\mathcal{S}_{2,2,2,2,2,2} &\lesssim& \Vert \Delta f \Vert_{L^{4}}  \int \frac{\Vert \nabla_x \delta^{\pm}_{y}f \Vert_{L^{4}}}{\vert y \vert^{3/2}} \frac{\Vert \Delta\delta_y f \Vert_{L^{2}}}{\vert y \vert^{3/2}} \ dy \nonumber\\
&\lesssim& \Vert \Delta f \Vert_{L^{4}} \Vert  f \Vert_{\dot B^{3/2}_{4,2}} \Vert f \Vert_{\dot H^{5/2}} \nonumber\\
&\lesssim& \Vert f \Vert^{2}_{\dot H^{5/2}} \Vert  f \Vert_{\dot H^{2}} 
\end{eqnarray}

Hence combining \eqref{e1} and \eqref{e2}, one finds

\begin{equation*} 
\mathcal{S}_{2,2,2,2,2} \lesssim  \Vert f \Vert^{2}_{\dot H^{5/2}} \Vert f \Vert_{\dot H^{2}}
\end{equation*}

\noindent $\bullet$ \ {\bf{Estimate of}} $\mathcal{S}_{2,2,2,2,3}$ \\

We split this term using identity \eqref{dsy}, we find that

\begin{eqnarray*}
\mathcal{S}_{2,2,2,2,3}&=& \frac{1}{8}\int \Delta f\int  \nabla_x (s_{y}f+2d_y f) . \frac{y}{\vert y \vert^3}  \nabla_{x} S_{y} f    \int_{0}^{\infty} \gamma e^{-\gamma} \nabla_{y}S_{y}f \\ 
&& \times\cos(\frac{\gamma}{2} S_y f) \sin(\frac{\gamma}{2} D_y f)   \sin(\frac{1}{2}(\arctan(\Delta_{y} f )-\arctan(\bar\Delta_{y} f ))) \\
&& \times \int_{0}^{\infty} e^{-k} \cos(\frac{k}{2} S_y f) \cos(\frac{k}{2} D_y f)   \ d\gamma \ dk  \ dy \ dx \\
\end{eqnarray*}
So that,
\begin{eqnarray*}
\mathcal{S}_{2,2,2,2,3}&=&\frac{1}{8}\int \Delta f\int  \nabla_x (s_{y}f+2d_y f) . \frac{y}{\vert y \vert^3}  \nabla_{x} S_{y} f . \nabla_{y} \left(\frac{1}{\vert y \vert}\right) \ s_y f \\
&& \times \int_{0}^{\infty} \gamma e^{-\gamma}   \ \cos(\frac{\gamma}{2} S_y f) \sin(\frac{\gamma}{2} D_y f)   \sin(\frac{1}{2}(\arctan(\Delta_{y} f )-\arctan(\bar\Delta_{y} f ))) \\
&&\int_{0}^{\infty} e^{-k} \cos(\frac{k}{2} S_y f) \cos(\frac{k}{2} D_y f)   \ d\gamma \ dk  \ dy \ dx \\
&-& \frac{1}{8}\int \Delta f\int  \nabla_x (s_{y}f+2d_y f) . \frac{y}{\vert y \vert^3}  \nabla_{x} S_{y} f \ \nabla_{x} D_{y}f    \int_{0}^{\infty} \gamma e^{-\gamma}  \ \cos(\frac{\gamma}{2} S_y f)  \\
&& \times \sin(\frac{\gamma}{2} D_y f)  \sin(\frac{1}{2}(\arctan(\Delta_{y} f )-\arctan(\bar\Delta_{y} f )))\\
&& \times \int_{0}^{\infty} e^{-k} \cos(\frac{k}{2} S_y f) \cos(\frac{k}{2} D_y f)   \ d\gamma \ dk  \ dy \ dx \\
&=& \mathcal{S}_{2,2,2,2,3,1} + \mathcal{S}_{2,2,2,2,3,2}
\end{eqnarray*}

For term $\mathcal{S}_{2,2,2,2,3,1}$, it suffices  to use that $\dot H^{\eta+1/2} \hookrightarrow \dot B^{k}_{4,4}$ and for $\eta=3/2$ and $\eta=2$, hence 
\begin{eqnarray} \label{r131}
\mathcal{S}_{2,2,2,2,3,1}&=&\frac{1}{8}\int \Delta f\int  \nabla_x (s_{y}f+2d_y f) . \frac{y}{\vert y \vert^3}  \nabla_{x} S_{y} f  s_y f\nabla_{y} \left(\frac{1}{\vert y \vert}\right) \int_{0}^{\infty}\gamma \\
&&   \times \ e^{-\gamma}   \cos(\frac{\gamma}{2} S_y f) \sin(\frac{\gamma}{2} D_y f)    \sin(\frac{1}{2}(\arctan(\Delta_{y} f )-\arctan(\bar\Delta_{y} f )))\nonumber\\
&& \times \int_{0}^{\infty} e^{-k} \cos(\frac{k}{2} S_y f) \cos(\frac{k}{2} D_y f)   \ d\gamma \ dk  \ dy \ dx \nonumber\\
&\lesssim& \Vert \Delta f \Vert_{L^{2}} \int \frac{\Vert \nabla_x\delta^{\pm}_{y} f\Vert_{L^{4}}}{\vert y \vert} \frac{\Vert \nabla_x s_{y} f\Vert_{L^{4}}}{\vert y \vert^{3/2}} \frac{\Vert s_{y} f \Vert_{L^{\infty}}}{\vert y \vert^{5/2}} \ dy \\
&\lesssim& \Vert \Delta f \Vert_{L^{2}} \Vert f \Vert_{\dot B^{3/2}_{\infty,2}} \left(\int \frac{\Vert \nabla_x\delta^{\pm}_{y} f\Vert^4_{L^{4}}}{\vert y \vert^4}  \ dy \int \frac{\Vert \nabla_x s_{y} f\Vert^4_{L^{4}}}{\vert y \vert^{6}}  \ dy\right)^{1/4} \nonumber \\
&\lesssim&\Vert f \Vert_{\dot H^{2}} \Vert f \Vert_{\dot B^{3/2}_{\infty,2}} \Vert  f  \Vert_{\dot B^{3/2}_{4,4}} \Vert  f  \Vert_{\dot B^{2}_{4,4}} \nonumber\\
&\lesssim&\Vert f \Vert^{2}_{\dot H^{5/2}} \Vert  f  \Vert^{2}_{\dot H^{2}} \nonumber
\end{eqnarray}

$\mathcal{S}_{2,2,2,2,3,2}$ is estimated as follows, using that $  \dot B^{5/3}_{6,3} \hookleftarrow \dot H^{7/3}= [\dot H^{5/2}, \dot H^{2}]_{\frac{2}{3},\frac{1}{3}}$, one finds that
\begin{eqnarray} \label{r132}
\mathcal{S}_{2,2,2,2,3,2}&=&- \frac{1}{8}\int \Delta f\int  \nabla_x (s_{y}f+2d_y f) . \frac{y}{\vert y \vert^3}  \nabla_{x} S_{y} f \ \nabla_{x} D_{y}f   \\
&\times& \int_{0}^{\infty} \gamma e^{-\gamma}  \ \cos(\frac{\gamma}{2} S_y f) \sin(\frac{\gamma}{2} D_y f)  \sin(\frac{1}{2}(\arctan(\Delta_{y} f )-\arctan(\bar\Delta_{y} f )))\nonumber \\
&& \int_{0}^{\infty} e^{-k} \cos(\frac{k}{2} S_y f) \cos(\frac{k}{2} D_y f)   \ d\gamma \ dk  \ dy \ dx \nonumber  \\
&\lesssim&  \Vert \Delta f \Vert_{L^{2}} \int \frac{\Vert \delta^{\pm}_{y}\nabla_x f \Vert^{3}_{L^{6}}}{\vert y \vert^{4}}  \ dy \\
&\lesssim&  \Vert  f \Vert_{\dot H^{2}} \Vert  f \Vert^{3}_{\dot B^{7/3}_{6,3}}  \nonumber \\
&\lesssim&\Vert f \Vert^{2}_{\dot H^{5/2}} \Vert  f  \Vert^{2}_{\dot H^{2}} \nonumber
\end{eqnarray}

Hence, combining \eqref{r131} and \eqref{r132}, we find that
\begin{eqnarray*}
\mathcal{S}_{2,2,2,2,3} \lesssim \Vert f \Vert^{2}_{\dot H^{5/2}} \Vert  f  \Vert^{2}_{\dot H^{2}}
\end{eqnarray*}

$\bullet$ \ Estimate of $\mathcal{S}_{2,2,2,2,4}$ \\

Using the identity \eqref{ddiff}, we may decompose  $\mathcal{S}_{2,2,2,2,4}$ as follows
\begin{eqnarray} \label{r14}
\mathcal{S}_{2,2,2,2,4}&=& \frac{1}{8}\int \Delta f\int  \nabla_x (s_{y}f+2d_y f) . \frac{y}{\vert y \vert^3}  \nabla_{x} S_{y} f    \int_{0}^{\infty} \gamma e^{-\gamma}  \\
&&\nabla_{y}D_{y}f \ \sin(\frac{\gamma}{2} S_y f) \cos(\frac{\gamma}{2} D_y f)\sin(\frac{1}{2}(\arctan(\Delta_{y} f )-\arctan(\bar\Delta_{y} f ))) \nonumber\\
&&  \int_{0}^{\infty} e^{-k} \cos(\frac{k}{2} S_y f) \cos(\frac{k}{2} D_y f)   \ d\gamma \ dk  \ dy \ dx \nonumber\\
&=&\frac{1}{8}\int \Delta f\int  \nabla_x (s_{y}f+2d_y f)  \frac{1}{\vert y \vert^3}  \nabla_{x} S_{y} f  \frac{1}{\vert y \vert} \\
&& \sin(\frac{1}{2}(\arctan(\Delta_{y} f )-\arctan(\bar\Delta_{y} f ))) \nonumber\\
&&  \int_{0}^{\infty} e^{-k} \cos(\frac{k}{2} S_y f) \cos(\frac{k}{2} D_y f)   \ dr \ d\gamma \ dk  \ dy \ dx \nonumber\\
&-&\frac{1}{8}\int \Delta f\int  \nabla_x (s_{y}f+2d_y f) . \frac{1}{\vert y \vert^3}  \nabla_{x} S_{y} f  \frac{y}{\vert y \vert} .{\nabla_x s_{y}f} \\
&&  \int_{0}^{\infty} \gamma e^{-\gamma} \ \sin(\frac{\gamma}{2} S_y f) \cos(\frac{\gamma}{2} D_y f)\sin(\frac{1}{2}(\arctan(\Delta_{y} f )-\arctan(\bar\Delta_{y} f ))) \nonumber\\
&&  \int_{0}^{\infty} e^{-k} \cos(\frac{k}{2} S_y f) \cos(\frac{k}{2} D_y f)   \ d\gamma \ dk  \ dy \ dx \nonumber \\
&=&\mathcal{S}_{2,2,2,2,4,1}+\mathcal{S}_{2,2,2,2,4,2}
\end{eqnarray}
In order to estimate $\mathcal{S}_{2,2,2,2,4,1}$ one uses an easy scaling argument for the integral in $r$, so that
\begin{eqnarray} \label{r141}
\mathcal{S}_{2,2,2,2,4,1}&=&\frac{1}{8}\int \Delta f\int  \nabla_x (s_{y}f+2d_y f)  \frac{1}{\vert y \vert^3}  \nabla_{x} S_{y} f  \frac{1}{\vert y \vert} \\
&&\int_{0}^{1}    \int_{0}^{\infty} y. s_{(r-1)y}  \nabla_{x}f\gamma e^{-\gamma}  \ \sin(\frac{\gamma}{2} S_y f) \cos(\frac{\gamma}{2} D_y f) \nonumber\\
&&  \sin(\frac{1}{2}(\arctan(\Delta_{y} f )-\arctan(\bar\Delta_{y} f ))) \nonumber \\
&&\int_{0}^{\infty} e^{-k} \cos(\frac{k}{2} S_y f) \cos(\frac{k}{2} D_y f)   \ dr \ d\gamma \ dk  \ dy \ dx \nonumber \\
&\lesssim&  \Vert  f \Vert_{\dot H^{2}}\int \frac{\Vert \delta^{\pm}_{y}\nabla_x f \Vert^{3}_{L^{6}}}{\vert y \vert^{4}}  \ dy \\
&\lesssim&  \Vert  f \Vert_{\dot H^{2}} \Vert  f \Vert^{3}_{\dot B^{7/3}_{6,3}}  \nonumber \\
&\lesssim&\Vert f \Vert^{2}_{\dot H^{5/2}} \Vert  f  \Vert^{2}_{\dot H^{2}} \nonumber
\end{eqnarray}
where we used again that $\dot B^{5/3}_{3,6} \hookleftarrow\dot H^{7/3}=[\dot H^{5/2}, \dot H^{2}]_{\frac{2}{3},\frac{1}{3}}.$ \\

The estimate of $\mathcal{S}_{2,2,2,2,4,2}$ is relatively easy, indeed, it suffices to observes that is it as regular as $\mathcal{S}_{2,2,2,2,4,1}$. More precisely, we have that
\begin{eqnarray*}
\mathcal{S}_{2,2,2,2,4,2} &\lesssim& \Vert \Delta f \Vert_{L^{2}} \int \frac{\Vert s_y \nabla_{x}f \Vert^{2}_{L^{4}}}{\vert y \vert^{5/2}}  \frac{ \Vert  s_y\nabla_{x}f \Vert_{\infty} }{\vert y \vert^{3/2}} \ dy \\  
&\lesssim&  \Vert  \Delta f \Vert_{L^{2}}\int \frac{\Vert \delta^{\pm}_{y}\nabla_x f \Vert^{3}_{L^{6}}}{\vert y \vert^{4}}  \ dy \\
&\lesssim&\Vert f \Vert^{2}_{\dot H^{5/2}} \Vert  f  \Vert^{2}_{\dot H^{2}} \nonumber
\end{eqnarray*}

$\bullet$ \ {\bf{Estimate of}} $\mathcal{S}_{2,2,2,2,5}$ \\

Recall that,

\begin{eqnarray*}
\mathcal{S}_{2,2,2,2,5}&=&-\frac{1}{4} \int \Delta f\int  \nabla_x (s_{y}f+2d_y f) .\frac{y}{\vert y \vert^3} \nabla_{x} S_{y} f \\
&&\int_{0}^{\infty} e^{-\gamma} \sin(\frac{\gamma}{2} S_y f) \sin(\frac{\gamma}{2} D_y f)  \nabla_{y}(\arctan(\Delta_{y} f )-\arctan(\bar\Delta_{y} f )) \\
&& \cos(\frac{1}{2}(\arctan(\Delta_{y} f )-\arctan(\bar\Delta_{y} f ))) \\
&&\int_{0}^{\infty} e^{-k} \cos(\frac{k}{2} S_y f) \cos(\frac{k}{2} D_y f)   \ d\gamma \ dk  \ dy \ dx \\
\end{eqnarray*}

Using formula \eqref{formuleA-}, we may decompose 
$\mathcal{S}_{2,2,2,2,5}$ as follows

\begin{eqnarray*}
\mathcal{S}_{2,2,2,2,5}&=&-\frac{1}{8} \int \Delta f\int  \nabla_x (s_{y}f+2d_y f) .\frac{y}{\vert y \vert^3} \nabla_{x} S_{y} f \int_{0}^{\infty} e^{-\gamma} \\
&&\sin(\frac{\gamma}{2} S_y f) \sin(\frac{\gamma}{2} D_y f)  S_y f D_y f K_y f \bar K_y f\nabla_{y}S_{y} f  \\
&& \cos(\frac{1}{2}(\arctan(\Delta_{y} f )-\arctan(\bar\Delta_{y} f )))\\ &&\int_{0}^{\infty} e^{-k} \cos(\frac{k}{2} S_y f) \cos(\frac{k}{2} D_y f)   \ d\gamma \ dk  \ dy \ dx \\
&& - \frac{1}{8} \int \Delta f\int  \nabla_x (s_{y}f+2d_y f) .\frac{y}{\vert y \vert^3} \nabla_{x} S_{y} f  \\
&&\int_{0}^{\infty} e^{-\gamma} \sin(\frac{\gamma}{2} S_y f) \sin(\frac{\gamma}{2} D_y f)  \left(K_{y}f+ \bar K_{y}f\right){\nabla_{y}D_{y} f} \\
&& \cos(\frac{1}{2}(\arctan(\Delta_{y} f )-\arctan(\bar\Delta_{y} f ))) \\
&&\int_{0}^{\infty} e^{-k} \cos(\frac{k}{2} S_y f) \cos(\frac{k}{2} D_y f)   \ d\gamma \ dk  \ dy \ dx. \\
&=&\mathcal{S}_{2,2,2,2,5,1}+\mathcal{S}_{2,2,2,2,5,2}
\end{eqnarray*}
By means of inequality \eqref{borne2}, one may write that
\begin{eqnarray}\label{p6}
\mathcal{S}_{2,2,2,2,5,1} \lesssim \int \vert \Delta f \vert \int 
\frac{\vert \delta^{\pm}_{y}\nabla_{x}f\vert}{\vert y \vert^{3/2}} \frac{\vert s_{y}\nabla_{x}f\vert}{\vert y \vert^{5/2}} \vert y. \nabla_{y} S_{y}f \vert \ dy \ dx.
\end{eqnarray}
Then, using the formula \eqref{som} one immediately finds that
\begin{eqnarray} \label{p5}
\mathcal{S}_{2,2,2,2,5,1} &\lesssim&\int \vert \Delta f \vert \int
 \frac{\vert \delta^{\pm}_{y}\nabla_{x}f\vert}{\vert y \vert^{3/2}} \frac{\vert s_{y}\nabla_{x}f\vert}{\vert y \vert^{5/2}}  \vert  s_y f(x) \vert \ dy \ dx \\
 && + \int \vert \Delta f \vert \int
  \frac{\vert \delta^{\pm}_{y}\nabla_{x}f\vert}{\vert y \vert^{3/2}} \frac{\vert s_{y}\nabla_{x}f\vert}{\vert y \vert^{5/2}} \vert \nabla_{x}\delta^{\pm}_{y}f \vert \ dy \ dx \nonumber\\
 &\lesssim& \Vert \Delta f \Vert_{L^{2}} \left(\int \frac{\Vert \nabla_x\delta^{\pm}_{y} f\Vert_{L^{4}}}{\vert y \vert} \frac{\Vert \nabla_x s_{y} f\Vert_{L^{4}}}{\vert y \vert^{3/2}} \frac{\Vert s_{y} f \Vert_{L^{\infty}}}{\vert y \vert^{5/2}} \ dy  \right. \nonumber \\
 && \hspace{3cm}+ \left. \int \frac{\Vert \delta^{\pm}_{y}\nabla_x f \Vert^{3}_{L^{6}}}{\vert y \vert^{4}}  \ dy \right) 
\end{eqnarray}
To control \eqref{p5}  one may follow the same steps as \eqref{r131} and \eqref{r141} and therefore

\begin{eqnarray} \label{p7}
\mathcal{S}_{2,2,2,2,5,1}&\lesssim&\Vert f \Vert^{2}_{\dot H^{5/2}} \Vert  f  \Vert^{2}_{\dot H^{2}}
\end{eqnarray}
As for $\mathcal{S}_{2,2,2,2,5,2}$, using that $\vert K_{y} f + \bar K_{y} f \vert\leq2$, we may write that
\begin{eqnarray*}
\mathcal{S}_{2,2,2,2,5,1} &\lesssim& \int \vert \Delta f \vert \int \frac{\vert s_{y}\nabla_{x}f\vert^{2}}{\vert y \vert^{4}} \vert y. \nabla_{y} D_{y}f \vert \ dy \ dx. \\
\end{eqnarray*}

Using formula \eqref{ddiff}, one finds

\begin{eqnarray*}
\mathcal{S}_{2,2,2,2,5,1} &\lesssim& \int \vert \Delta f \vert \int \frac{\vert \delta^{\pm}_{y}\nabla_{x}f\vert}{\vert y \vert}  \frac{\vert s_{y}\nabla_{x}f\vert}{\vert y \vert^{3}} \int_{0}^{1}\left( \vert   s_{(r-1)y}  \nabla_{x}f \vert + \vert \nabla_{x} s_{y}f \vert \right) dr \ dy \ dx \\
&\lesssim&  \Vert \Delta f \Vert_{L^{2}} \int \frac{\Vert \delta^{\pm}_y \nabla_{x}f \Vert^{3}_{L^{6}}}{\vert y \vert^{4}}  \ dy,
\end{eqnarray*}
where in the last inequality we used the same steps as \eqref{r141}
hence, we have

\begin{eqnarray*}
\mathcal{S}_{2,2,2,2,5} &\lesssim& \Vert f \Vert^{2}_{\dot H^{5/2}} \Vert  f  \Vert^{2}_{\dot H^{2}}
\end{eqnarray*}

$\bullet$ \ {\bf{Estimate of}} $\mathcal{S}_{2,2,2,2,6}$ \\

Using identity \eqref{som}, one finds that,

\begin{eqnarray*}
\mathcal{S}_{2,2,2,2,6}&=& - \frac{1}{8}\int \Delta f\int  \nabla_x (s_{y}f+2d_y f) .\frac{y}{\vert y \vert^3}  \nabla_{x} S_{y} f . \nabla_{y}S_{y}f   \\
&& \int_{0}^{\infty} e^{-\gamma}  \ \sin(\frac{\gamma}{2} S_y f) \sin(\frac{\gamma}{2} D_y f) \\
&&  \sin(\frac{1}{2}(\arctan(\Delta_{y} f )-\arctan(\bar\Delta_{y} f )))\\&&\int_{0}^{\infty} ke^{-k} \sin(\frac{k}{2} S_y f) \cos(\frac{k}{2} D_y f)   \ d\gamma \ dk  \ dy \ dx \\
&\lesssim& \int \vert \Delta f \vert \int  \frac{\vert \delta^{\pm}_{y}\nabla_{x}f\vert}{\vert y \vert}  \frac{\vert s_{y}\nabla_{x}f\vert}{\vert y \vert^{3}} \vert y. \nabla_{y} S_{y}f \vert \ dy \ dx. \\
\end{eqnarray*}
Therefore, following the same steps as  \eqref{p6}, hence we obtain the same control as \eqref{p7}, that is
\begin{eqnarray*}
\mathcal{S}_{2,2,2,2,6} &\lesssim&\Vert f \Vert^{2}_{\dot H^{5/2}} \Vert  f  \Vert^{2}_{\dot H^{2}}
\end{eqnarray*}
It remains to estimate the last term, that is,
\begin{eqnarray*}
\mathcal{S}_{2,2,2,2,7}&=& - \frac{1}{8}\int \Delta f\int  \nabla_x (s_{y}f+2d_y f) . \frac{y}{\vert y \vert^3}  \nabla_{x} S_{y} f.\nabla_{y}D_{y}f   \\
&& \int_{0}^{\infty} e^{-\gamma}  \ \sin(\frac{\gamma}{2} S_y f) \sin(\frac{\gamma}{2} D_y f) \sin(\frac{1}{2}(\arctan(\Delta_{y} f )-\arctan(\bar\Delta_{y} f ))) \\
&&  \int_{0}^{\infty} ke^{-k} \cos(\frac{k}{2} S_y f) \sin(\frac{k}{2} D_y f)   \ d\gamma \ dk  \ dy \ dx. \\
\end{eqnarray*}

Up to some bounded harmless terms, $\mathcal{S}_{2,2,2,2,7}$  is analogous to $\mathcal{S}_{2,2,2,2,4}$ (see \eqref{r14}) and therefore we may directly conclude that
\begin{eqnarray*}
\mathcal{S}_{2,2,2,2,7} \lesssim \Vert f \Vert^{2}_{\dot H^{5/2}} \Vert  f  \Vert^{2}_{\dot H^{2}}
\end{eqnarray*}
Finally, we have obtained that
\begin{eqnarray*}
\mathcal{S}_{2,2,2,2} \lesssim \Vert f \Vert^{2}_{\dot H^{5/2}} \Vert  f  \Vert^{2}_{\dot H^{2}}.
\end{eqnarray*}
Hence, combining all the previous estimates, we conclude that
\begin{eqnarray*}
\mathcal{S}_{2,2} &\lesssim& \Vert f \Vert^{2}_{\dot H^{5/2}} \Vert f \Vert_{\dot H^{2}}
\end{eqnarray*}

\subsubsection{{\bf{Estimate of}} $\mathcal{S}_{3,3}$ }
It remains to estimate  $\mathcal{S}_{3,3}$, we have
\begin{eqnarray*}
\mathcal{S}_{3,3}&=&-\int \Delta f\int \nabla_{x}\Delta f(x-y) . \frac{y}{\vert y \vert^3} \sin(\frac{1}{2}(\arctan(\Delta_{y} f )+\arctan(\bar\Delta_{y} f ))) \\ 
&& \sin(\frac{1}{2}(\arctan(\Delta_{y} f )-\arctan(\bar\Delta_{y} f ))) \int_{0}^{\infty} e^{-k} \cos(\frac{k}{2} S_y f) \cos(\frac{k}{2} D_y f)   \ dk  \ dy  \\
\end{eqnarray*}

 {  {Unlike $\mathcal{S}_{2,1}$ and $\mathcal{S}_{2,2}$ there are no derivative in $x$ in the oscillatory terms, so it cannot be treated in the same way as these terms. It is rather clear that the term  $\nabla_x\Delta f(x-y)$ is quite problematic. We would need a term of the kind $f(x-y)-f(x)=-\delta_y f$ in stead of $f(x-y)$. By using the fact that $\Delta \nabla_x f(x-y)=-\Delta \nabla_y  \delta_y f$, one may integrate by parts in $y$ and obtain a kind of regularization of this term. More precisely, we have that by integrating by parts in $y$}}
 \begin{eqnarray*}
\mathcal{S}_{2,3}&=&-\int \Delta f\int \Delta \delta_y f  \ \nabla_{y}. \left(\frac{y}{\vert y \vert^3} \sin(\frac{1}{2}(\arctan(\Delta_{y} f )+\arctan(\bar\Delta_{y} f ))) \right. \\ 
&& \left. \sin(\frac{1}{2}(\arctan(\Delta_{y} f )-\arctan(\bar\Delta_{y} f ))) \int_{0}^{\infty} e^{-k} \cos(\frac{k}{2} S_y f) \cos(\frac{k}{2} D_y f)   \ dk  \ dy \right)  
\end{eqnarray*}
This term may be controlled exactly the same way as $\mathcal{S}_{1}$ in \eqref{t1}. Indeed, the operator $s_y$ in $\Delta s_y f$ may be replaced by  $\Delta \delta_y f$. This is because of the fact that even if we would like to use the maximal regularity of $\Delta s_y f$ the operator $s_y$ would not be helpful. Recall that $\dot H^{5/2}$ is the maximale regularity one can afford. Hence, if we replace $\Delta s_y f$ by $\Delta \delta_y f$ it will give the same outcome. Moreover, the action of the differential operator $\nabla_{y}.$ when one integrates by parts will give rise to the same terms up to some harmless bounded functions (essentially trigonometric functions and Gamma functions evaluated in special values). Therefore, we have the same control as \eqref{t2} namely
\begin{eqnarray} \label{s33}
\mathcal{S}_{2,3} \lesssim \Vert f \Vert^{2}_{\dot H^{5/2}} \Vert f \Vert_{\dot H^{2}}
\end{eqnarray}

Hence, we have proved that 
\begin{eqnarray} 
\mathcal{S}_{2} \lesssim \Vert f \Vert^{2}_{\dot H^{5/2}} \left(\Vert f \Vert_{\dot H^{2}}+\Vert f \Vert^2_{\dot H^{2}} \right).
\end{eqnarray}

{  {Therefore, the proof of Lemma \ref{s2} is complete.
\qed
}}

\subsection{{\bf{Estimate of}} $\mathcal{S}_{3}$ }

The estimate of $\mathcal{S}_{3}$  is analogous to 
$\mathcal{S}_{2}$ that is the following Lemma holds.
  {  {The term $\mathcal{S}_2$ will be decomposed into several terms which involve the second finite order differences. The goal will be to prove the following estimate.   
      
      \begin{lemma} \label{s3}
      The term $\mathcal{S}_{3}$ is estimated as follows
      \begin{eqnarray*}
      \mathcal{S}_{3} \lesssim \Vert f \Vert^{2}_{\dot H^{5/2}} (\Vert f \Vert_{\dot H^{2}}+\Vert f \Vert^2_{\dot H^{2}})
\end{eqnarray*}
      \end{lemma}
      
      \noindent{{\bf{Proof of Lemma}}} \ref{s3}
      }}

 Indeed, recall that we have

\begin{eqnarray*}
\mathcal{S}_{3}&=&\int \nabla_{x} \Delta\bar\Delta_{y} f . \frac{y}{\vert y \vert^2} 
\cos(\frac{1}{2}(\arctan(\bar\Delta_{y} f + \arctan(\Delta_{y} f )) \times \\
&&\cos(\frac{1}{2}(\arctan(\bar\Delta_{y} f - \arctan(\Delta_{y} f ))\int_{0}^{\infty} e^{-k} \sin(\frac{k}{2} S_y f) \sin(\frac{k}{2} D_y f)  \ dk  \ dy. 
\end{eqnarray*}
If we do the change of variable $y \leftarrow -y$, then
\begin{eqnarray*}
\mathcal{S}_{3}&=&\int \nabla_{x} \Delta\Delta_{y} f . \frac{y}{\vert y \vert^2} 
\cos(\frac{1}{2}(\arctan(\bar\Delta_{y} f + \arctan(\Delta_{y} f )) \times \\
&&\cos(\frac{1}{2}(\arctan(\bar\Delta_{y} f - \arctan(\Delta_{y} f ))\int_{0}^{\infty} e^{-k} \sin(\frac{k}{2} S_y f) \sin(\frac{k}{2} D_y f)  \ dk  \ dy. 
\end{eqnarray*}
 Recall also that $\mathcal{S}_{3}$  is 
\begin{eqnarray*}
\mathcal{S}_{3}&=& \int \nabla_{x} \Delta\Delta_{y} f . \frac{y}{\vert y \vert^2} \sin(\frac{1}{2}(\arctan(\Delta_{y} f )+\arctan(\bar\Delta_{y} f )))    \times \\
&&\sin(\frac{1}{2}(\arctan(\Delta_{y} f )-\arctan(\bar\Delta_{y} f )))\int_{0}^{\infty} e^{-k} \cos(\frac{k}{2} S_y f) \cos(\frac{k}{2} D_y f)   \ dk  \ dy.  
\end{eqnarray*}

It is clear that they are equal up to interchanging the role of the sine and cosine functions. The role played by the oscillatory terms (that is all terms involving cosine and sine) in the estimate of $\mathcal{S}_{2}$ was not important since we finally estimated these terms by 1. Also, one notice that importantly, 
$\mathcal{S}_{3}$ and $\mathcal{S}_{2}$ have the same symmetry properties, that is, they are left invariant by the transformation $y \rightarrow -y$. Hence we may directly follow the same steps as  the control of $\mathcal{S}_{2}$ for the term $\mathcal{S}_{3}$. We deduce that,
\begin{eqnarray}
\mathcal{S}_{3} \lesssim \Vert f \Vert^{2}_{\dot H^{5/2}}\left(\Vert f \Vert_{\dot H^{2}}+\Vert f \Vert^{2}_{\dot H^{2}}\right).
\end{eqnarray}
{  {which is the desired estimated.

\qed}}

\subsection{Estimate of  $\mathcal{S}_{4}$}

{  {This term is fundamental in the sense that it contains the dissipation term. In order to extract the diffusive term, we have not only to linearize the oscillatory integrals but also to keep track of the directional derivative in the singular integral. Recall that,}}
\begin{eqnarray*}
\mathcal{S}_{4}&=& \frac{1}{2} \int \Delta f \int  \nabla_{x}\Delta D_{y} f . \frac{y}{\vert y \vert^2} \cos(\frac{1}{2}(\arctan(\Delta_{y} f )+\arctan(\bar\Delta_{y} f )))  \times  \\
&& \cos(\frac{1}{2}(\arctan(\Delta_{y} f )-\arctan(\bar\Delta_{y} f )))\int_{0}^{\infty} e^{-k} \cos(\frac{k}{2}(D_y f ))\cos(\frac{k}{2}(S_y f ))  \ dk  \ dy \ dx.\\
\end{eqnarray*}
{  {We are going to prove the following Lemma.
\begin{lemma} \label{s4} The term $\mathcal{S}_{4}$ is controlled as follows
\begin{eqnarray} \label{grat}
 \mathcal{S}_{4} \lesssim -\frac{1}{2}\frac{1}{(1+K(t)^{2})^{3/2}}\Vert f \Vert^2_{\dot H^{5/2}} +  \Vert f \Vert^{2}_{\dot H^{5/2}}\Vert f \Vert_{\dot H^{2}},
 \end{eqnarray}
 where, $K(t)=\displaystyle\sup_{x\in \mathbb R^2} \vert \nabla_x f  \vert_{L^\infty}(t)$.

\end{lemma}

\noindent {\bf{Proof of Lemma}} \ref{s4}

}}

In order to linearize we use the fact that $\cos(x)-1=-2\sin^{2}(x/2)$ twice, hence we may write
\begin{eqnarray*}
\mathcal{S}_{4}&=& \frac{1}{2} \int \Delta f \int  \nabla_{x}\Delta D_{y} f . \frac{y}{\vert y \vert^2} \cos(\frac{1}{2}(\arctan(\Delta_{y} f )+\arctan(\bar\Delta_{y} f )))   \times  \\
&&\cos(\frac{1}{2}(\arctan(\Delta_{y} f ))-(\arctan(\bar\Delta_{y} f )))\int_{0}^{\infty} e^{-k} (\cos(\frac{k}{2}(D_y f ))  \ dk  \ dy \ dx\\
&-&   \int \Delta f \int  \nabla_{x}\Delta D_{y} f . \frac{y}{\vert y \vert^2} \cos(\frac{1}{2}(\arctan(\Delta_{y} f )+\arctan(\bar\Delta_{y} f )))    \\
&&\cos(\frac{1}{2}(\arctan(\Delta_{y} f )-\arctan(\bar\Delta_{y} f )))\int_{0}^{\infty} e^{-k} \cos(\frac{k}{2}(D_y f )\sin^{2}(\frac{k}{4}(S_y f )  \ dk  \ dy \ dx\\
&=&- \int \Delta f \int  \nabla_{x}\Delta D_{y} f . \frac{y}{\vert y \vert^2} \sin^{2}(\frac{1}{4}(\arctan(\Delta_{y} f )+\arctan(\bar\Delta_{y} f )))     \\
&&\cos(\frac{1}{2}(\arctan(\Delta_{y} f )-\arctan(\bar\Delta_{y} f ))) \times \int_{0}^{\infty} e^{-k} (\cos(\frac{k}{2}(D_y f ))  \ dk  \ dy  \ dx\\
&-& \int \Delta f \int  \nabla_{x}\Delta D_{y} f . \frac{y}{\vert y \vert^2} \cos(\frac{1}{2}(\arctan(\Delta_{y} f )+\arctan(\bar\Delta_{y} f )))    \times  \\
&&\cos(\frac{1}{2}(\arctan(\Delta_{y} f )-\arctan(\bar\Delta_{y} f )))\int_{0}^{\infty} e^{-k} \cos(\frac{k}{2}(D_y f )\sin^{2}(\frac{k}{4}(S_y f )  \ dk  \ dy \ dx\\
&+& \frac{1}{2}\int \Delta f \int  \nabla_{x}\Delta D_{y} f . \frac{y}{\vert y \vert^2}   \cos(\frac{1}{2}(\arctan(\Delta_{y} f )-\arctan(\bar\Delta_{y} f ))) \int_{0}^{\infty} e^{-k} (\cos(\frac{k}{2}(D_y f ))  \ dk  \ dy \ dx\\
&=& \mathcal{S}_{4,1}+\mathcal{S}_{4,2}+\mathcal{S}_{4,3}.
\end{eqnarray*}

{  {In the sequel, we are going to estimate each of the $S_{4,i}$ and with a special attention on the term $\mathcal{S}_{4,3}$ which contains the elliptic component, the other terms being remainders. One the main difficulty in estimating the term  $\mathcal{S}_{4,3}$ is to have estimate of the singular integral which does not depend on the direction.}}

\subsubsection{{\bf{Estimate of}} $\mathcal{S}_{4,1}$}

 By integrating by parts, we find
\begin{eqnarray*}
\mathcal{S}_{4,1}&=& \int \Delta f \int \frac{\Delta s_{y} f }{\vert y \vert^3} \sin^{2}(\frac{1}{4}(\arctan(\Delta_{y} f )+\arctan(\bar\Delta_{y} f )))  \times  \\
&& \cos(\frac{1}{2}(\arctan(\Delta_{y} f )-\arctan(\bar\Delta_{y} f ))) \int_{0}^{\infty} e^{-k} (\cos(\frac{k}{2}(D_y f ))  \ dk  \ dy \ dx\\
&-&\int \Delta f \int \frac{\Delta s_{y} f }{\vert y \vert^3} y.\nabla_y(\arctan(\Delta_{y} f )+\arctan(\bar\Delta_{y} f )) \\ && \times\sin(\frac{1}{2}(\arctan(\Delta_{y} f )+\arctan(\bar\Delta_{y} f )))   \\
&&  \cos(\frac{1}{2}(\arctan(\Delta_{y} f )-\arctan(\bar\Delta_{y} f ))) \int_{0}^{\infty} e^{-k} (\cos(\frac{k}{2}(D_y f ))  \ dk  \ dy \ dx\\
&+& \int \Delta f \int \frac{\Delta s_{y} f }{\vert y \vert^3} \sin^{2}(\frac{1}{4}(\arctan(\Delta_{y} f )+\arctan(\bar\Delta_{y} f ))) \\
&& \times y.\nabla_y(\arctan(\Delta_{y} f )+\arctan(\bar\Delta_{y} f )) \times  \\
&&\sin((\arctan(\Delta_{y} f )-\arctan(\bar\Delta_{y} f ))) \int_{0}^{\infty} e^{-k} (\cos(\frac{k}{2}(D_y f ))  \ dk  \ dy \ dx\\
&+&\frac{1}{2}\int \Delta f \int \frac{\Delta s_{y} f }{\vert y \vert^3} \sin^{2}(\frac{1}{4}(\arctan(\Delta_{y} f )+\arctan(\bar\Delta_{y} f )))    
 \\
&& \times\cos(\frac{1}{2}(\arctan(\Delta_{y} f )-\arctan(\bar\Delta_{y} f ))) y.\nabla_y(D_y f)\\ 
&& \times \int_{0}^{\infty} ke^{-k} \sin(\frac{k}{2}(D_y f ))  \ dk  \ dy \ dx\\
&=&\mathcal{S}_{4,1,1}+\mathcal{S}_{4,1,2}+\mathcal{S}_{4,1,3}+\mathcal{S}_{4,1,4}.
\end{eqnarray*}

In order to estimate $\mathcal{S}_{4,1,1}$ we use the embedding $\dot H^{5/2} \hookrightarrow \dot B^{3/2}_{\infty,2}$, hence we get that

\begin{eqnarray*}
S_{4,1,1} &\lesssim& \Vert f \Vert_{\dot H^{2}} \int \frac{\Vert \Delta s_{y} f \Vert_{L^{2}}}{\vert y \vert^{3}} \frac{\Vert s_y f \Vert_{L^{\infty}}}{\vert y \vert} \ dy \\
&\lesssim&\Vert f \Vert_{\dot H^{2}} \Vert \Delta f\Vert_{\dot B^{1/2}_{2,2}}
   \Vert f \Vert_{\dot B^{3/2}_{\infty,2}} \\
   &\lesssim&\Vert f\Vert^2_{\dot H^{5/2}} \Vert f \Vert_{\dot H^{2}}
\end{eqnarray*}

The estimate of $\mathcal{S}_{4,1,2}$ is not difficult since, it suffices for instance to use the formula \eqref{formuleA+} we get that,
\begin{eqnarray*}
S_{4,1,2}&\lesssim& \int \vert \Delta f  \vert \int \frac{\vert s_y \Delta f \vert}{\vert y \vert^3} \vert y.\nabla_{y}D_{y} f \vert \left \vert  {S_y f D_y f} \ K_y f \bar K_y f  \right \vert    \ dy \ dx.
\end{eqnarray*}
Using the same step as \eqref{relou} we finally find that

      \begin{eqnarray*}
S_{4,1,2}\lesssim  \Vert f\Vert^2_{\dot H^{5/2}} \Vert f \Vert_{\dot H^{2}}
\end{eqnarray*}
It is not difficult to check that as well  for $i=2,3,4$ we have

  \begin{eqnarray*}
S_{4,1,i}\lesssim  \Vert f\Vert^2_{\dot H^{5/2}} \Vert f \Vert_{\dot H^{2}}
\end{eqnarray*}

\subsubsection{{\bf{Estimate of}} $\mathcal{S}_{4,2}$}

Recall that,

\begin{eqnarray*}
S_{4,2}&=&- \int \Delta f \int  \nabla_{x}\Delta D_{y} f . \frac{y}{\vert y \vert^2} \cos(\frac{1}{2}(\arctan(\Delta_{y} f )+\arctan(\bar\Delta_{y} f )))    \times  \\
&&\hspace{-0.5cm}\cos(\frac{1}{2}(\arctan(\Delta_{y} f )-\arctan(\bar\Delta_{y} f ))) \int_{0}^{\infty} e^{-k} \cos(\frac{k}{2}(D_y f )\sin^{2}(\frac{k}{4}(S_y f )  \ dk  \ dy \ dx\\
\end{eqnarray*}
So that by integration by parts, it is easy to estimate
\begin{eqnarray*}
S_{4,2}&\lesssim& \Vert f\Vert^2_{\dot H^{5/2}} \Vert f \Vert_{\dot H^{2}}
\end{eqnarray*}

\subsubsection{{\bf{Estimate of}} $\mathcal{S}_{4,3}$}
Recall  that 
   \begin{eqnarray} \label{dissi}
\mathcal{S}_{4,3}&=& \frac{1}{2} \int \Delta f \int  \nabla_{x}\Delta D_{y} f . \frac{y}{\vert y \vert^2}   \cos(\frac{1}{2}(\arctan(\Delta_{y} f )-\arctan(\bar\Delta_{y} f )))  \nonumber \\
&& \times \int_{0}^{\infty} e^{-k} (\cos(\frac{k}{2}(D_y f ))  \ dk  \ dy \ dx 
\end{eqnarray}
This term is absolutely fundamental since it plays a central role in the analysis of the Cauchy problem in the critical Sobolev space. Indeed, it contains the competition between the elliptic term and the diffusion. Of course, to see this competition one has to go through the term via the actions of "symmetrization" operators giving rise to sub-principal terms and the wanted ellipticity versus dissipative term. More precisely, one start by noticing that
\begin{eqnarray*}
\mathcal{S}_{4,3}&=&-\frac{1}{2}\int \Delta f  \int    \frac{\Delta\delta_{y} f}{\vert y \vert^3}   \left(\cos((\arctan(\frac{y}{\vert y \vert}. \nabla f(x)))-\cos((\arctan(\frac{y}{\vert y \vert}. \nabla f(x-y)))\right) \\
&& \times \int_{0}^{\infty} e^{-k} \cos({k} \frac{y}{\vert y \vert}. \nabla f(x))     \ dk  \ dy \ dx \\
&& - \frac{1}{2}\int \Delta f  \int    \frac{\Delta\delta_{y} f}{\vert y \vert^3}   \left(\cos((\arctan(\frac{y}{\vert y \vert}. \nabla f(x)))+\cos((\arctan(\frac{y}{\vert y \vert}. \nabla f(x-y)))\right) \\
&&\times \int_{0}^{\infty} e^{-k} \cos({k} \frac{y}{\vert y \vert}. \nabla f(x))    \ dk  \ dy \ dx \\
\end{eqnarray*}
Then, we write 
\begin{eqnarray*}
\mathcal{S}_{4,3}&=&-\frac{1}{2}\int \Delta f  \int    \frac{\Delta\delta_{y} f}{\vert y \vert^3}   \left(\cos((\arctan(\frac{y}{\vert y \vert}. \nabla f(x)))-\cos((\arctan(\frac{y}{\vert y \vert}. \nabla f(x-y)))\right) \\
&& \times \int_{0}^{\infty} e^{-k} \cos({k} \frac{y}{\vert y \vert}. \nabla f(x))      \ dk  \ dy \ dx \\
&-&\frac{1}{4}\int \Delta f  \int    \frac{\Delta\delta_{y} f}{\vert y \vert^3}   \left(\cos((\arctan(\frac{y}{\vert y \vert}. \nabla f(x)))+\cos((\arctan(\frac{y}{\vert y \vert}. \nabla f(x-y)))\right) \\
&& \int_{0}^{\infty} e^{-k} \left(\cos({k} \frac{y}{\vert y \vert}. \nabla f(x))-\cos({k} \frac{y}{\vert y \vert}. \nabla f(x-y)) \right)     \ dk  \ dy \ dx \\
&-&\frac{1}{4}\int \Delta f  \int    \frac{\Delta\delta_{y} f}{\vert y \vert^3}   \left(\cos((\arctan(\frac{y}{\vert y \vert}. \nabla f(x)))+\cos((\arctan(\frac{y}{\vert y \vert}. \nabla f(x-y)))\right) \\
&& \int_{0}^{\infty} e^{-k} \left(\cos({k} \frac{y}{\vert y \vert}. \nabla f(x))+\cos({k} \frac{y}{\vert y \vert}. \nabla f(x-y)) \right)      \ dk  \ dy \ dx \\
&=& \mathcal{S}_{4,3,1}+\mathcal{S}_{4,3,2}+\mathcal{S}_{4,3,3}
      \end{eqnarray*} 
      
We first remark that $\mathcal{S}_{4,3,i}$, for $i=1,2$ are easy to control. Indeed, it suffices to see that
\begin{eqnarray*}
\mathcal{S}_{4,3,i} &\lesssim& \Vert f \Vert_{\dot H^{2}} \int \frac{\Vert \delta_{y} \Delta f \Vert_{L^{2}}}{\vert y \vert^{3/2}}\frac{\Vert \delta_{y} \nabla f \Vert_{L^{\infty}}}{\vert y \vert^{3/2}} \ dy \\
&\lesssim& \Vert f \Vert_{\dot H^{2}} \Vert f \Vert_{\dot H^{5/2}} \Vert f \Vert_{\dot B^{3/2}_{\infty,2}} \\
&\lesssim&  \Vert f \Vert^2_{\dot H^{5/2}}\Vert f \Vert_{\dot H^{2}}
\end{eqnarray*}

As for $\mathcal{S}_{4,3,3}$, we need to extract the dissipation {\it{via}} several symmetrizations. More precisely, one writes that
      
      \begin{eqnarray*}
    \mathcal{S}_{4,3,3}  &=&-\frac{1}{8}\int   \int    \frac{\vert \Delta\delta_{y} f\vert^{2}}{\vert y \vert^3}   \left(\cos((\arctan(\frac{y}{\vert y \vert}. \nabla f(x)))+\cos((\arctan(\frac{y}{\vert y \vert}. \nabla f(x-y)))\right) \\
&& \int_{0}^{\infty} e^{-k} \left(\cos({k} \frac{y}{\vert y \vert}. \nabla f(x))+\cos({k} \frac{y}{\vert y \vert}. \nabla f(x-y)) \right)     \ dk  \ dy \ dx \\
&=&-\frac{1}{8}\int  \int    \frac{\vert \Delta\delta_{y} f\vert^{2}}{\vert y \vert^3} \left( \frac{1}{\sqrt{1+(\frac{y}{\vert y \vert}. \nabla f(x))^2}} + \frac{1}{\sqrt{1+(\frac{y}{\vert y \vert}. \nabla f(x-y))^2}} \right) \\
&& \times \left(\frac{1}{1+(\frac{y}{\vert y \vert}. \nabla f(x))^{2}} + \frac{1}{1+(\frac{y}{\vert y \vert}. \nabla f(x-y))^{2}} \right) \ dx \ dy \\
&=&\frac{1}{8}\int  \int    \frac{\vert \Delta\delta_{y} f\vert^{2}}{\vert y \vert^3} \\
&&\left(-\frac{1}{1+(\frac{y}{\vert y \vert}. \nabla f(x-y))^{2})^{3/2}} - \frac{1}{1+(\frac{y}{\vert y \vert}. \nabla f(x))^{2})^{3/2}} \right.\\
&& \left. \hspace{4cm}  -2 \frac{1}{\sqrt{1+(\frac{y}{\vert y \vert}. \nabla f(x-y))^2}} \frac{1}{1+(\frac{y}{\vert y \vert}. \nabla f(x))^{2}} \right) \ dk  \ dy  \\
&=&\frac{1}{8}\int  \int    \frac{\vert \Delta\delta_{y} f\vert^{2}}{\vert y \vert^3} \\
&&\times \left(-4+4-\frac{1}{1+(\frac{y}{\vert y \vert}. \nabla f(x-y))^{2})^{3/2}} - \frac{1}{1+(\frac{y}{\vert y \vert}. \nabla f(x))^{2})^{3/2}} \right. \\
&& \left. \hspace{4cm} -2 \frac{1}{\sqrt{1+(\frac{y}{\vert y \vert}. \nabla f(x-y))^2}} \frac{1}{1+(\frac{y}{\vert y \vert}. \nabla f(x))^{2}} \right) \ dy \ dx \\
&=& -\frac{1}{2} \Vert f \Vert^2_{\dot H^{5/2}} + \frac{1}{8}\int  \int    \frac{\vert \Delta\delta_{y} f\vert^{2}}{\vert y \vert^3} \\
&&\times \left(4-\frac{1}{1+(\frac{y}{\vert y \vert}. \nabla f(x-y))^{2})^{3/2}} - \frac{1}{1+(\frac{y}{\vert y \vert}. \nabla f(x))^{2})^{3/2}} \right. \\
&&\left. \hspace{4cm} -2 \frac{1}{\sqrt{1+(\frac{y}{\vert y \vert}. \nabla f(x-y))^2}} \frac{1}{1+(\frac{y}{\vert y \vert}. \nabla f(x))^{2}} \right) \ dy \ dx \\
&\leq&\underbrace{-\frac{1}{2} \Vert f \Vert^2_{\dot H^{5/2}} + \frac{1}{8} \Vert f \Vert^2_{\dot H^{5/2}} \left(4-4 \frac{1}{(1+K^{2})^{3/2}}\right)}_{=-\frac{1}{2} \Vert f \Vert^2_{\dot H^{5/2}} + \frac{1}{2} \Vert f \Vert^2_{\dot H^{5/2}} \left(1- \frac{1}{(1+K^{2})^{3/2}}\right)}
      \end{eqnarray*}
      Hence, one finally finds 
        \begin{eqnarray} \label{top}
\mathcal{S}_{4,3}\leq  -\frac{1}{2}\frac{1}{(1+K(t)^{2})^{3/2}}\Vert f \Vert^2_{\dot H^{5/2}} 
      \end{eqnarray}
      where, $K(t)=\displaystyle\sup_{x\in \mathbb R^2} \vert \nabla_x f  \vert_{L^\infty}(t)$.
      
      Hence, we gathering the estimates, we 
      \begin{eqnarray}
       \mathcal{S}_{4} \lesssim -\frac{1}{2}\frac{1}{(1+K(t)^{2})^{3/2}}\Vert f \Vert^2_{\dot H^{5/2}} +  \Vert f \Vert^{2}_{\dot H^{5/2}}\Vert f \Vert_{\dot H^{2}},
 \end{eqnarray}have proved that

      \begin{remark} It is crucial to note that estimate \eqref{top} above shows the parabolic character of the Muskat problem whenever the slope does not blow-up. Indeed, when $K(t)\rightarrow +\infty$ the regularizing  effect disappear (as it was also observed in the 2D case \cite{CL}).
      \end{remark}

\section{Sobolev energy inequality}

From the Section \ref{ti} (less singular terms) and from  Section \ref{sing} (most singular terms) we have proved respectively inequality \eqref{flop} and inequality \eqref{grat} (see in particular \eqref{top}). Hence, combining all these estimates lead to  
\begin{eqnarray}
\frac{1}{2} \partial_{t} \Vert f \Vert^{2}_{\dot H^{2}} + \frac{1}{2(1+K(t)^{2})^{3/2}} \Vert f \Vert^{2}_{\dot H^{5/2}}  \lesssim \Vert f \Vert^{2}_{\dot H^{5/2}}\left(\Vert f \Vert_{\dot H^{2}}+\Vert f \Vert^{2}_{\dot H^{2}}\right).
\end{eqnarray}
Integrating in time $s\in[0,T]$, and multiplying by 2 one finds
\begin{eqnarray*}
  \Vert f(.,T) \Vert^{2}_{\dot H^{2}} + \frac{1}{(1+K^{2})^{3/2}} \int_0^T \Vert f \Vert^{2}_{\dot H^{5/2}} \ ds   \lesssim  \Vert f_0 \Vert^{2}_{\dot H^{2}}  + \int_0^T \Vert f \Vert^{2}_{\dot H^{5/2}}\left(\Vert f \Vert_{\dot H^{2}}+\Vert f \Vert^{2}_{\dot H^{2}}\right) \ ds,
\end{eqnarray*}
where $K=\displaystyle \sup_{t\geq0} \sup_{x\in \mathbb R^2} \vert \nabla_x f(x,t)  \vert.$
{  {Finally, we have proved that
\begin{eqnarray} \label{grat}
\sum_{i=1}^{4} \mathcal{S}_{i} \lesssim -\frac{1}{2}\frac{1}{(1+K(t)^{2})^{3/2}}\Vert f \Vert^2_{\dot H^{5/2}} +  \Vert f \Vert^{2}_{\dot H^{5/2}}\left(\Vert f \Vert_{\dot H^{2}}+\Vert f \Vert^{2}_{\dot H^{2}}\right).
\end{eqnarray}
which ends the proof of Lemma \ref{s4}

\qed

Collecting all the estimates proved in Lemma \ref{s1}, \ref{s2}, \ref{s3} and \ref{s4} we have finally obtained that

\begin{eqnarray} \label{grat}
 \sum_{i=1}^{4}\mathcal{S}_{i} \lesssim -\frac{1}{2}\frac{1}{(1+K(t)^{2})^{3/2}}\Vert f \Vert^2_{\dot H^{5/2}} +  \Vert f \Vert^{2}_{\dot H^{5/2}}\Vert f \Vert_{\dot H^{2}},
 \end{eqnarray}
 
 }}

\section{\hspace{-0,08cm}Slope control and uniform bound using control of slope}

In this section we show how to control $\|\nabla f\|_{L^\infty}$ in terms of critical Sobolev norms. 
{  { More precisely, we have the following Lemma
\begin{lemma}
Let $f$ be a solution to the 3D Muskat equation with initial data $f_0 \in \dot H^{2} \cap \dot W^{1,\infty}$, then one has the following control of the Lipschitz semi-norm
\begin{equation}\label{bf} 
\Vert \nabla f \Vert^2_{L^\infty}(t)\leq \|\nabla f_0\|^2_{L^\infty}+\int_0^t\|f\|^2_{\dot{H}^{5/2}}(s)ds.
\end{equation} 
\end{lemma}
}}
Recall that the Muskat problem can be written as follows 
\begin{equation}\label{slopecontrol}
\partial_t f (t,x) =  P.V.\int \frac{\nabla f(x)\cdot y-(f(x)-f(x-y))}{\vert y \vert^3}  \frac{dy}{(1+\Delta^{2}_{y}f)^{3/2}}
\end{equation}
 
By taking one derivative in equation \eqref{slopecontrol} one finds
\begin{align*}
\partial_{j}f_{t} (x) =& \nabla\partial_{j}f(x)\cdot P.V.\int \frac{ y}{\vert y \vert^3}  \frac{dy}{(1+\Delta^{2}_{y}f(x))^{3/2}}- P.V.\int \frac{(\partial_j f(x)-\partial_j f(x-y))dy}{\vert y \vert^3(1+\Delta^{2}_{y}f(x))^{3/2}}\\
&-3 P.V.\int \frac{\nabla f(x)\cdot \frac{y}{|y|}-\Delta_{y}f(x)}{|y|}\Delta_{y}\partial_jf(x) \frac{\Delta_{y}f(x)}{(1+\Delta^{2}_{y}f(x))^{5/2}}\frac{dy}{|y|}.
\end{align*}
Set $M(t)=\displaystyle\sup_{x \in \mathbb R^2} \partial_{j}f(x,t)$. Since we are considering a regular solution {\it{e.g.}} $f(t,.) \in \mathcal{C}^2$, we have that $M(t)=\displaystyle\sup_{x \in \mathbb R^2} \partial_{j}f(x,t)= \partial_{j}f(x_t,t)$ and that $M'(t)=\partial_{j}f_t(x_t,t)$ are differentiable almost every time $t$ (thanks to Rademacher's theorem). By evaluating the above evolution equation at $x=x_t$ one finds that the first term on the right is zero and the second has a  sign. Omitting to write the $p.v$ for simplicity, we find\begin{eqnarray*}
M'(t)&\leq& -3 \int \frac{\nabla f(x_t)\cdot \frac{y}{|y|}-\Delta_{y}f(x_t)}{|y|}\Delta_{y}\partial_jf(x_t) \frac{\Delta_{y}f(x_t)}{(1+\Delta^{2}_{y}f(x_t))^{5/2}}\frac{dy}{|y|} \\
&\lesssim&  \int \frac{\Vert \nabla f(x_t). \frac{y}{\vert y\vert}-\Delta_{y}f(x_t) \Vert_{L^{\infty}}}{\vert y \vert} \frac{\Vert \partial_j \delta_{y}f(x_t) \Vert_{L^{\infty}}}{\vert y \vert^2} \ dy \\
&\lesssim& \left( \int \frac{\Vert \nabla f(x_t). \frac{y}{\vert y\vert}-\Delta_{y}f(x_t) \Vert^2_{L^{\infty}}}{\vert y \vert^3} \ dy \right)^{1/2} \left( \int \frac{\Vert \partial_j \delta_{y}f(x_t) \Vert^2_{L^{\infty}}}{\vert y \vert^3} \ dy \right)^{1/2}.
\end{eqnarray*}
Hence,
$$
M'(t)\lesssim \Vert \nabla f\Vert^2_{\dot{B}_{\infty,2}^{1/2}}\lesssim \Vert f \Vert^2_{\dot{H}^{5/2}}.
$$
Analogously, the same holds for the evolution of the minimum $m(t)$, so that by integrating in time 
\begin{equation}\label{bf}
\|\nabla f\|^2_{L^\infty}(t)\leq \|\nabla f_0\|^2_{L^\infty}+\int_0^t\|f\|^2_{\dot{H}^{5/2}}(s)ds.
\end{equation}

From the Sobolev energy inequality of the previous section, we have 
$$
\partial_{t}\Vert f \Vert^2_{\dot{H}^2}(t)+\frac{\Vert f \Vert ^2_{\dot{H}^{5/2}}}{(1+\|\nabla f\|_{L^\infty}^2)^{3/2}}\leq C\Vert f\Vert^2_{\dot{H}^{5/2}}\left(\Vert f \Vert_{\dot{H}^2}(t)+ \Vert f \Vert ^2_{\dot{H}^2}\right),
$$
where $C>0$ is a fixed constant. Since $(1+x^2+\mathcal{D}(t))^{-3/2}\leq(1+x^2)^{-3/2}$ for any  $\mathcal{D}(t)\geq0$, then from inequality \eqref{bf} we obtain that
\begin{equation}\label{global1iteration}
\partial_{t}\Vert f\Vert^2_{\dot{H}^2}+\frac{\Vert f\Vert^2_{\dot{H}^{5/2}}}{(1+\Vert \nabla f_0\Vert_{L^\infty}^2+\mathcal{D}(t))^{3/2}}\leq   C\Vert f\Vert^2_{\dot{H}^{5/2}}\left(\Vert f\Vert_{\dot{H}^2}(t)+\Vert f\Vert^2_{\dot{H}^2}(t)\right),
\end{equation}
where
$$
\mathcal{D}(t)=\int_0^t \Vert f \Vert^2_{\dot{H}^{5/2}}(s)ds,\quad \mbox{with}  \quad \mathcal{D}(0)=0.
$$
We consider the smallness conditions (to get control of the $L^{2}H^{5/2}$ semi-norm) for $ \Vert f_0 \Vert_{\dot{H}^2}$ given by
\begin{equation}\label{smallcondition1}
C(\Vert f_0 \Vert_{\dot{H}^2}+ \Vert f_0 \Vert^2_{\dot{H}^2})<\frac{1}{(2+ \Vert \nabla f_0 \Vert_{L^\infty}^2)^{3/2}}, 
\end{equation}
together with
\begin{equation}
\frac{\Vert f_0 \Vert^2_{\dot{H}^2}(2+\Vert\nabla f_0\Vert_{L^\infty}^2)^{3/2}}{1-C(\Vert f_0 \Vert_{\dot{H}^2}+ \Vert f_0 \Vert^2_{\dot{H}^2})(2+ \Vert \nabla f_0 \Vert_{L^\infty}^2)^{3/2}}<1.
\end{equation}
Therefore, after a short amount of time 
$$
\partial_t \Vert f \Vert^2_{\dot{H}^2}<0,\quad\mbox{together with}\quad \mathcal{D}(t)<1.
$$ 
By integrating in time,
$$
 \Vert f \Vert^2_{\dot{H}^2}(t)+\Big(\frac{1}{(2+ \Vert \nabla f_0 \Vert_{L^\infty}^2)^{3/2}}- C\left( \Vert f_0 \Vert_{\dot{H}^2}+ \Vert f_0 \Vert^2_{\dot{H}^2}\right)\Big)\mathcal{D}(t)\leq \Vert f_0 \Vert^2_{\dot{H}^2},
$$
so that bootstrapping the argument, we are able to find 
above identity for all time $t>0$, so that
$$
\|f\|^2_{\dot{H}^2}(t)\leq \|f_0\|^2_{\dot{H}^2},\quad\mbox{together with}\quad \mathcal{D}(t)<1.
$$
Assuming that there exists a first time $t^*$ such that $\mathcal{D}(t^*)=1,$ gives a contradiction. 
\qed
{  {
\section{Uniqueness}

We are going to prove the following Lemma which will imply the uniqueness.
\begin{lemma} Let $f$ and $g$ be two solutions with the same initial data. These two solutions are in the space
$\mathcal{C}([0,T], \dot W^{1,\infty} \cap \dot H^2)  \cap L^{2}([0,T],\dot H^2).$
 Then, if we set $\mathcal{U}:=f-g$,  $\mathcal{U}$ verifies the following Gronwall's type inequality
 \begin{eqnarray*}
\Vert \mathcal{U} \Vert_{L^{\infty} \dot H^1} \leq \Vert \mathcal{U}_{0} \Vert_{\dot H^1} \exp\left({c(K)( \Vert f \Vert^{2}_{L^\infty\dot H^2} + \Vert g \Vert^{2}_{L^\infty \dot H^2}) \left(\Vert f \Vert^2_{L^2_{T}\dot H^2}+\Vert g \Vert^2_{{L^2}_{T}\dot H^2}\right)} \right).
 \end{eqnarray*} 
\end{lemma}
Let $f$ and $g$ be two solutions of the 3D Muskat equation with the same initial data. Let $\mathcal{U}=f-g$, then $\mathcal{U}$ verifies
 \begin{eqnarray*}
\partial_{t} \mathcal{U} &=& \int \Delta_{y} \nabla_{x}\mathcal{U}. \frac{y}{\vert y \vert^{2}} \int_{0}^{\infty} e^{-k} \cos(k \Delta_{y} f) \cos(\arctan(\Delta_{y} f)) \ dk \ dy  \\
&+&  \int \Delta_{y} \nabla_{x}g. \frac{y}{\vert y \vert^{2}} \int_{0}^{\infty} e^{-k} \left[ \cos(k \Delta_{y} f) \cos(\arctan(\Delta_{y} f))-\cos(k \Delta_{y} g) \cos(\arctan(\Delta_{y} g)) \right]\ dk \  dy \\
\end{eqnarray*}
 
We shall do estimates in $H^{1}$ on $\mathcal{U}$. That is, we dot multiply the gradient of the evolution equation with $\nabla\mathcal{U}$ and integrate with respect to the space variable. We obtain
 \begin{eqnarray*}
\frac{1}{2}\partial_{t} \Vert \mathcal{U} \Vert^{2}_{\dot H^{1}} &=&\int \nabla\mathcal{U} . \int \nabla \left(\Delta_{y}  \nabla\mathcal{U}. \frac{y}{\vert y \vert^{2}} \int_{0}^{\infty} e^{-k} \cos(k \Delta_{y} f) \cos(\arctan(\Delta_{y} f)) \right) \ dk \ dy \ dx   \\
&+& \int \nabla\mathcal{U}. \int  \nabla \left(\Delta_{y} \nabla g. \frac{y}{\vert y \vert^{2}} \right.  \\ &&
\left. \times  \int_{0}^{\infty} e^{-k}\left[ \cos(k \Delta_{y} f) \cos(\arctan(\Delta_{y} f))-\cos(k \Delta_{y} g) \cos(\arctan(\Delta_{y} g)) \right]\ dk \  dy \ dx \right) \\
&=& \mathcal{A}_1 + \mathcal{A}_2
\end{eqnarray*}

\subsection{Estimate of $ \mathcal{A}_1$}

We first notice that the most singular term is when the gradient hits $\Delta_{y}  \nabla\mathcal{U}$ and a reminder which corresponds to the term where the gradient hits the oscillatory integrals. We start estimating the most singular term. To do so, we first write that

\begin{eqnarray*}
\mathcal{A}_1 &=& \int \nabla\mathcal{U} . \nabla \left(\int \Delta_{y}  \nabla\mathcal{U}. \frac{y}{\vert y \vert^{2}} \int_{0}^{\infty} e^{-k} \cos(k \Delta_{y} f) \cos(\arctan(\Delta_{y} f)) \ dk \ dy \ dx \right)  \\
&=&  \int \nabla\mathcal{U}.  \left(\int \nabla\Delta_{y}  \nabla\mathcal{U}. \frac{y}{\vert y \vert^{2}} \int_{0}^{\infty} e^{-k} \cos(k \Delta_{y} f) \cos(\arctan(\Delta_{y} f)) \ dk \ dy \ dx \right)  \\
&+&\int \nabla\mathcal{U}.  \int \left(\Delta_{y}  \nabla\mathcal{U}. \frac{y}{\vert y \vert^{2}} \right) \int_{0}^{\infty} e^{-k} \nabla\left(\cos(k \Delta_{y} f) \cos(\arctan(\Delta_{y} f))\right) \ dk \ dy \ dx \\
&=& \mathcal{A}_{1,1}+\mathcal{A}_{1,2}
\end{eqnarray*}

We start by estimating the more singular term, that is $\mathcal{A}_{1,1}$. We use the {\it{a priori estimates}}  in $\dot H^2$ obtained previously and  replace the first two $\Delta f$ by $\nabla\mathcal{U}$. We immediately find that 

\begin{eqnarray*}
\mathcal{A}_{1,1}&=&\underbrace{\int \nabla\mathcal{U} \ \int  \nabla\Delta_{y}  \nabla\mathcal{U} . \frac{y}{\vert y \vert^2} \cos(\arctan(\Delta_{y} f )) \int_{0}^{\infty} e^{-k} \cos(k\Delta_{y} f )  \ dk  \ dy  \ dx}\\ 
\end{eqnarray*}
Again, since the first Laplacian operator $\Delta$ does not play any role in the proof of the Lemma. We may replace it by the  nabla operator $\nabla$ and get that
\begin{eqnarray*}
\mathcal{A}_{1,1}&=&\frac{1}{2} \int  \nabla\mathcal{U}.  \int \nabla D_{y} \nabla\mathcal{U} . \frac{y}{\vert y \vert^2} \sin(\frac{1}{2}(\arctan(\Delta_{y} f )+\arctan(\bar\Delta_{y} f )))    \times \\
 &&\sin(\frac{1}{2}(\arctan(\Delta_{y} f )-\arctan(\bar\Delta_{y} f ))) \int_{0}^{\infty} e^{-k} \sin(\frac{k}{2} S_y f) \sin(\frac{k}{2} D_y f)    \ dk  \ dy  \\
  &+&\int  \nabla\mathcal{U}. \int \nabla\bar\Delta_{y}  \nabla\mathcal{U}  \frac{y}{\vert y \vert^2} 
\cos(\frac{1}{2}(\arctan(\bar\Delta_{y} f) + \arctan(\Delta_{y} f ))\times \\
&&\cos(\frac{1}{2}(\arctan(\bar\Delta_{y} f )- \arctan(\Delta_{y} f )) \int_{0}^{\infty} e^{-k} \sin(\frac{k}{2} S_y f) \sin(\frac{k}{2} D_y f)  \ dk  \ dy \\
&+& \int  \nabla\mathcal{U}\int \nabla\Delta_{y}  \nabla\mathcal{U} . \frac{y}{\vert y \vert^2} \sin(\frac{1}{2}(\arctan(\Delta_{y} f )+\arctan(\bar\Delta_{y} f )))    \times \\
&&\sin(\frac{1}{2}(\arctan(\Delta_{y} f )-\arctan(\bar\Delta_{y} f )))\int_{0}^{\infty} e^{-k} \cos(\frac{k}{2} S_y f) \cos(\frac{k}{2} D_y f)   \ dk  \ dy  \\
&+& \frac{1}{2} \int  \nabla\mathcal{U} \int  \nabla D_{y}  \nabla\mathcal{U} . \frac{y}{\vert y \vert^2} \cos(\frac{1}{2}(\arctan(\Delta_{y} f )+\arctan(\bar\Delta_{y} f ))) \int_{0}^{\infty} e^{-k}   \\
&&\cos(\frac{1}{2}(\arctan(\Delta_{y} f )-\arctan(\bar\Delta_{y} f )))  (\cos(\frac{k}{2}(D_y f ))(\cos(\frac{k}{2}(S_y f ))  \ dk  \ dy\\
&:=& \sum_{i=1}^4 \mathcal{A}_{1,1,i}(t) \\
\end{eqnarray*}

To estimates $\mathcal{A}_{1,1,1}$, we integrate by parts in $y$ and then estimate. The first term is when we differentiate the kernel, that is
\begin{eqnarray*}
\mathcal{A}_{1,1,1} &=&\frac{1}{2}\int  \nabla\mathcal{U} \int \left( s_{y} \nabla\mathcal{U} \right)\left(\nabla.\frac{y}{\vert y \vert^3}\right) \sin(\frac{1}{2}(\arctan(\Delta_{y} f )+\arctan(\bar\Delta_{y} f )))     \\
  &&     \sin(\frac{1}{2}(\arctan(\Delta_{y} f )-\arctan(\bar\Delta_{y} f ))) \int_{0}^{\infty} e^{-k} \sin(\frac{k}{2}S_y f) \sin(\frac{k}{2}D_y f)  \ dk  \ dy \ dx. \\
  &+&\frac{1}{2}\int  \nabla\mathcal{U} \int \left( s_{y} \nabla\mathcal{U} \right)\frac{1}{\vert y \vert^3} \ y.\nabla_{y}\left(\sin(\frac{1}{2}(\arctan(\Delta_{y} f )+\arctan(\bar\Delta_{y} f )))  \right.   \\
  &&  \left.   \sin(\frac{1}{2}(\arctan(\Delta_{y} f )-\arctan(\bar\Delta_{y} f ))) \int_{0}^{\infty} e^{-k} \sin(\frac{k}{2}S_y f) \sin(\frac{k}{2}D_y f) \right) \ dk  \ dy \ dx\\
  &=&  \mathcal{A}_{1,1,1,1}+\mathcal{A}_{1,1,1,2}
  \end{eqnarray*}
  Since $\left\vert \nabla.\frac{y}{\vert y \vert^3} \right \vert \lesssim \frac{1}{\vert y \vert^3} $ then we find  that,
  \begin{eqnarray*}
\mathcal{A}_{1,1,1,1} &\lesssim&\Vert \mathcal{U} \Vert^{2}_{L^{2}} \int \frac{\Vert s_{y} f \Vert^{2}_{L^{\infty}}}{\vert y \vert^{5}} \ dy \\
&\lesssim& \Vert \nabla\mathcal{U} \Vert^{2}_{L^{2}} \Vert f \Vert^{2}_{\dot B^{3/2}_{\infty,2}} \\
&\lesssim& \Vert \mathcal{U} \Vert^{2}_{\dot H^{1}} \Vert f \Vert^{2}_{\dot H^{5/2}}
  \end{eqnarray*}
  Otherwise, then we differentiate one of the oscillatory terms. In this case, we use Holder ($L^2-L^2-L^\infty-L^ \infty $) where one of the ${L^\infty}$ norm will necessary be in a term of order $\nabla_{x} \delta_{\alpha} f$ and the other one in any of the $s_{\alpha} f$. So, one finds

  \begin{eqnarray*}
\mathcal{A}_{1,1,1,2} &\lesssim&\Vert \nabla\mathcal{U} \Vert^2_{L^2} \int \frac{\Vert s_{y} f \Vert_{L^{\infty}} \Vert \delta_{y} \nabla f \Vert_{L^{\infty}}}{\vert y \vert^{4}} \ dy \\
 &\lesssim&\Vert \nabla\mathcal{U} \Vert^{2}_{L^{2}} \left(\int \frac{\Vert s_{y} f \Vert^2_{L^{\infty}}}{\vert y \vert^5} \ dy  \int \frac{ \Vert \delta_{y} \nabla f \Vert^{2}_{L^{\infty}}}{\vert y \vert^{3}} \ dy \right)^{1/2} \\
&\lesssim& \Vert \nabla\mathcal{U} \Vert^{2}_{L^{2}} \Vert f \Vert_{\dot B^{3/2}_{\infty,2}} \Vert \nabla f \Vert_{\dot B^{1/2}_{\infty,2}} \\
&\lesssim& \Vert \mathcal{U} \Vert^{2}_{\dot H^{1}} \Vert f \Vert^{2}_{\dot H^{5/2}}
  \end{eqnarray*}
  
  \subsubsection{Estimate of $\mathcal{A}_{1,2}$}
  
  Now, we estimates $\mathcal{A}_{1,1,2}$. We use the decomposition previously proved (see \eqref{A2}). We analogously find that 
  
  \begin{eqnarray} 
\mathcal{A}_{1,1,2}&=&-\frac{1}{2} \int \nabla \mathcal{U}\int   \delta_{y} \nabla\mathcal{U}
\frac{y}{\vert y \vert^3} .\nabla_{x}\left(\sin(\frac{1}{2}(\arctan(\Delta_{y} f )+\arctan(\bar\Delta_{y} f )))    \right. \times \nonumber \\
&& \left.\sin(\frac{1}{2}(\arctan(\Delta_{y} f )-\arctan(\bar\Delta_{y} f )))\int_{0}^{\infty} e^{-k} \cos(\frac{k}{2} S_y f) \cos(\frac{k}{2} D_y f) \right)  \ dk  \ dy \ dx  \nonumber \\
&&-\frac{1}{2} \int  \nabla\mathcal{U}\int  \nabla\mathcal{U}(x-y)  \frac{y}{\vert y \vert^3} .\nabla_{x}\left(\sin(\frac{1}{2}(\arctan(\Delta_{y} f )+\arctan(\bar\Delta_{y} f )))    \right. \times \nonumber \\
&& \left. \sin(\frac{1}{2}(\arctan(\Delta_{y} f )-\arctan(\bar\Delta_{y} f )))\int_{0}^{\infty} e^{-k} \cos(\frac{k}{2} S_y f) \cos(\frac{k}{2} D_y f) \right)   \ dk  \ dy \ dx  \nonumber  \\
&&-\int  \nabla\mathcal{U}\int \nabla_{x}  \nabla\mathcal{U}(x-y) . \frac{y}{\vert y \vert^3} \sin(\frac{1}{2}(\arctan(\Delta_{y} f )+\arctan(\bar\Delta_{y} f )))   \times \nonumber \\
&&\sin(\frac{1}{2}(\arctan(\Delta_{y} f )-\arctan(\bar\Delta_{y} f )))  \int_{0}^{\infty} e^{-k} \cos(\frac{k}{2} S_y f) \cos(\frac{k}{2} D_y f)   \ dk  \ dy \ dx \nonumber  \\
&=& \mathcal{A}_{1,2,1} + \mathcal{A}_{1,2,2} + \mathcal{A}_{1,2,3}.   
\end{eqnarray}

To estimate $\mathcal{A}_{1,1,2,1}$, we write
  \begin{eqnarray*}
  \mathcal{A}_{1,1,2,1} &\lesssim& \Vert \nabla\mathcal{U}  \Vert_{L^{2}} \int \frac{\Vert \delta_{y} \nabla\mathcal{U} \Vert_{L^{2}}}{\vert y \vert^{3/2}} \frac{\Vert \nabla \delta^{\pm}_{y} f \Vert_{L^{\infty}}}{\vert y \vert^{3/2}} \ dy  \\
  &\lesssim& \Vert \nabla\mathcal{U}  \Vert_{\dot H^{2}}  \Vert \nabla \mathcal{U} \Vert_{\dot B^{1/2}_{2,2}} \Vert \nabla f \Vert_{\dot B^{1/2}_{\infty,2}} \\
  &\lesssim&  \Vert \mathcal{U}  \Vert_{\dot H^1} \Vert \mathcal{U}  \Vert_{\dot H^{3/2}} \Vert f \Vert_{\dot H^{5/2}},  
   \end{eqnarray*} 
where $\delta^{\pm}_{y} f:=f(x)-f(x\pm y)$. \\

The estimate of $\mathcal{A}_{1,1,2,2}$ is done by using the decomposition 
\begin{eqnarray*} 
 \mathcal{A}_{1,1,2,2}&=&-\frac{1}{2} \int \nabla\mathcal{U} \int  \nabla\mathcal{U}(x-y)  \frac{y}{\vert y \vert^3} .\nabla_{x}\left(\sin(\frac{1}{2}(\arctan(\Delta_{y} f )+\arctan(\bar\Delta_{y} f )))    \right. \times \\
&& \left.  \hspace{-0.5cm}\sin(\frac{1}{2}(\arctan(\Delta_{y} f )-\arctan(\bar\Delta_{y} f )))\int_{0}^{\infty} e^{-k} \cos(\frac{k}{2} S_y f) \cos(\frac{k}{2} D_y f)   \ dk  \ dy \ dx \right)  \\
&=&-\frac{1}{2} \int \nabla\mathcal{U} \int  \nabla\mathcal{U}(x-y)  \frac{y}{\vert y \vert^3} . \sin(\frac{1}{2}(\arctan(\Delta_{y} f )+\arctan(\bar\Delta_{y} f ))) \times \\
&& \hspace{-1cm}\nabla_{x}\left(\sin(\frac{1}{2}(\arctan(\Delta_{y} f )-\arctan(\bar\Delta_{y} f )))\int_{0}^{\infty} e^{-k} \cos(\frac{k}{2} S_y f) \cos(\frac{k}{2} D_y f)   \ dk  \ dy \ dx \right) \\
&&-\frac{1}{2} \int \nabla\mathcal{U} \int \nabla \mathcal{U}(x-y)  \frac{y}{\vert y \vert^3} .\nabla_{x}\left(\sin(\frac{1}{2}(\arctan(\Delta_{y} f )+\arctan(\bar\Delta_{y} f )))\right)  \\
&&  \sin(\frac{1}{2}(\arctan(\Delta_{y} f )-\arctan(\bar\Delta_{y} f )))\int_{0}^{\infty} e^{-k} \cos(\frac{k}{2} S_y f) \cos(\frac{k}{2} D_y f)   \ dk  \ dy  \ dx \\
&=&\mathcal{A}_{1,1,2,2,1}+\mathcal{A}_{1,1,2,2,2}.
   \end{eqnarray*} 
   $\mathcal{A}_{1,1,2,2,1}$ is easy to estimate, indeed, we have
    \begin{eqnarray*}
\mathcal{A}_{1,1,2,2,1} &\lesssim&\Vert \nabla\mathcal{U} \Vert^2_{L^2} \int \frac{\Vert s_{y} f \Vert_{L^{\infty}} \Vert \delta^{\pm}_{y}  \nabla f  \Vert_{L^{\infty}}}{\vert y \vert^{4}} \ dy \\
 &\lesssim&\Vert \nabla\mathcal{U} \Vert^{2}_{L^{2}} \left(\int \frac{\Vert s_{y} f \Vert^2_{L^{\infty}}}{\vert y \vert^5} \ dy  \int \frac{ \Vert \delta^{\pm}_{y}  \nabla f \Vert^{2}_{L^{\infty}}}{\vert y \vert^{3}} \ dy \right)^{1/2} \\
&\lesssim& \Vert \nabla\mathcal{U} \Vert^{2}_{L^{2}} \Vert f \Vert_{\dot B^{3/2}_{\infty,2}} \Vert \nabla f \Vert_{\dot B^{1/2}_{\infty,2}} \\
&\lesssim& \Vert\mathcal{U} \Vert^{2}_{\dot H^1} \Vert f \Vert^{2}_{\dot H^{5/2}}
  \end{eqnarray*}
As for $\mathcal{A}_{1,1,2,2,2}$, we use the fact that $\nabla_{x}\left(\sin(\frac{1}{2}(\arctan(\Delta_{y} f )+\arctan(\bar\Delta_{y} f )))\right)$ may be written as follows (see \eqref{reste})
\begin{equation}
     \frac{\nabla_{x}\Delta_{y} f}{1+\Delta^{2}_{y} f}+ \frac{\nabla_{x}\bar\Delta_{y} f}{1+\bar\Delta^{2}_{y} f}={\frac{\nabla_{x}S_{y} f}{1+\Delta^{2}_{y} f}}+{\nabla_{x} D_{y}f \frac{S_{y} f \ D_{y} f}{(1+\Delta^{2}_{y} f)(1+\bar\Delta^{2}_{y} f)}}.
    \end{equation}
    The latter identity gives two terms which we call $\mathcal{A}_{1,1,2,2,2,1}$ and $\mathcal{A}_{1,1,2,2,2,2}$, namely

      \begin{eqnarray*}
\mathcal{A}_{1,1,2,2,2,1}&=&-\frac{1}{2} \int \nabla\mathcal{U} \int  \nabla\mathcal{U}(x-y)  \frac{y}{\vert y \vert^3} .{\frac{\nabla_{x}S_{y} f}{1+\Delta^{2}_{y} f}}  \\
&&  \sin(\frac{1}{2}(\arctan(\Delta_{y} f )-\arctan(\bar\Delta_{y} f )))\int_{0}^{\infty} e^{-k} \cos(\frac{k}{2} S_y f) \cos(\frac{k}{2} D_y f)   \ dk  \ dy \ dx \\
 \end{eqnarray*}
 and 
  \begin{eqnarray*}
\mathcal{A}_{1,1,2,2,2,2}&=&-\frac{1}{2} \int \nabla\mathcal{U} \int  \nabla\mathcal{U}(x-y)  \frac{y}{\vert y \vert^3} .{\nabla_{x} D_{y}f \frac{S_{y} f \ D_{y} f}{(1+\Delta^{2}_{y} f)(1+\bar\Delta^{2}_{y} f)}}  \\
&&  \sin(\frac{1}{2}(\arctan(\Delta_{y} f )-\arctan(\bar\Delta_{y} f )))\int_{0}^{\infty} e^{-k} \cos(\frac{k}{2} S_y f) \cos(\frac{k}{2} D_y f)   \ dk  \ dy \ dx \\
 \end{eqnarray*}
 To estimate $\mathcal{A}_{1,1,2,2,2,1}$ we use the following lemma whose proof is completely analogous to Lemma (see \eqref{R})

\begin{lemma}  The term $\mathcal{A}_{1,1,2,2,2,2}$ may be rewritten as follows,
\begin{eqnarray*}
\mathcal{A}_{1,1,2,2,2,2}&=&\frac{1}{4} \int \nabla\mathcal{U} \int  (\nabla\mathcal{U}(x-y)-\nabla\mathcal{U}(x+y))  \frac{y}{\vert y \vert^3} .{\nabla_{x}S_{y} f}  \\&&\times\sin(\frac{1}{2}(\arctan(\Delta_{y} f )-\arctan(\bar\Delta_{y} f )))  \int_{0}^{\infty} e^{-\gamma}\sin(\frac{\gamma}{2} S_y f)\\
&&  \sin(\frac{\gamma}{2} D_y f) \int_{0}^{\infty} e^{-k} \cos(\frac{k}{2} S_y f) \cos(\frac{k}{2} D_y f) \ d\gamma \ dk  \ dy \ dx \\
&-&\frac{1}{2} \int \nabla\mathcal{U} \int  \nabla\mathcal{U}(x+y)  \frac{y}{\vert y \vert^3} .{\nabla_{x}S_{y} f}  \sin(\frac{1}{2}(\arctan(\Delta_{y} f )-\arctan(\bar\Delta_{y} f ))) \times\\
&&  \hspace{-0.5cm}\int_{0}^{\infty} e^{-\gamma} 
\cos(\frac{\gamma}{2} S_y f) \cos(\frac{\gamma}{2} D_y f) \int_{0}^{\infty} e^{-k} 
\cos(\frac{k}{2} S_y f) \cos(\frac{k}{2} D_y f) \ d\gamma \ dk  \ dy \ dx \\
&=&\mathcal{A}_{1,1,2,2,2,2,1} + \mathcal{A}_{1,1,2,2,2,2,2}
\end{eqnarray*}
\end{lemma}
The first term $\mathcal{A}_{1,1,2,2,2,2,1}$ is easy to estimated, indeed, it suffices to observe that 

 \begin{eqnarray*}
\mathcal{A}_{1,1,2,2,2,2,1} &\lesssim&\Vert \nabla\mathcal{U} \Vert^2_{L^2} \int \frac{\Vert s_{y} f \Vert_{L^{\infty}} \Vert \delta_{y} \nabla f \Vert_{L^{\infty}}}{\vert y \vert^{4}} \ dy \\
 &\lesssim&\Vert \nabla\mathcal{U} \Vert^{2}_{L^{2}} \left(\int \frac{\Vert s_{y} f \Vert^2_{L^{\infty}}}{\vert y \vert^5} \ dy  \int \frac{ \Vert \delta_{y} \nabla f \Vert^{2}_{L^{\infty}}}{\vert y \vert^{3}} \ dy \right)^{1/2} \\
&\lesssim& \Vert \nabla\mathcal{U} \Vert^{2}_{L^{2}} \Vert f \Vert_{\dot B^{3/2}_{\infty,2}} \Vert \nabla f \Vert_{\dot B^{1/2}_{\infty,2}} \\
&\lesssim& \Vert \mathcal{U} \Vert^{2}_{\dot H^{1}} \Vert f \Vert^{2}_{\dot H^{5/2}}
  \end{eqnarray*}

As for the second one, we first observe that by using the change of variables $y \rightarrow -y$, one obtains
\begin{eqnarray*}
\mathcal{A}_{1,1,2,2,2,2,2}&=&-\frac{1}{4} \int \nabla\mathcal{U} \int  \nabla(\mathcal{U}(x+y)+\mathcal{U}(x-y))  \frac{y}{\vert y \vert^3} .{\nabla_{x}S_{y} f}  \sin(\frac{1}{2}(\arctan(\Delta_{y} f )-\arctan(\bar\Delta_{y} f ))) \times\\
&&  \hspace{-0.5cm}\int_{0}^{\infty} e^{-\gamma} 
\cos(\frac{\gamma}{2} S_y f) \cos(\frac{\gamma}{2} D_y f) \int_{0}^{\infty} e^{-k} 
\cos(\frac{k}{2} S_y f) \cos(\frac{k}{2} D_y f) \ d\gamma \ dk  \ dy \ dx 
\end{eqnarray*}

Then, by using the fact that $\nabla(\mathcal{U}(x-y)+U(x+y))=\nabla_{y} s_{y}\mathcal{U}$, we may integrate by parts in $y$ and find that 
\begin{eqnarray*}
\mathcal{A}_{1,1,2,2,2,2,2}&=&-\frac{1}{4} \int \nabla\mathcal{U} \int  \nabla(\mathcal{U}(x+y)+\mathcal{U}(x-y))  \frac{y}{\vert y \vert^3} .{\nabla_{x}S_{y} f}  \sin(\frac{1}{2}(\arctan(\Delta_{y} f )-\arctan(\bar\Delta_{y} f ))) \times\\
&&  \hspace{-0.5cm}\int_{0}^{\infty} e^{-\gamma} 
\cos(\frac{\gamma}{2} S_y f) \cos(\frac{\gamma}{2} D_y f) \int_{0}^{\infty} e^{-k} 
\cos(\frac{k}{2} S_y f) \cos(\frac{k}{2} D_y f) \ d\gamma \ dk  \ dy \ dx  \\
&\lesssim& \Vert \nabla\mathcal{U} \Vert_{L^{2}} \left(\int \frac{\Vert s_{y}\mathcal{U} \Vert_{L^{\infty}}}{\vert y \vert^{3/2}} \frac{ \Vert \nabla s_y f \Vert_{L^{2}}}{\vert y \vert^{5/2}} \ dy +\int \frac{\Vert s_{y}\mathcal{U} \Vert_{L^{\infty}}}{\vert y \vert^{3/2}} \frac{ \Vert \Delta s_y f \Vert_{L^{2}}}{\vert y \vert^{3/2}} \ dy \right) \\
&\lesssim& \Vert \mathcal{U} \Vert_{\dot H^{1}}\left(\Vert \mathcal{U} \Vert_{\dot B^{1/2}_{\infty,2}}
\Vert \nabla f \Vert_{\dot B^{3/2}_{2,2}} + \Vert \mathcal{U} \Vert_{\dot B^{1/2}_{\infty,2}} \Vert \Delta f \Vert_{\dot B^{1/2}_{2,2}}\right) \\
&\lesssim& \Vert \mathcal{U} \Vert_{\dot H^{1}}\Vert \mathcal{U} \Vert_{\dot H^{3/2}}
\Vert f \Vert_{\dot H^{5/2}} 
\end{eqnarray*}
The control of $\mathcal{A}_{1,1,3}$ is the same as the one of   $\mathcal{A}_{1,1,2}$ since they are the equal up to interchanging the role of one sine and cosine (they are just bounded by 1 in all the steps). Hence, we have
\begin{eqnarray*}
\mathcal{A}_{1,3} \lesssim  \Vert \mathcal{U} \Vert_{\dot H^{1}}\Vert \mathcal{U} \Vert_{\dot H^{3/2}}
\Vert f \Vert_{\dot H^{5/2}}  + \Vert \mathcal{U} \Vert^{2}_{\dot H^{1}} \Vert f \Vert^{2}_{\dot H^{5/2}}
\end{eqnarray*}
The last term, that $\mathcal{A}_{1,4}$ contains the dissipation. More precisely, we have that
\begin{eqnarray*}
\mathcal{A}_{1,1,4}&=&\frac{1}{2} \int  \nabla\mathcal{U} \int  \nabla D_{y}  \nabla\mathcal{U} . \frac{y}{\vert y \vert^2} \cos(\frac{1}{2}(\arctan(\Delta_{y} f )+\arctan(\bar\Delta_{y} f ))) \int_{0}^{\infty} e^{-k}   \\
&&\cos(\frac{1}{2}(\arctan(\Delta_{y} f )-\arctan(\bar\Delta_{y} f )))  (\cos(\frac{k}{2}(D_y f ))(\cos(\frac{k}{2}(S_y f ))  \ dk  \ dy \ dx
\end{eqnarray*}
We start by linearizing, 
\begin{eqnarray*}
\mathcal{A}_{1,1,4}&=& \frac{1}{2} \int \nabla \mathcal{U} \int   \nabla D_{y} \nabla \mathcal{U} . \frac{y}{\vert y \vert^2} \cos(\frac{1}{2}(\arctan(\Delta_{y} f )+\arctan(\bar\Delta_{y} f )))   \times  \\
&&\cos(\frac{1}{2}(\arctan(\Delta_{y} f ))-(\arctan(\bar\Delta_{y} f )))\int_{0}^{\infty} e^{-k} (\cos(\frac{k}{2}(D_y f ))  \ dk  \ dy \ dx\\
&-&   \int \nabla \mathcal{U} \int  \nabla D_{y} \nabla  \mathcal{U} . \frac{y}{\vert y \vert^2} \cos(\frac{1}{2}(\arctan(\Delta_{y} f )+\arctan(\bar\Delta_{y} f )))    \\
&&\cos(\frac{1}{2}(\arctan(\Delta_{y} f )-\arctan(\bar\Delta_{y} f )))\int_{0}^{\infty} e^{-k} \cos(\frac{k}{2}(D_y f )\sin^{2}(\frac{k}{4}(S_y f )  \ dk  \ dy \ dx\\
&=&- \int \nabla \mathcal{U} \int  \nabla D_{y} \nabla \mathcal{U} . \frac{y}{\vert y \vert^2} \sin^{2}(\frac{1}{4}(\arctan(\Delta_{y} f )+\arctan(\bar\Delta_{y} f )))     \\
&&\cos(\frac{1}{2}(\arctan(\Delta_{y} f )-\arctan(\bar\Delta_{y} f ))) \times \int_{0}^{\infty} e^{-k} (\cos(\frac{k}{2}(D_y f ))  \ dk  \ dy  \ dx\\
&-& \int \nabla \mathcal{U} \int  \nabla D_{y} \nabla  \mathcal{U} . \frac{y}{\vert y \vert^2} \cos(\frac{1}{2}(\arctan(\Delta_{y} f )+\arctan(\bar\Delta_{y} f )))    \times  \\
&&\cos(\frac{1}{2}(\arctan(\Delta_{y} f )-\arctan(\bar\Delta_{y} f )))\int_{0}^{\infty} e^{-k} \cos(\frac{k}{2}(D_y f )\sin^{2}(\frac{k}{4}(S_y f )  \ dk  \ dy \ dx\\
&+& \frac{1}{2}\int \nabla \mathcal{U} \int  \nabla D_{y} \nabla  \mathcal{U} . \frac{y}{\vert y \vert^2}   \cos(\frac{1}{2}(\arctan(\Delta_{y} f )-\arctan(\bar\Delta_{y} f ))) \\
&\times&\int_{0}^{\infty} e^{-k} (\cos(\frac{k}{2}(D_y f ))  \ dk  \ dy \ dx\\
&=& \mathcal{A}_{1,1,4,1}+\mathcal{A}_{1,1,4,2}+\mathcal{A}_{1,1,4,3}.
\end{eqnarray*}
To estimate $\mathcal{A}_{1,1,4,1}$, we balance the derivative in $x$ by using again the fact that $\delta_{y} \nabla_x d_y=-\delta_{y} \nabla_y s_y$, then integrating by parts in $y$ gives 
\begin{eqnarray*}
\mathcal{A}_{1,1,4,1}&=& \int \nabla \mathcal{U} \int \frac{\nabla s_{y} \mathcal{U} }{\vert y \vert^3} \sin^{2}(\frac{1}{4}(\arctan(\Delta_{y} f )+\arctan(\bar\Delta_{y} f )))  \times  \\
&& \cos(\frac{1}{2}(\arctan(\Delta_{y} f )-\arctan(\bar\Delta_{y} f ))) \int_{0}^{\infty} e^{-k} (\cos(\frac{k}{2}(D_y f ))  \ dk  \ dy \ dx\\
&-&\int \nabla \mathcal{U} \int \frac{\nabla s_{y} \mathcal{U} }{\vert y \vert^3} y.\nabla_y(\arctan(\Delta_{y} f )+\arctan(\bar\Delta_{y} f )) \\ && \times\sin(\frac{1}{2}(\arctan(\Delta_{y} f )+\arctan(\bar\Delta_{y} f )))   \\
&&  \cos(\frac{1}{2}(\arctan(\Delta_{y} f )-\arctan(\bar\Delta_{y} f ))) \int_{0}^{\infty} e^{-k} (\cos(\frac{k}{2}(D_y f ))  \ dk  \ dy \ dx\\
&+& \int \nabla \mathcal{U} \int \frac{\nabla s_{y} \mathcal{U} }{\vert y \vert^3} \sin^{2}(\frac{1}{4}(\arctan(\Delta_{y} f )+\arctan(\bar\Delta_{y} f ))) \\
&& \times y.\nabla_y(\arctan(\Delta_{y} f )+\arctan(\bar\Delta_{y} f )) \times  \\
&&\sin((\arctan(\Delta_{y} f )-\arctan(\bar\Delta_{y} f ))) \int_{0}^{\infty} e^{-k} (\cos(\frac{k}{2}(D_y f ))  \ dk  \ dy \ dx\\
&+&\frac{1}{2}\int \nabla \mathcal{U} \int \frac{\nabla s_{y} \mathcal{U} }{\vert y \vert^3} \sin^{2}(\frac{1}{4}(\arctan(\Delta_{y} f )+\arctan(\bar\Delta_{y} f )))    
 \\
&& \times\cos(\frac{1}{2}(\arctan(\Delta_{y} f )-\arctan(\bar\Delta_{y} f ))) y.\nabla_y(D_y f)\\ 
&& \times \int_{0}^{\infty} ke^{-k} \sin(\frac{k}{2}(D_y f ))  \ dk  \ dy \ dx\\
&=&\mathcal{A}_{1,1,4,1,1}+\mathcal{A}_{1,1,4,1,2}+\mathcal{A}_{1,1,4,1,3}+\mathcal{A}_{1,1,4,1,4}.
\end{eqnarray*}
The estimate of $\mathcal{A}_{1,1,4,1,1}$ is easy, indeed, it suffices to write that
\begin{eqnarray*}
\mathcal{A}_{1,1,4,1,1} &\lesssim&\Vert \nabla\mathcal{U} \Vert^2_{L^2} \int \frac{\Vert s_{y} f \Vert^2_{L^{\infty}} }{\vert y \vert^{5}} \ dy \\
&\lesssim& \Vert \mathcal{U} \Vert^{2}_{\dot H^{1}} \Vert f \Vert^{2}_{\dot H^{5/2}}
  \end{eqnarray*}
  Using Lemma \eqref{astuce}  	 \eqref{som} and \eqref{ddiff} for $i=2,3,4$ one has
  \begin{eqnarray*}
\mathcal{A}_{1,1,4,1,i} &\lesssim&\Vert \nabla\mathcal{U} \Vert^2_{L^2} \int \frac{\Vert s_{y} f \Vert_{L^{\infty}} \Vert \delta^{\pm}_{y}  \nabla f  \Vert_{L^{\infty}}}{\vert y \vert^{4}} \ dy \\
 &\lesssim&\Vert \nabla\mathcal{U} \Vert^{2}_{L^{2}} \left(\int \frac{\Vert s_{y} f \Vert^2_{L^{\infty}}}{\vert y \vert^5} \ dy  \int \frac{ \Vert \delta^{\pm}_{y}  \nabla f \Vert^{2}_{L^{\infty}}}{\vert y \vert^{3}} \ dy \right)^{1/2} \\
&\lesssim& \Vert \nabla\mathcal{U} \Vert^{2}_{L^{2}} \Vert f \Vert_{\dot B^{3/2}_{\infty,2}} \Vert \nabla f \Vert_{\dot B^{1/2}_{\infty,2}} \\
&\lesssim& \Vert\mathcal{U} \Vert^{2}_{\dot H^1} \Vert f \Vert^{2}_{\dot H^{5/2}}
  \end{eqnarray*}
  Then, notice that $\mathcal{A}_{1,1,4,1}$ and $\mathcal{A}_{1,4,2}$ have exactly the same regularity in the sense that the terms $\sin^{2}(\frac{1}{4}(\arctan(\Delta_{y} f )+\arctan(\bar\Delta_{y} f )))$ and $\sin^{2}(\frac{k}{4}(S_y f )$ have the same regularity. Indeed they are both bounded by $c \vert S_{\alpha} f \vert^2$ where $c>0$ is a constant. Hence, we conclude that
  \begin{eqnarray*}
\mathcal{A}_{1,1,4,2} \lesssim \Vert\mathcal{U} \Vert^{2}_{\dot H^1} \Vert f \Vert^{2}_{\dot H^{5/2}}
\end{eqnarray*}

The dissipation comes from the term $\mathcal{A}_{1,1,4,3}$ which is analogous to the term $\mathcal{S}_{4,3}$ (see \eqref{dissi}). By replacing the first two $\Delta f$ in $\mathcal{S}_{4,3}$ by $\nabla \mathcal{U}$ we immediately find that

\begin{eqnarray*}
\mathcal{A}_{1,1,4,3}&=&-\frac{1}{2}\int \nabla \mathcal{U}  \int    \frac{\nabla\delta_{y} \mathcal{U}}{\vert y \vert^3}   \left(\cos((\arctan(\frac{y}{\vert y \vert}. \nabla f(x)))-\cos((\arctan(\frac{y}{\vert y \vert}. \nabla f(x-y)))\right) \\
&& \times \int_{0}^{\infty} e^{-k} \cos({k} \frac{y}{\vert y \vert}. \nabla f(x))      \ dk  \ dy \ dx \\
&-&\frac{1}{4}\int \nabla \mathcal{U} \int    \frac{\nabla\delta_{y} \mathcal{U}}{\vert y \vert^3}   \left(\cos((\arctan(\frac{y}{\vert y \vert}. \nabla f(x)))+\cos((\arctan(\frac{y}{\vert y \vert}. \nabla f(x-y)))\right) \\
&& \int_{0}^{\infty} e^{-k} \left(\cos({k} \frac{y}{\vert y \vert}. \nabla f(x))-\cos({k} \frac{y}{\vert y \vert}. \nabla f(x-y)) \right)     \ dk  \ dy \ dx \\
&-&\frac{1}{4}\int \nabla \mathcal{U}  \int    \frac{\nabla\delta_{y} \mathcal{U}}{\vert y \vert^3}   \left(\cos((\arctan(\frac{y}{\vert y \vert}. \nabla f(x)))+\cos((\arctan(\frac{y}{\vert y \vert}. \nabla f(x-y)))\right) \\
&& \int_{0}^{\infty} e^{-k} \left(\cos({k} \frac{y}{\vert y \vert}. \nabla f(x))+\cos({k} \frac{y}{\vert y \vert}. \nabla f(x-y)) \right)      \ dk  \ dy \ dx \\
&=& \mathcal{A}_{1,1,4,3,1}+\mathcal{A}_{1,1,4,3,2}+\mathcal{A}_{1,1,4,3,3}
      \end{eqnarray*} 
      The first two terms are easy to control, indeed, one has, for $i=1,2$
      \begin{eqnarray*}
\mathcal{A}_{1,1,4,3,i} &\lesssim& \Vert \mathcal{U} \Vert_{\dot H^{1}} \int \frac{\Vert \delta_{y} \nabla \mathcal{U} \Vert_{L^{2}}}{\vert y \vert^{3/2}}\frac{\Vert \delta_{y} \nabla f \Vert_{L^{\infty}}}{\vert y \vert^{3/2}} \ dy \\
&\lesssim& \Vert \mathcal{U} \Vert_{\dot H^{1}} \Vert \mathcal{U} \Vert_{\dot H^{3/2}} \Vert f \Vert_{\dot B^{3/2}_{\infty,2}} \\
&\lesssim&  \Vert \mathcal{U} \Vert_{\dot H^{1}} \Vert \mathcal{U} \Vert_{\dot H^{3/2}} \Vert f \Vert_{\dot H^{2}}
\end{eqnarray*}
The term $\mathcal{A}_{1,1,4,3,3}$ is the dissipative term. Following the same step as $\mathcal{S}_{4,3,3}$, one finds

  \begin{eqnarray*}
    \mathcal{A}_{1,1,4,3,3}  &=&\frac{1}{8}\int  \int    \frac{\vert \nabla\delta_{y} \mathcal{U}\vert^{2}}{\vert y \vert^3} \\
&&\times \left(-4+4-\frac{1}{1+(\frac{y}{\vert y \vert}. \nabla f(x-y))^{2})^{3/2}} - \frac{1}{1+(\frac{y}{\vert y \vert}. \nabla f(x))^{2})^{3/2}} \right. \\
&& \left. \hspace{4cm} -2 \frac{1}{\sqrt{1+(\frac{y}{\vert y \vert}. \nabla f(x-y))^2}} \frac{1}{1+(\frac{y}{\vert y \vert}. \nabla f(x))^{2}} \right) \ dy \ dx \\
&=& -\frac{1}{2} \Vert \mathcal{U} \Vert^2_{\dot H^{3/2}} + \frac{1}{8}\int  \int    \frac{\vert \nabla\delta_{y} \mathcal{U}\vert^{2}}{\vert y \vert^3} \\
&&\times \left(4-\frac{1}{1+(\frac{y}{\vert y \vert}. \nabla f(x-y))^{2})^{3/2}} - \frac{1}{1+(\frac{y}{\vert y \vert}. \nabla f(x))^{2})^{3/2}} \right. \\
&&\left. \hspace{4cm} -2 \frac{1}{\sqrt{1+(\frac{y}{\vert y \vert}. \nabla f(x-y))^2}} \frac{1}{1+(\frac{y}{\vert y \vert}. \nabla f(x))^{2}} \right) \ dy \ dx \\
&\leq&-\frac{1}{2} \Vert \mathcal{U} \Vert^2_{\dot H^{3/2}} + \frac{1}{2} \Vert \mathcal{U} \Vert^2_{\dot H^{3/2}} \left(1- \frac{1}{(1+K^{2})^{3/2}}\right)
      \end{eqnarray*}

  Now we estimate $\mathcal{A}_{1,2}$ that is
  \begin{eqnarray*}
\mathcal{A}_{1,2}&=&\int \nabla\mathcal{U}.  \int \left(\Delta_{y}  \nabla\mathcal{U}. \frac{y}{\vert y \vert^{2}} \right) \int_{0}^{\infty} e^{-k} \nabla\left(\cos(k \Delta_{y} f) \cos(\arctan(\Delta_{y} f))\right) \ dk \ dy \ dx \\
\end{eqnarray*}
Note that, it suffices to treat the case where the gradient hits $\cos(k \Delta_{y} f)$ since the other term is analogous.
\begin{eqnarray*}
\mathcal{A}_{1,2}&\lesssim& \Vert \nabla\mathcal{U} \Vert_{L^{2}} \int \frac{\Vert \delta_{y} \mathcal{U} \Vert_{L^{2}}}{ \vert y\vert^{3/2} } \frac{\Vert \delta_{y} f \Vert_{L^{\infty}}}{\vert y \vert^{3/2}} \ dy \\
&\lesssim&  \Vert \mathcal{U} \Vert_{\dot H^{1}}  \left(\int \frac{\Vert \delta_{y} \mathcal{U} \Vert^2_{L^{2}}}{ \vert y\vert^{3} } \ dy \right)^{1/2} \left(\int \frac{\Vert \delta_{y} \nabla f \Vert^2_{L^{\infty}}}{ \vert y\vert^{3} } \ dy \right)^{1/2}  \\
&\lesssim& \Vert \mathcal{U} \Vert_{\dot H^{1}}\Vert \mathcal{U} \Vert_{\dot H^{3/2}}
\Vert f \Vert_{\dot H^{5/2}}
\end{eqnarray*}

  \subsection{Estimate of $\mathcal{A}_2$}
For $\mathcal{A}_2$, introduce the operator
\begin{eqnarray*}
\mathcal{S}(f,g):=\int \nabla\mathcal{U} . \nabla \left(\int   \Delta_{y}  \nabla f. \frac{y}{\vert y \vert^2} \cos(\arctan(\Delta_{y} g ))\int_{0}^{\infty} e^{-k} \cos(k\Delta_{y} g ) \ dk  \ dy \ dx  \right)
\end{eqnarray*}
One easily notices that we have $\mathcal{A}_{2}=\mathcal{S}(g,f)-\mathcal{S}(g,g)$ and therefore, as a direct application of the Lemma (see \eqref{dec}) we may write that
\begin{eqnarray*}
\mathcal{S}(g,f)-\mathcal{S}(g,g)&=&\frac{1}{8} \int \nabla \mathcal{U}  . \nabla\left(\int  \Delta_{y} \nabla g-\bar\Delta_{y} \nabla g) . \frac{y}{\vert y \vert^2} \left(\cos(\arctan(\Delta_{y} f ))- \cos(\arctan(\bar\Delta_{y} f ))\right) \right. \\
  && \left. \times \int_{0}^{\infty} e^{-k} \left(\cos(k\Delta_{y} f )-\cos(k\bar\Delta_{y} f ) \right) \ dk  \ dy \ dx \right)  \\
  &-&\frac{1}{8} \int \nabla \mathcal{U} .  \nabla \left(\int  (\Delta_{y} \nabla g- \bar\Delta_{y} \nabla g) . \frac{y}{\vert y \vert^2} (\cos(\arctan(\Delta_{y} g ))- \cos(\arctan(\bar\Delta_{y} g ))) \right.\\
  && \left. \times \int_{0}^{\infty} e^{-k} (\cos(k\Delta_{y} g )-\cos(k\bar\Delta_{y} g) ) \ dk  \ dy \right) dx   \\
  &+& \frac{1}{8} \int \nabla \mathcal{U}  .\nabla \int ( \Delta_{y} \nabla g-\nabla \bar\Delta_{y} g) . \frac{y}{\vert y \vert^2}  \left(\cos(\arctan(\Delta_{y} f ))+\cos(\arctan(\bar\Delta_{y} f ))\right) \\
  && \int_{0}^{\infty} e^{-k}( \cos(k\bar\Delta_{y} f )+\cos(k\Delta_{y} f ) ) \ dk  \ dy \ dx  \\
  &-& \frac{1}{8} \int \nabla \mathcal{U} . \nabla  \int (\Delta_{y} \nabla g- \bar\Delta_{y} \nabla g) . \frac{y}{\vert y \vert^2}  \left(\cos(\arctan(\Delta_{y} g ))+\cos(\arctan(\bar\Delta_{y} g ))\right) \\
  && \times \int_{0}^{\infty} e^{-k}( \cos(k\bar\Delta_{y} g )+\cos(k\Delta_{y} g ) ) \ dk  \ dy \ dx  \\
  &-&\frac{1}{4} \int \nabla \mathcal{U}. \nabla  \int\bar\Delta_{y} \nabla g . \frac{y}{\vert y \vert^2} \left(\cos(\arctan(\bar\Delta_{y} f ))+\cos(\arctan(\Delta_{y} f ))\right) \times \\
  && \int_{0}^{\infty} e^{-k} (\cos(k\bar\Delta_{y} f )-\cos(k\Delta_{y} f ))  \ dk  \ dy \ dx \\
  &+&\frac{1}{4} \int \nabla \mathcal{U} . \nabla \int\bar\Delta_{y} \nabla g . \frac{y}{\vert y \vert^2} \left(\cos(\arctan(\bar\Delta_{y} g ))+\cos(\arctan(\Delta_{y} g ))\right) \times \\
  && \int_{0}^{\infty} e^{-k} (\cos(k\bar\Delta_{y} g )-\cos(k\Delta_{y} g ))  \ dk  \ dy\ dx \\
&-&\frac{1}{4}\int \nabla \mathcal{U} . \nabla \int\Delta_{y} \nabla g . \frac{y}{\vert y \vert^2} (\cos(\arctan(\bar\Delta_{y} f ))-\cos(\arctan(\Delta_{y} f )))\times\\
&& \int_{0}^{\infty} e^{-k} (\cos(k\Delta_{y} f) + \cos(k\bar\Delta_{y} f)) \ dk  \ dy \ dx \\
&+&\frac{1}{4}\int \nabla \mathcal{U} . \nabla\int \Delta_{y} \nabla g . \frac{y}{\vert y \vert^2} (\cos(\arctan(\bar\Delta_{y} g ))-\cos(\arctan(\Delta_{y} g )))\times\\
&&\int_{0}^{\infty} e^{-k} (\cos(k\Delta_{y} g) + \cos(k\bar\Delta_{y} g)) \ dk  \ dy \ dx   \\
&=&\sum_{i=1}^{8} \mathcal{A}_{2,1,{i}}
\end{eqnarray*}
We shall consider $\mathcal{A}_{2,1,i}$ and  $\mathcal{A}_{2,1,i+1}$ for $i=1...7$ and find some nice cancellations. More precisely, we write that
\begin{eqnarray*}
\mathcal{A}_{2,1,{1}}&=&\frac{1}{8} \int \nabla \mathcal{U}  . \nabla\left(\int  \Delta_{y} \nabla g-\bar\Delta_{y} \nabla g) . \frac{y}{\vert y \vert^2} \left(\cos(\arctan(\Delta_{y} f ))- \cos(\arctan(\bar\Delta_{y} f ))\right) \right. \\
  && \left. \times \int_{0}^{\infty} e^{-k} \left(\cos(k\Delta_{y} f )-\cos(k\bar\Delta_{y} f ) \right) \ dk  \ dy \ dx \right)  \\
  &=&\frac{1}{8} \int \nabla \mathcal{U}  . \nabla  \int  (\nabla\Delta_{y} g- \nabla \bar\Delta_{y} g) . \frac{y}{\vert y \vert^2} (\cos(\arctan(\Delta_{y} f ))-\cos(\arctan(\Delta_{y} g ))\\
  && +\cos(\arctan(\bar\Delta_{y} g ))- \cos(\arctan(\bar\Delta_{y} f )))  \int_{0}^{\infty} e^{-k} (\cos(k\Delta_{y} f )-\cos(k\bar\Delta_{y} f ) ) \ dk  \ dy \ dx  \\
  &+&\frac{1}{8}   \int \nabla \mathcal{U}  . \nabla    \int (\nabla  \Delta_{y} g- \nabla \bar\Delta_{y} g) . \frac{y}{\vert y \vert^2} (\cos(\arctan(\Delta_{y} g ))-\cos(\arctan(\bar\Delta_{y} g ))) \\
  && \times \int_{0}^{\infty} e^{-k} (\cos(k\Delta_{y} f )-\cos(k\bar\Delta_{y} f ) ) \ dk  \ dy \ dx  \\
\end{eqnarray*}
On the other hand, for $\mathcal{A}_{2,1,2}$, we may write
\begin{eqnarray*}
\mathcal{A}_{2,1,2}&=&-\frac{1}{8} \int \nabla \mathcal{U}. \nabla \int  (\nabla\Delta_{y} g-\nabla \bar\Delta_{y} g) . \frac{y}{\vert y \vert^2} (\cos(\arctan(\Delta_{y} g ))- \cos(\arctan(\bar\Delta_{y} g ))) \\
  && \times \int_{0}^{\infty} e^{-k} (\cos(k\Delta_{y} g )-\cos(k\Delta_{y} f )+\cos(k\bar\Delta_{y} f )-\cos(k\bar\Delta_{y} g) ) \ dk  \ dy \ dx  \\
  &-&\frac{1}{8}  \int \nabla  \mathcal{U} . \nabla \int  (\nabla\Delta_{y} g-\nabla \bar\Delta_{y} g) . \frac{y}{\vert y \vert^2} (\cos(\arctan(\Delta_{y} g ))- \cos(\arctan(\bar\Delta_{y} g ))) \\
  && \times \int_{0}^{\infty} e^{-k} (\cos(k\Delta_{y} f )-\cos(k\bar\Delta_{y} f ) ) \ dk  \ dy \ dx  \\
  \end{eqnarray*}
  Hence, noticing that the second term in $\mathcal{A}_{2,1,{1}}$ and $\mathcal{A}_{2,1,{2}}$ cancels out, one finds
  \begin{eqnarray*}
\mathcal{A}_{2,1,{1}}+\mathcal{A}_{2,1,{2}}&=&\frac{1}{8} \int \nabla \mathcal{U}. \nabla  \int  (\nabla_{x}\Delta_{y} g-\nabla_{x} \bar\Delta_{y} g) . \frac{y}{\vert y \vert^2} (\cos(\arctan(\Delta_{y} f ))-\cos(\arctan(\Delta_{y} g ))\\
  &+& \cos(\arctan(\bar\Delta_{y} g ))- \cos(\arctan(\bar\Delta_{y} f )))  \int_{0}^{\infty} e^{-k} (\cos(k\Delta_{y} f )-\cos(k\bar\Delta_{y} f ) ) \ dk  \ dy \ dx  \\
  &-&\frac{1}{8}  \int \nabla \mathcal{U}. \nabla \int  (\nabla_{x}\Delta_{y} g-\nabla_{x} \bar\Delta_{y} g) . \frac{y}{\vert y \vert^2} (\cos(\arctan(\Delta_{y} g ))- \cos(\arctan(\bar\Delta_{y} g ))) \\
  && \times \int_{0}^{\infty} e^{-k} (\cos(k\Delta_{y} g )-\cos(k\Delta_{y} f )+\cos(k\bar\Delta_{y} f )-\cos(k\bar\Delta_{y} g) ) \ dk  \ dy \ dx
  \end{eqnarray*}
  Now, we may estimate $\mathcal{A}_{2,1,{1}}+\mathcal{A}_{2,1,{2}}$.  By integrating by parts and by using \eqref{astuce}  	 \eqref{som} and \eqref{ddiff} together with the mean value theorem and classical Besov embeddings, one finds
  \begin{eqnarray*}
\mathcal{A}_{2,1,{1}}+\mathcal{A}_{2,1,{2}} &\lesssim&  \Vert \mathcal{U} \Vert_{\dot H^{1}} \int \frac{\Vert g \Vert_{\dot H^{2}} }{\vert y \vert^{3/2}} \frac{\Vert \mathcal{U} \Vert_{L^{\infty}}}{\vert y \vert^{3/2}} \ dy +  \Vert \mathcal{U} \Vert_{\dot H^{1}} \int \frac{\Vert \nabla g \Vert_{L^{\infty}} }{\vert y \vert^{3/2}} \frac{\Vert \nabla\mathcal{U} \Vert_{L^{2}}}{\vert y \vert^{3/2}} \ dy  \\
&& + \ \Vert \mathcal{U} \Vert_{\dot H^{1}} \int \frac{\Vert  \nabla \delta_{y} g \Vert_{L^{4}} }{\vert y \vert^{5/4}} \frac{\Vert \nabla \delta_{y}f+ \nabla \delta_{y}g \Vert_{L^{4}} }{\vert y \vert^{5/4}}  \frac{\Vert \mathcal{U} \Vert_{L^{\infty}}}{\vert y \vert^{3/2}} \ dy  \\
&\lesssim& \Vert \mathcal{U} \Vert_{\dot H^{1}} \int \frac{\Vert g \Vert_{\dot H^{2}} }{\vert y \vert^{3/2}} \frac{\Vert \mathcal{U} \Vert_{L^{\infty}}}{\vert y \vert^{3/2}} \ dy +  \Vert \mathcal{U} \Vert_{\dot H^{1}} \int \frac{\Vert \nabla g \Vert_{L^{\infty}} }{\vert y \vert^{3/2}} \frac{\Vert \nabla\mathcal{U} \Vert_{L^{2}}}{\vert y \vert^{3/2}} \ dy  \\
&+&\Vert \mathcal{U} \Vert_{\dot H^{1}} \left(\int \frac{\Vert  \nabla \delta_y g \Vert^4_{L^{4}} }{\vert y \vert^{5}} dy\right)^{1/4} \left(\int \frac{\Vert \nabla \delta_y f+ \nabla \delta_y g \Vert^4_{L^{4}} }{\vert y \vert^{5}} dy\right)^{1/4}   \left(\int \frac{\Vert \mathcal{U} \Vert^{2}_{L^{\infty}}}{\vert y \vert^{3}} \ dy \right)^{1/2} \\
&\lesssim&   \Vert \mathcal{U} \Vert_{\dot H^{1}}\Vert \mathcal{U} \Vert_{\dot H^{3/2}} \left(
\Vert g \Vert_{\dot H^{5/2}} + \Vert g \Vert^2_{\dot H^{9/4}} + \Vert f  \Vert_{\dot H^{9/4}} \Vert g \Vert_{\dot H^{9/4}}  \right)\\
\end{eqnarray*}

  Using the fact that $\dot H^{9/4}=[\dot H^{2}, \dot H^{5/2}]_{1/2}$ one finally gets
  
   \begin{eqnarray*}
   \mathcal{A}_{2,1,{1}}+\mathcal{A}_{2,1,{2}} &\lesssim&  \Vert \mathcal{U} \Vert_{\dot H^{1}}\Vert \mathcal{U} \Vert_{\dot H^{3/2}} \left( \Vert g \Vert_{\dot H^{5/2}} +\Vert g \Vert_{\dot H^{5/2}} \Vert g \Vert_{\dot H^{2}} +\Vert g \Vert^{1/2}_{\dot H^{2}} \Vert g \Vert^{1/2}_{\dot H^{5/2}} \Vert f \Vert^{1/2}_{\dot H^{2}} \Vert f \Vert^{1/2}_{\dot H^{5/2}} \right)
\end{eqnarray*}

Then, we consider $\mathcal{A}_{2,1,{3}}+\mathcal{A}_{2,1,{4}}$, we have
\begin{eqnarray*}
\mathcal{A}_{2,1,3}&=& \frac{1}{8} \int \nabla \mathcal{U}  .\nabla \int ( \Delta_{y} \nabla g-\nabla \bar\Delta_{y} g) . \frac{y}{\vert y \vert^2}  \left(\cos(\arctan(\Delta_{y} f ))+\cos(\arctan(\bar\Delta_{y} f ))\right) \\
  && \int_{0}^{\infty} e^{-k}( \cos(k\bar\Delta_{y} f )+\cos(k\Delta_{y} f ) ) \ dk  \ dy \ dx  \\
  &=&\frac{1}{8} \int \nabla \mathcal{U}  .\nabla \int ( \Delta_{y} \nabla g-\nabla \bar\Delta_{y} g) . \frac{y}{\vert y \vert^2}  \\
  &&\left(\cos(\arctan(\Delta_{y} f ))-\cos(\arctan(\Delta_{y} g))-\cos(\arctan(\bar\Delta_{y} g))+\cos(\arctan(\bar\Delta_{y} f ))\right) \\
  && \int_{0}^{\infty} e^{-k}( \cos(k\bar\Delta_{y} f )+\cos(k\Delta_{y} f ) ) \ dk  \ dy \ dx  \\
  &+&\frac{1}{8} \int \nabla \mathcal{U}  .\nabla \int ( \Delta_{y} \nabla g-\nabla \bar\Delta_{y} g) . \frac{y}{\vert y \vert^2}  \left(\cos(\arctan(\Delta_{y} g))+\cos(\arctan(\bar\Delta_{y} g)))\right) \\
  && \int_{0}^{\infty} e^{-k}( \cos(k\bar\Delta_{y} f )+\cos(k\Delta_{y} f ) ) \ dk  \ dy \ dx  \\
  &=&\frac{1}{8} \int \nabla \mathcal{U}  .\nabla \int ( \Delta_{y} \nabla g-\nabla \bar\Delta_{y} g) . \frac{y}{\vert y \vert^2}  \\
  &&\left(\cos(\arctan(\Delta_{y} f ))-\cos(\arctan(\Delta_{y} g))-\cos(\arctan(\bar\Delta_{y} g))+\cos(\arctan(\bar\Delta_{y} f ))\right) \\
  && \int_{0}^{\infty} e^{-k}( \cos(k\bar\Delta_{y} f )+\cos(k\Delta_{y} f ) ) \ dk  \ dy \ dx  \\
  &+&\frac{1}{8} \int \nabla \mathcal{U}  .\nabla \int ( \Delta_{y} \nabla g-\nabla \bar\Delta_{y} g) . \frac{y}{\vert y \vert^2}  \left(\cos(\arctan(\Delta_{y} g))+\cos(\arctan(\bar\Delta_{y} g)))\right) \\
  && \int_{0}^{\infty} e^{-k}( \cos(k\bar\Delta_{y} f )- \cos(k\bar\Delta_{y} g )- \cos(k\Delta_{y} g )+\cos(k\Delta_{y} f ) ) \ dk  \ dy \ dx  \\
   &+&\frac{1}{8} \int \nabla \mathcal{U}  .\nabla \int ( \Delta_{y} \nabla g-\nabla \bar\Delta_{y} g) . \frac{y}{\vert y \vert^2}  \left(\cos(\arctan(\Delta_{y} g))+\cos(\arctan(\bar\Delta_{y} g)))\right) \\
  && \int_{0}^{\infty} e^{-k}(  \cos(k\bar\Delta_{y} g )+ \cos(k\Delta_{y} g ) ) \ dk  \ dy \ dx  
\end{eqnarray*}
The last term cancels out with $\mathcal{A}_{2,1,4}$ and therefore,
\begin{eqnarray*}
\mathcal{A}_{2,1,3}+\mathcal{A}_{2,1,4}&=&\frac{1}{8} \int \nabla \mathcal{U}  .\nabla \int ( \Delta_{y} \nabla g-\nabla \bar\Delta_{y} g) . \frac{y}{\vert y \vert^2}  \\
  &&\left(\cos(\arctan(\Delta_{y} f ))-\cos(\arctan(\Delta_{y} g))-\cos(\arctan(\bar\Delta_{y} g))+\cos(\arctan(\bar\Delta_{y} f ))\right) \\
  && \int_{0}^{\infty} e^{-k}( \cos(k\bar\Delta_{y} f )+\cos(k\Delta_{y} f ) ) \ dk  \ dy \ dx  \\
  &+&\frac{1}{8} \int \nabla \mathcal{U}  .\nabla \int ( \Delta_{y} \nabla g-\nabla \bar\Delta_{y} g) . \frac{y}{\vert y \vert^2}  \left(\cos(\arctan(\Delta_{y} g))+\cos(\arctan(\bar\Delta_{y} g)))\right) \\
  && \int_{0}^{\infty} e^{-k}( \cos(k\bar\Delta_{y} f )- \cos(k\bar\Delta_{y} g )- \cos(k\Delta_{y} g )+\cos(k\Delta_{y} f ) ) \ dk  \ dy \ dx  
  \end{eqnarray*}
  
  It has the same structure as the term $\mathcal{A}_{2,1,{1}}+\mathcal{A}_{2,1,{2}}$, by integrating by parts, following the same steps we easily get that 
  
  \begin{eqnarray*}
\mathcal{A}_{2,1,3}+\mathcal{A}_{2,1,4}&\lesssim&  \Vert \mathcal{U} \Vert_{\dot H^{1}}\Vert \mathcal{U} \Vert_{\dot H^{3/2}} \left( \Vert g \Vert_{\dot H^{5/2}} +\Vert g \Vert_{\dot H^{5/2}} \Vert g \Vert_{\dot H^{2}} +\Vert g \Vert^{1/2}_{\dot H^{2}} \Vert g \Vert^{1/2}_{\dot H^{5/2}} \Vert f \Vert^{1/2}_{\dot H^{2}} \Vert f \Vert^{1/2}_{\dot H^{5/2}} \right)
\end{eqnarray*}

  We now estimate $\mathcal{A}_{2,1,5}+\mathcal{A}_{2,1,6}$, we first notice that 
  
   \begin{eqnarray*}
\mathcal{A}_{2,1,5}&=&-\frac{1}{4} \int \nabla \mathcal{U}. \nabla  \int\bar\Delta_{y} \nabla g . \frac{y}{\vert y \vert^2} \left(\cos(\arctan(\bar\Delta_{y} f ))+\cos(\arctan(\Delta_{y} f ))\right) \times \\
  && \int_{0}^{\infty} e^{-k} (\cos(k\bar\Delta_{y} f )-\cos(k\Delta_{y} f ))  \ dk  \ dy \ dx \\
  &=&-\frac{1}{4} \int \nabla \mathcal{U}. \nabla  \int\bar\Delta_{y} \nabla g . \frac{y}{\vert y \vert^2} \\
  &&\left(\cos(\arctan(\bar\Delta_{y} f ))-\cos(\arctan(\bar\Delta_{y} g ))-\cos(\arctan(\Delta_{y} g ))+\cos(\arctan(\Delta_{y} f ))\right) \times \\
  && \int_{0}^{\infty} e^{-k} (\cos(k\bar\Delta_{y} f )-\cos(k\Delta_{y} f ))  \ dk  \ dy \ dx \\
  &-&\frac{1}{4} \int \nabla \mathcal{U}. \nabla  \int\bar\Delta_{y} \nabla g . \frac{y}{\vert y \vert^2} \left(\cos(\arctan(\bar\Delta_{y} g ))+\cos(\arctan(\Delta_{y} g ))\right) \times \\
  && \int_{0}^{\infty} e^{-k} (\cos(k\bar\Delta_{y} f )-\cos(k\Delta_{y} f ))  \ dk  \ dy \ dx \\
\end{eqnarray*}
On the other hand, we have

 \begin{eqnarray*}
\mathcal{A}_{2,1,6}&=&\frac{1}{4} \int \nabla \mathcal{U} . \nabla \int\bar\Delta_{y} \nabla g . \frac{y}{\vert y \vert^2} \left(\cos(\arctan(\bar\Delta_{y} g ))+\cos(\arctan(\Delta_{y} g ))\right) \times \\
  && \int_{0}^{\infty} e^{-k} (\cos(k\bar\Delta_{y} g )-\cos(k\Delta_{y} g ))  \ dk  \ dy\ dx \\
  &=&\frac{1}{4} \int \nabla \mathcal{U} . \nabla \int\bar\Delta_{y} \nabla g . \frac{y}{\vert y \vert^2} \left(\cos(\arctan(\bar\Delta_{y} g ))+\cos(\arctan(\Delta_{y} g ))\right) \times \\
  && \int_{0}^{\infty} e^{-k} (\cos(k\bar\Delta_{y} g )-\cos(k\bar\Delta_{y} f )+\cos(k\Delta_{y} f )-\cos(k\Delta_{y} g ))  \ dk  \ dy\ dx \\
  &+&\frac{1}{4} \int \nabla \mathcal{U} . \nabla \int\bar\Delta_{y} \nabla g . \frac{y}{\vert y \vert^2} \left(\cos(\arctan(\bar\Delta_{y} g ))+\cos(\arctan(\Delta_{y} g ))\right) \times \\
  && \int_{0}^{\infty} e^{-k} (\cos(k\bar\Delta_{y} f )-\cos(k\Delta_{y} f ))  \ dk  \ dy\ dx \\
  \end{eqnarray*}
One notices that the last two terms in $\mathcal{A}_{2,1,5}$ and $\mathcal{A}_{2,1,6}$ cancel out. Hence,
\begin{eqnarray*}
\mathcal{A}_{2,1,5}+\mathcal{A}_{2,1,6}&=&-\frac{1}{4} \int \nabla \mathcal{U}. \nabla  \int\bar\Delta_{y} \nabla g . \frac{y}{\vert y \vert^2} \\
  &&\left(\cos(\arctan(\bar\Delta_{y} f ))-\cos(\arctan(\bar\Delta_{y} g ))-\cos(\arctan(\Delta_{y} g ))+\cos(\arctan(\Delta_{y} f ))\right) \times \\
  && \int_{0}^{\infty} e^{-k} (\cos(k\bar\Delta_{y} f )-\cos(k\Delta_{y} f ))  \ dk  \ dy \ dx \\
  &+&\frac{1}{4} \int \nabla \mathcal{U} . \nabla \int\bar\Delta_{y} \nabla g . \frac{y}{\vert y \vert^2} \left(\cos(\arctan(\bar\Delta_{y} g ))+\cos(\arctan(\Delta_{y} g ))\right) \times \\
  && \int_{0}^{\infty} e^{-k} (\cos(k\bar\Delta_{y} g )-\cos(k\bar\Delta_{y} f )+\cos(k\Delta_{y} f )-\cos(k\Delta_{y} g ))  \ dk  \ dy\ dx 
  \end{eqnarray*}
  
  Again here, we integrate by parts and we notice that this term can be estimated in a similar manner as $\mathcal{A}_{2,1,{1}}+\mathcal{A}_{2,1,{2}}$ and therefore,
  \begin{eqnarray*}
\mathcal{A}_{2,1,5}+\mathcal{A}_{2,1,6}&\lesssim& \Vert \mathcal{U} \Vert_{\dot H^{1}}\Vert \mathcal{U} \Vert_{\dot H^{3/2}} \left( \Vert g \Vert_{\dot H^{5/2}} +\Vert g \Vert_{\dot H^{5/2}} \Vert g \Vert_{\dot H^{2}} +\Vert g \Vert^{1/2}_{\dot H^{2}} \Vert g \Vert^{1/2}_{\dot H^{5/2}} \Vert f \Vert^{1/2}_{\dot H^{2}} \Vert f \Vert^{1/2}_{\dot H^{5/2}} \right)
\end{eqnarray*}

  It remains to estimate $\mathcal{A}_{2,1,7}+\mathcal{A}_{2,1,8}$. To do so, one first writes that
  
  \begin{eqnarray*}
  \mathcal{A}_{2,1,7}&=&-\frac{1}{4}\int \nabla \mathcal{U} . \nabla \int\Delta_{y} \nabla g . \frac{y}{\vert y \vert^2} \\
  && (\cos(\arctan(\bar\Delta_{y} f ))-\cos(\arctan(\bar\Delta_{y} g ))+\cos(\arctan(\Delta_{y} g ))-\cos(\arctan(\Delta_{y} f )))\times\\
&& \int_{0}^{\infty} e^{-k} (\cos(k\Delta_{y} f) + \cos(k\bar\Delta_{y} f)) \ dk  \ dy \ dx \\
&+&\frac{1}{4}\int \nabla \mathcal{U} . \nabla \int\Delta_{y} \nabla g . \frac{y}{\vert y \vert^2} (\cos(\arctan(\bar\Delta_{y} g ))-\cos(\arctan(\Delta_{y} g )))\times\\
&& \int_{0}^{\infty} e^{-k} (\cos(k\Delta_{y} f) + \cos(k\bar\Delta_{y} f)) \ dk  \ dy \ dx
  \end{eqnarray*}

  Then, we may rewrite $\mathcal{A}_{2,1,8}$, as
  \begin{eqnarray*}
\mathcal{A}_{2,1,8}&=&\frac{1}{4}\int \nabla \mathcal{U} . \nabla\int \Delta_{y} \nabla g . \frac{y}{\vert y \vert^2} (\cos(\arctan(\bar\Delta_{y} g ))-\cos(\arctan(\Delta_{y} g )))\times\\
&&\int_{0}^{\infty} e^{-k} (\cos(k\Delta_{y} g)-\cos(k\Delta_{y} f)-\cos(k\Delta_{y} g)) + \cos(k\bar\Delta_{y} g)) \ dk  \ dy \ dx   \\
&-&\frac{1}{4}\int \nabla \mathcal{U} . \nabla\int \Delta_{y} \nabla g . \frac{y}{\vert y \vert^2} (\cos(\arctan(\bar\Delta_{y} g ))-\cos(\arctan(\Delta_{y} g )))\times\\
&&\int_{0}^{\infty} e^{-k} (\cos(k\Delta_{y} f)+\cos(k\Delta_{y} f))  \ dk  \ dy \ dx   \\
  \end{eqnarray*}
  Then, we again notice that the last terms cancel out and one finds that
  \begin{eqnarray*}
  \mathcal{A}_{2,1,7}+\mathcal{A}_{2,1,8}&=&-\frac{1}{4} \int \nabla \mathcal{U}. \nabla  \int\bar\Delta_{y} \nabla g . \frac{y}{\vert y \vert^2} \\
  &&\left(\cos(\arctan(\bar\Delta_{y} f ))-\cos(\arctan(\bar\Delta_{y} g ))-\cos(\arctan(\Delta_{y} g ))+\cos(\arctan(\Delta_{y} f ))\right) \\
  && \times \int_{0}^{\infty} e^{-k} (\cos(k\bar\Delta_{y} f )-\cos(k\Delta_{y} f ))  \ dk  \ dy \ dx  \\
  &+&\frac{1}{4}\int \nabla \mathcal{U} . \nabla\int \Delta_{y} \nabla g . \frac{y}{\vert y \vert^2} (\cos(\arctan(\bar\Delta_{y} g ))-\cos(\arctan(\Delta_{y} g )))\times\\
&&\int_{0}^{\infty} e^{-k} (\cos(k\Delta_{y} g)-\cos(k\Delta_{y} f)-\cos(k\Delta_{y} g)) + \cos(k\bar\Delta_{y} g)) \ dk  \ dy \ dx   
  \end{eqnarray*}
  
  By integrating by parts, we again notice that this term is similar to $\mathcal{A}_{2,1,{1}}+\mathcal{A}_{2,1,{2}}$, hence we have
  \begin{eqnarray*}
\mathcal{A}_{2,1,7}+\mathcal{A}_{2,1,8}&\lesssim& \Vert \mathcal{U} \Vert_{\dot H^{1}}\Vert \mathcal{U} \Vert_{\dot H^{3/2}} \left( \Vert g \Vert_{\dot H^{5/2}} +\Vert g \Vert_{\dot H^{5/2}} \Vert g \Vert_{\dot H^{2}} +\Vert g \Vert^{1/2}_{\dot H^{2}} \Vert g \Vert^{1/2}_{\dot H^{5/2}} \Vert f \Vert^{1/2}_{\dot H^{2}} \Vert f \Vert^{1/2}_{\dot H^{5/2}} \right)
\end{eqnarray*}
Finally,  
\begin{eqnarray*}
\mathcal{A}_{2} &\lesssim& C(K) \Vert \mathcal{U} \Vert^2_{\dot H^{1}} \Vert g \Vert^2_{\dot H^{5/2}}+ \frac{1}{100(1+K^2)^{3/2}} \Vert \mathcal{U} \Vert^2_{\dot H^{3/2}} \\
&+&C(K) \left(\sup_{t \in [0,T]} \Vert g\Vert_{\dot H^{2}}\right)^2 \Vert \mathcal{U} \Vert^2_{\dot H^{1}} \Vert g \Vert^2_{\dot H^{5/2}}+ \frac{1}{100(1+K^2)^{3/2}} \Vert \mathcal{U} \Vert^2_{\dot H^{3/2}} \\
&+&C(K) \left(\sup_{t \in [0,T]} \Vert g\Vert_{\dot H^{2}}\right) \left(\sup_{t \in [0,T]} \Vert f\Vert_{\dot H^{2}}\right)  \Vert \mathcal{U} \Vert^2_{\dot H^{1}} \Vert g \Vert_{\dot H^{5/2}}\Vert f \Vert_{\dot H^{5/2}}+ \frac{1}{100(1+K^2)^{3/2}} \Vert \mathcal{U} \Vert^2_{\dot H^{3/2}} \\
\end{eqnarray*}

 Combining the latter inequality with the estimate obtained for $\mathcal{A}_{1}$, one finally finds
\begin{eqnarray*}
\frac{1}{2} \partial_{t} \Vert \nabla\mathcal{U} \Vert^{2}_{L^{2}} + \frac{1}{2(1+K(t)^{2})^{3/2}} \Vert \mathcal{U}\Vert^{2}_{\dot H^{3/2}}  &\lesssim& C(K)\Vert \nabla \mathcal{U} \Vert^{2}_{L^{2}} \left(\left(\sup_{t \in [0,T]} \Vert f\Vert_{\dot H^{2}}\right)^2 +\left(\sup_{t \in [0,T] } \Vert g\Vert_{\dot H^{2}}\right)^2 \right) \\
&\times& \left(\Vert g \Vert^2_{H^{5/2}}+\Vert f \Vert^{2}_{H^{5/2}}\right).
\end{eqnarray*}

Hence, integrating in time $s \in [0,T]$ and using Gronwall's inequality one finally concludes that 

\begin{eqnarray*}
\sup_{t \in [0,T]} \Vert \mathcal{U} (t)\Vert^{2}_{\dot H^{1}}   &\leq& \Vert \mathcal{U}_{0} \Vert_{\dot H^1} \exp\left(C(K) \left(\left(\sup_{t \in [0,T]} \Vert f\Vert_{\dot H^{2}}\right)^2+\left(\sup_{t \in [0,T]} \Vert g\Vert_{\dot H^{2}}\right)^2  \right) \right.\\
 && \times \left. \int_{0}^{T} \left(\Vert g (s) \Vert^2_{H^{5/2}}+\Vert f(s) \Vert^{2}_{H^{5/2}}\right) \ ds \right).
\end{eqnarray*}
Which readily gives uniqueness.
\qed

}}
\section*{Acknowledgments}
\noindent  Both  F.G. and O.L. were supported by the ERC through the Starting Grant project H2020-EU.1.1.-639227. F.G. were partially supported by the grant MTM2017-89976-P (Spain).

\vspace{3mm}

\noindent\textbf{Francisco Gancedo,} Departamento de An\'alisis 
Matem\'atico \& IMUS Universidad de Sevilla , Spain. 

\vspace{3mm}

\noindent\textbf{Omar Lazar,} Departamento de An\'alisis Matem\'atico \& IMUS, Universidad de Sevilla, Spain.

\end{document}